\documentclass[a4paper, 10pt, abstracton]{scrartcl}

\usepackage{amsmath}
\usepackage{amssymb}
\usepackage{graphicx}
\usepackage{epstopdf}
\usepackage{inputenc}
\usepackage{geometry}
\usepackage{a4wide}

\usepackage[font=small]{caption}
\usepackage{subcaption}
\usepackage{cite}
\usepackage{calc}
\usepackage{float}
\usepackage{grffile} 
\usepackage{enumitem}

\usepackage{color} 
\usepackage{xcolor}

\usepackage{mathtools}

\usepackage{comment} 

\newcommand{\U}{\mathcal{U}}

\newcommand{\R}{\mathbb{R}}
\newcommand{\X}{\mathcal{X}}

\newcommand{\Np}{N}
\newcommand{\pb}{\pmb{x}_*}
\newcommand{\gb}{\hat{\pmb{x}}_*}
\newcommand{\budget}{B^{\text{max}}}
\newcommand{\inner}{B^{\text{in}}}

\usepackage{algpseudocode}
\usepackage{algorithm}

\usepackage{multirow}
\usepackage{rotating}

\usepackage{authblk}

\usepackage{lscape} 

\usepackage[normalem]{ulem} 

\algnewcommand{\IIf}[1]{\State\algorithmicif\ #1\ \algorithmicthen} 
\algnewcommand{\EndIIf}{\unskip\ \algorithmicend\ \algorithmicif} 

\usepackage{mathrsfs} 
\usepackage{eufrak} 

\definecolor{Gray}{gray}{0.9} 

\begin{document}

\title{Particle Swarm Metaheuristics for Robust Optimisation with Implementation Uncertainty\thanks{Partially funded through EPSRC grants EP/L504804/1 and EP/M506369/1.}}
\author[1]{Martin Hughes\footnote{Corresponding author. Email: \texttt{m.hughes6@lancaster.ac.uk}} }
\author[2]{Marc Goerigk}
\author[1]{Trivikram Dokka}
\affil[1]{Department of Management Science, Lancaster University, United Kingdom}
\affil[2]{Network and Data Science Management, University of Siegen, Germany}

\date{} 

\maketitle


\abstract{We consider global non-convex optimisation problems under uncertainty. In this setting, it is not possible to implement a desired solution exactly. Instead, any other solution within some distance to the intended solution may be implemented. The aim is to find a robust solution, i.e., one where the worst possible solution nearby still performs as well as possible.

Problems of this type exhibit another maximisation layer to find the worst case solution within the minimisation level of finding a robust solution, which makes them harder to solve than classic global optimisation problems. So far, only few methods have been provided that can be applied to black-box problems with implementation uncertainty. We improve upon existing techniques by introducing a novel particle swarm based framework which adapts elements of previous approaches, combining them with new features in order to generate more effective techniques. In computational experiments, we find that our new method outperforms state of the art comparator heuristics in almost 80\% of cases.
}

\textbf{Keywords: } robust optimisation; implementation uncertainty; metaheuristics; global optimisation; particle swarm optimisation


\section{Introduction}
\label{sec:introduction}

Decision making in the face of uncertainty is a widespread challenge. In many real-world situations it is common practice to use models to support informed decision making. However model run times and the extent of the decision variable solution space may render an extensive assessment of the problem space computationally impractical. In such circumstances an efficient global optimisation search approach is needed. The consideration of uncertainty in the modelling process, reflecting uncertainty in the real-world problem, may impact model outputs and therefore the optimum objective function value. Thus uncertainty adds an additional feature into any global optimisation search. Whilst simply ignoring the uncertainty is one strategy, such an approach has been shown to produce sub-optimal results, see \cite{BenTalElGhaouiNemirovski2009, GoerigkSchobel2016}. An approach is required that can identify a solution that performs well over a range of scenarios as opposed to simply in the expected case.

The use of a model in the form of a mathematical program is preferable from an optimisation standpoint, as such models may be solved efficiently with the determination of global optima guaranteed. However such an approach necessitates that the problem at hand can be adequately expressed in the form of a mathematical program. In many real-world situations this is not possible. Rather some more general form of simulation model will be used, which from an optimisation perspective, may be considered a black-box: decision variables values are input and an objective value is output. For such black-box problems a more general search technique that can be applied to any model is required, such as a metaheuristic. 

Optimisation under uncertainty is typically approached using either stochastic or robust techniques. In stochastic optimisation the probability distributions of the uncertain parameters are assumed to be known and the fitness of any solution is determined by some statistical measure, see \cite{PaenkeBrankeJin2006, HomemdeMelloBayraksan2014}. By contrast robust optimisation only assumes that all uncertainty scenarios can be described by some set \cite{BenTalNemirovski1998}. A classic robust approach finds a solution that optimises its performance in the worst case. This is known as min max: at any point in the decision variable space an inner objective is employed to identify the maximal function value in the point's uncertainty neighbourhood, with an outer objective employed to identify the minimum maximal value.

Here we develop new metaheuristics for the robust global optimisation of black-box problems, including non-convex problems. These algorithms accommodate implementation uncertainty, where a desired solution may be somewhat perturbed in a real-world setting. We adopt the classic robust worst case approach. More formally, the general optimisation problem to be considered here is:
\begin{align*}
	\quad \min\ & f(\pmb{x}) \\
	\text{s.t. } & \pmb{x} \in \X
\end{align*}
where $f: \R^n \to \R$ is the objective function, $\pmb{x}=(x_{1}, x_{2}, \ldots, x_{n})^T$ is the $n$-dimensional vector of decision variables, and $\X \subseteq \R^n$ is the set of feasible solutions. We use the notation $[n]:=\{1,\ldots,n\}$ and assume box constraints $\X = \prod_{i\in[n]} [l_i,u_i]$. A penalty in the objective is assumed in the case of other feasibility constraints. Such a problem, without any consideration of uncertainty, is designated the nominal problem here.

\begin{figure}[htbp]
\centering
\begin{subfigure}{\textwidth}
\centering
\includegraphics[width=0.6\textwidth]{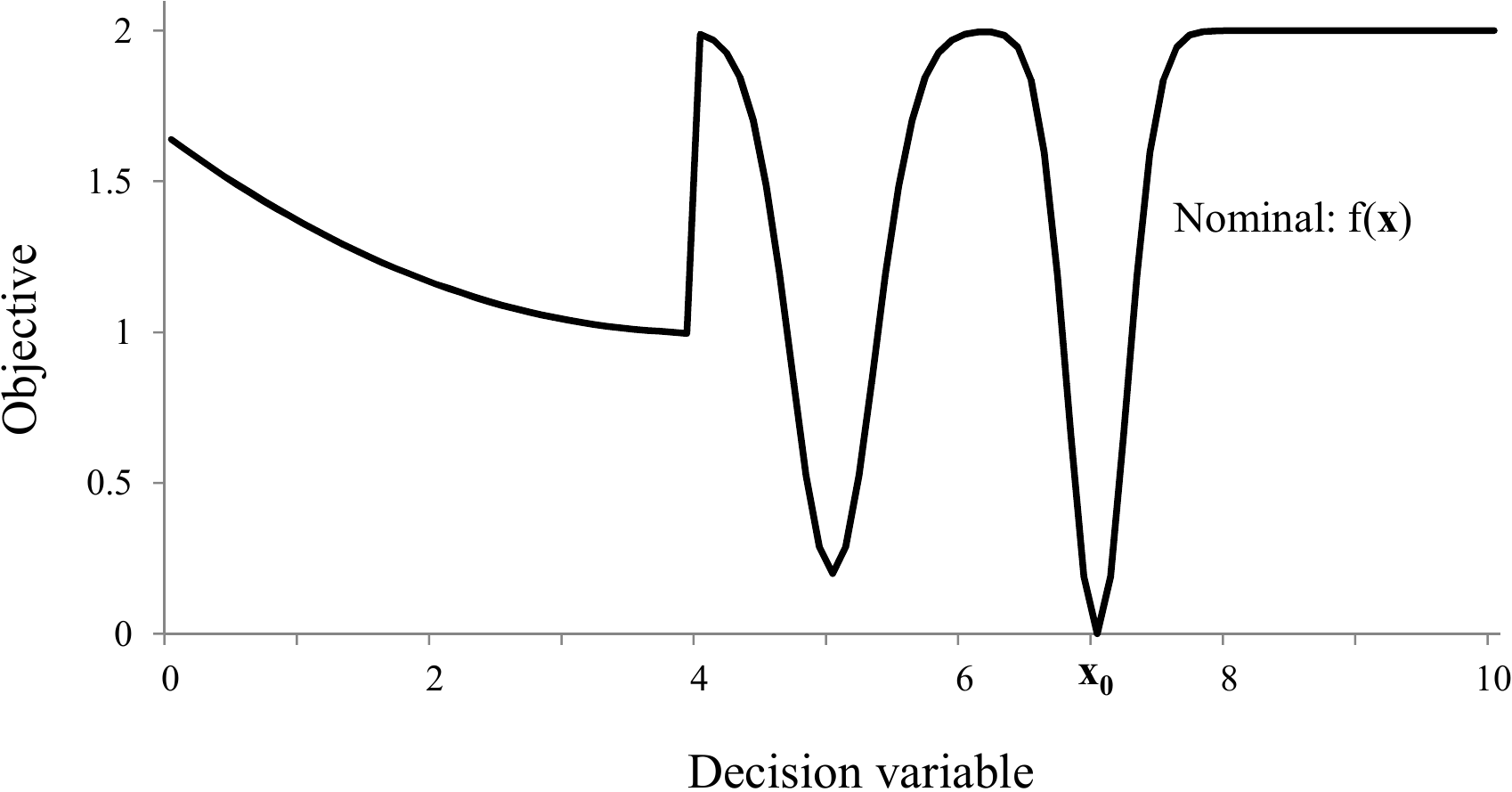} 
\caption{The nominal global optimum is at $\pmb{x}_0$.}
\label{fig:KruisselbrinkNominal}
\end{subfigure}%

\vspace{5mm}

\begin{subfigure}{\textwidth}
\centering
\includegraphics[width=0.6\textwidth]{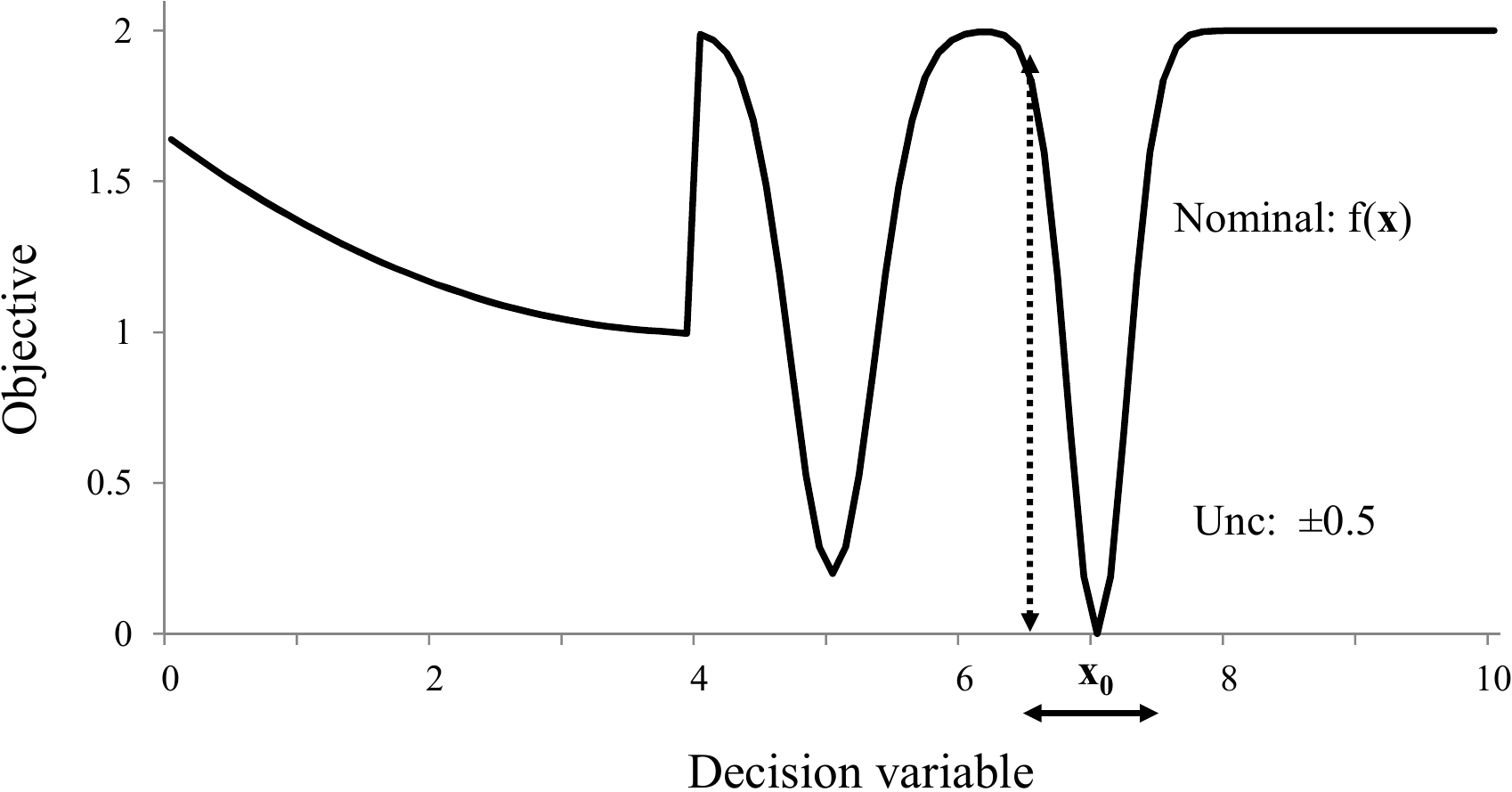} 
\caption{With uncertainty only a `close' solution may be realised, which is of particular concern in the region of $\pmb{x}_0$.}
\label{fig:KruisselbrinkWithUnc}
\end{subfigure}%

\vspace{5mm}

\begin{subfigure}{\textwidth}
\centering
\includegraphics[width=0.6\textwidth]{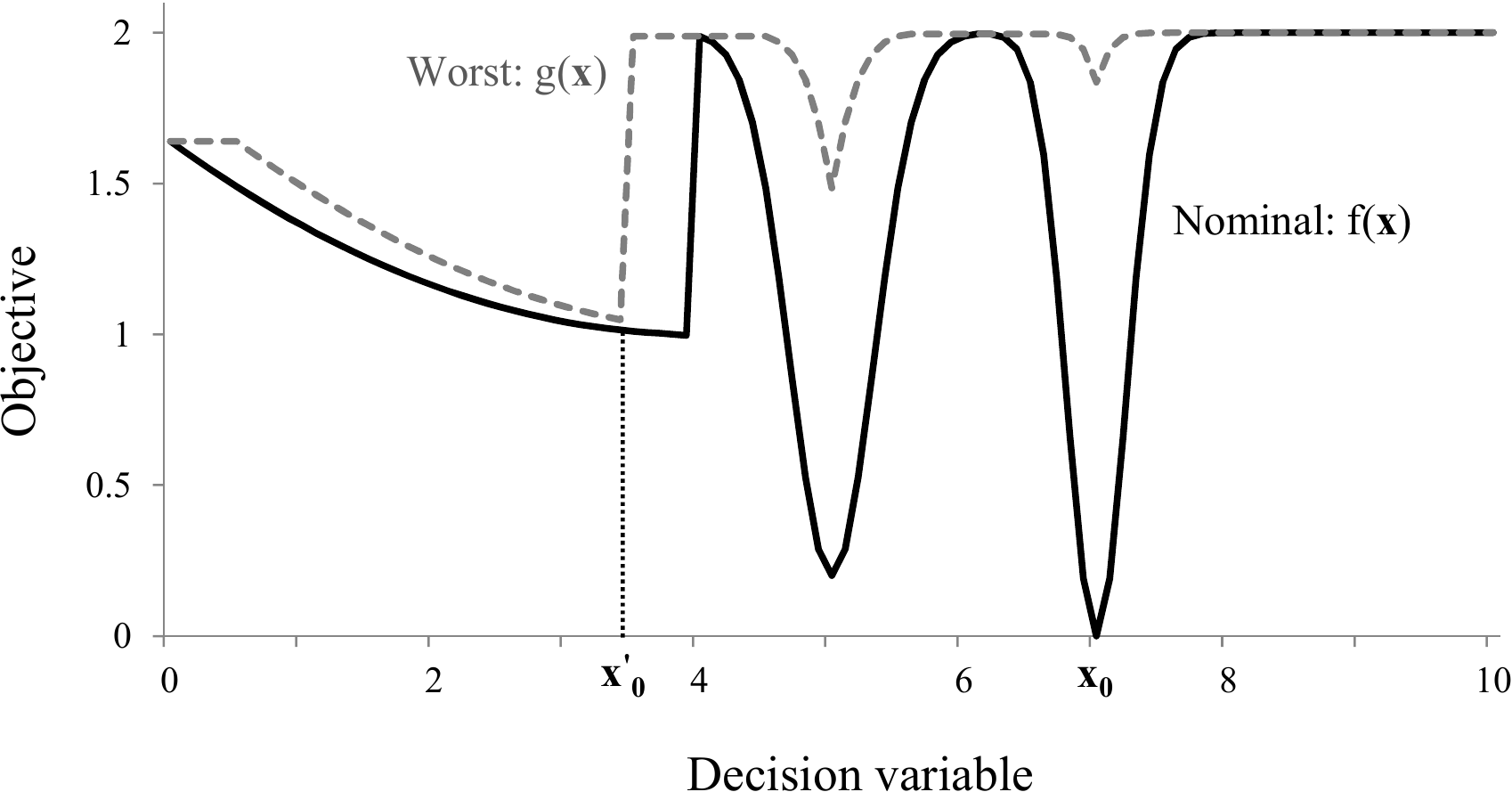} 
\caption{The worst case cost curve (dashed grey) is generated by determining the maximum objective value in the uncertainty neighbourhood around all points $\pmb{x}$ on the nominal (solid black) curve. Due to the uncertainty the global optimum shifts to $\pmb{x}'_0$.}
\label{fig:KruisselbrinkWorst}
\end{subfigure}%

\caption{One dimensional problem due to \cite{Kruisselbrink2012}.}
\end{figure}

As an example, consider a non-convex one dimensional problem due to \cite{Kruisselbrink2012}. For the nominal problem, shown in Figure~\ref{fig:KruisselbrinkNominal}, some standard metaheuristic could be used to locate the global minimum at $\pmb{x}_0$. However if the solution that a decision maker wants to implement is somewhat perturbed in practice, the potential impact on the identification of the global minimum needs to be taken into consideration. The sensitivity of the objective to variations in the region of $\pmb{x}_0$ is of particular concern, as highlighted in Figure~\ref{fig:KruisselbrinkWithUnc}.

We assume only a solution similar to the desired solution $\pmb{x}$, $\tilde{\pmb{x}}=\pmb{x}+\Delta\pmb{x}$ may be achieved. The classic robust approach is then to find a robust solution $\pmb{x}$ such that for any $\tilde{\pmb{x}}$ from the uncertainty neighbourhood of $\pmb{x}$, the worst performance is optimised. As common in the literature (see, e.g., \cite{BertsimasNohadaniTeo2010}), we consider the following uncertainty neighbourhood (also known as uncertainty set):
\[
	\U:=\{ \Delta \pmb{x}\in\R^n \mid \| \Delta \pmb{x} \| \leq \Gamma \}
\]
where $\Gamma > 0$ defines the radius of the uncertainty neighbourhood around a solution $\pmb{x}\in\X$ and $\|\cdot\|$ represents the Euclidean norm. The worst case costs of a solution $\pmb{x}$ are: 
\[
	g(\pmb{x}):=\max_{\Delta \pmb{x} \in \U} f(\pmb{x} + \Delta \pmb{x})
\]
The min max robust optimisation problem is then:
\[
	\min_{\pmb{x}\in\X} g(\pmb{x}) = \min_{\pmb{x}\in\X} \max_{\Delta\pmb{x} \in \U} f(\pmb{x} + \Delta \pmb{x}) \tag{MM}
\]
Finding the robust global optimum is based on an outer minimum worst case cost objective function value in $\X$, such that that minimum objective is based on an inner maximisation of the nominal objective function in the uncertainty neighbourhood around each solution $\pmb{x}\in\X$.

Note that we assume that $f$ is not restricted to $\X$, as $\pmb{x}+\Delta\pmb{x}$ may be outside of $\X$. Alternatively $\pmb{x}+\Delta\pmb{x}\in\X$ for all $\Delta\pmb{x}\in\U$ could be achieved through a reduction in the size of the original $\X$ by $\Gamma$.

In Figure~\ref{fig:KruisselbrinkWorst} the worst case cost $g(\pmb{x})$ (dashed grey curve) at any individual point $\pmb{x}$ can be determined by assessing the uncertainty neighbourhood around that point, in order to identify the maximum value within that uncertainty neighbourhood. Then within the global minimisation search the nominal objective is superseded by the worst case cost. It can be observed that the global optimum has shifted to $\pmb{x}'_{0}$. It should also be noted that if we were to ignore the implementation uncertainty and simply accept $\pmb{x}_{0}$ as the global optimum, which is a common approach in practice, then we risk the possibility of a very poor outcome, i.e., $g(\pmb{x}'_0) < g(\pmb{x}_0)$, whereas $f(\pmb{x}'_0) > f(\pmb{x}_0)$.

\paragraph{Contributions and outline.} 

In this paper we propose a particle swarm optimisation (PSO) framework encompassing new robust metaheuristics for black-box problems under implementation uncertainty. These techniques can be applied to general problems of reasonable dimension, where little if anything is known about the nature of the objective function surface, and under restrictions on the number of function evaluations (the budget) due to run-time considerations. The methods developed here are based on extensions of PSO, see \cite{KennedyEberhart1995, KennedyEberhartShi2001, Ghazali2009}. Specifically we employ a PSO frame, augmenting it with adapted elements of the robust local search descent directions (d.d.) approach due to \cite{BertsimasNohadaniTeo2007, BertsimasNohadaniTeo2010, BertsimasNohadaniTeo2010nonconvex}, and the robust global largest empty hypersphere (LEH) approach due to \cite{HughesGoerigkWright2019}, and introducing original features in order to generate novel techniques. We undertake a series of computational experiments comparing these new methods with a baseline robust PSO (rPSO), a global version of d.d.\ and LEH, see \cite{BertsimasNohadaniTeo2010, HughesGoerigkWright2019}. We find that our new metaheuristics considerably outperform these approaches on a large number of problem instances.

In Section~\ref{sec:sota} we review the current state of the art by discussing the relevant literature in Section~\ref{sec:literature} and details of the d.d.\ and LEH algorithms in Sections~\ref{sec:localRobustDD} and~\ref{sec:globalRobustLEH}. We discuss PSO in Section~\ref{sec:PSO}, including an extension of the nominal formulation to a baseline rPSO formulation. In Section~\ref{sec:baselineComparators} we provide an illustrative example. We then outline our new robust framework in Section~\ref{sec:EnhancedPSO}, including descriptions of our heuristic algorithms. In Section~\ref{sec:experimentsResults} we describe the experimental set up used to test our new heuristics, and present our results. We end with a summary and consideration of potential further work in Section~\ref{sec:concusionsFurtherWork}. The appendices cover descriptions of our experimental test functions, box plots of results, and a list of abbreviations.


\section{State of the art}
\label{sec:sota}

\subsection{Literature review}
\label{sec:literature}

The modern form of robust optimisation was first developed in \cite{KouvelisYu1997} and \cite{BenTalNemirovski1998}, since when the field has been strongly associated with mathematical programming through the use of appropriate uncertainty sets and the reformulation of mathematical programming problems, see e.g. the surveys \cite{bertsimas2011theory,gabrel2014recent,GoerigkSchobel2016}. In general, two types of uncertainty can be distinguished: parameter uncertainty, where the problem data is unknown; and implementation uncertainty, where the decision is subject to change during its implementation. In a mathematical programming context robust optimisation with parameter uncertainty has been widely applied to specific problems and formulations, while robust optimisation of general black-box problems with implementation uncertainty is much less widely addressed, see \cite{MarzatWalterPietLahanier2013, GoerigkSchobel2016, MarzatWalterPietLahanier2016}.

A worst case analysis can be approached by applying standard global metaheuristics to both the inner maximisation and the outer min. In a co-evolutionary approach inner and outer populations evolve separately, but the fitness of individuals in the outer min is determined by individuals in an inner maximisation, see \cite{Herrmann1999, ShiKrohling2002, Jensen2004, CramerSudhoffZivi2009, MasudaKuriharaAiyoshi2011}. However a completely brute force approach using full inner maximisation searches to inform the outer minimisation involves large numbers of function evaluations (model runs), see \cite{MarzatWalterPietLahanier2016}. Additional simplifications or assumptions are required to reduce the number of function evaluations in a co-evolutionary approach, see \cite{CramerSudhoffZivi2009, MasudaKuriharaAiyoshi2011}.

An alternative robust evolutionary approach, introduced by \cite{TsutsuiGhosh1997}, is based around the idea of `genetic algorithms with a robust solution searching scheme' (GAs/RS$^3$) \cite{BeyerSendhoff2007}. Uncertainty is added to the individuals in the population prior to the determination of the next generation; the next generation is then determined based on an assessment of the fitness of the extended (uncertain) population. \cite{OngNairLum2006} adopts such an approach in a min max robust design analysis. This work also employs an approach that can be considered more generally in robust analyses, the use of emulation (surrogates or meta-models) alongside true objective function evaluations to reduce the potential burden of computational run times and the number of model-function evaluations, see \cite{BeyerSendhoff2007, VuDAmbrosioHamadiLiberti2016}. \cite{OngNairLum2006} use surrogates for the inner optimisation local search. In \cite{ZhouZhang2010} the inner maximisation is tackled using a surrogate-assisted evolutionary algorithm, whilst \cite{MarzatWalterPietLahanier2013, urRehmanLangelaarvanKeulen2014, MarzatWalterPietLahanier2016, urRehmanLangelaar2017} all employ Kriging meta-modelling techniques. By contrast \cite{ChenLucierSingerSyrgkanis2017, SandersEversonFieldsendRahat2019} employ Bayesian emulation approaches. \cite{ChenLucierSingerSyrgkanis2017} uses a mathematical programming approach assuming the availability of a valid Bayesian oracle, whilst \cite{SandersEversonFieldsendRahat2019} employ a Bayesian approach for very expensive-to-evaluate functions, applying it to test problems of up to 10 dimensions using only small numbers of function evaluations. However current emulation based approaches suffer from the same limitation, in that they struggle when applied to problems other than those of relatively low dimension.

Of particular interest here are the single-solution descent directions \cite{BertsimasNohadaniTeo2007, BertsimasNohadaniTeo2010nonconvex, BertsimasNohadaniTeo2010} and largest empty hypersphere \cite{HughesGoerigkWright2019} robust (min max) metaheuristics, and standard (i.e. not robust) population based metaheuristic, particle swarm optimisation \cite{KennedyEberhart1995, KennedyEberhartShi2001, Ghazali2009}.

The d.d.\ approach is actually a robust local search. Given a start point in the decision variable space an inner maximisation is undertaken in the point's uncertainty neighbourhood. From this neighbourhood search undesirable 'high cost points' (hcps) are identified, and a direction which optimally points away from all of these hcps is determined by solving a quadratic program. A step is taken in this descent direction, to a new point where the process is repeated until no such direction can be found. This approach is considered in more detail in Section~\ref{sec:localRobustDD}. This work also informs the global approach outlined in \cite{BertsimasNohadani2010} where similar techniques are applied to the inner maximisation, but the outer minimization is tackled by simulated annealing. In \cite{HughesGoerigkWright2019} d.d.\ has been simplistically extended to a global search through the use of random re-starts.

The LEH metaheuristic is a relatively new robust global approach which extends the idea of locally moving away from undesirable hcps to a global setting by identifying regions of the feasible region devoid of hcps and moving to the centres of such regions. Hence this approach is exploration-focussed. This approach is considered in more detail in Section~\ref{sec:globalRobustLEH}.

The PSO approach is a population-based metaheuristic inspired by swarm intelligence; the description given here is based on \cite{KennedyEberhart1995, KennedyEberhartShi2001, Ghazali2009}. A swarm consists of multiple particles, moving through the solution space. The position of each particle represents a point visited in the decision variable space, and the objective function value at that location. An additional attribute of a particle is its velocity, which here represents the vector (direction and step size) of the particle's movement. Based on some combination of an individual particle's own information and the collective information from other particles in the swarm, each particle moves to new locations as the algorithms iterates. It is the intention that what emerges from such complex, self-organising systems of particles approximates an efficient search of the solution space to identify global optima.

In terms of a robust PSO approach, \cite{ShiKrohling2002, MasudaKuriharaAiyoshi2011} consider two-swarm co-evolutionary PSO techniques while \cite{HughesGoerigkWright2019} employs a baseline rPSO approach as a comparator test heuristic. \cite{Dippel2010} develops several PSO formulations, and includes material on topologies, memory (archive) and test functions, however the approaches considered are essentially stochastic and only very low dimensional problems are used in the testing. PSO and a baseline rPSO are considered in more detail in Section~\ref{sec:PSO}.

\subsection{Local robust search using descent directions}
\label{sec:localRobustDD}

Descent directions \cite{BertsimasNohadaniTeo2007, BertsimasNohadaniTeo2010, BertsimasNohadaniTeo2010nonconvex} is a local search technique for solving the robust optimization problem (MM), which uses the points evaluated in each inner maximisation local uncertainty neighbourhood analysis to inform a gradient descent approach. At a given point $\pmb{x}$ an inner maximisation search is performed to approximate the worst case cost $\tilde{g}(\pmb{x})\approx g(\pmb{x})$. In \cite{BertsimasNohadaniTeo2010} an extensive two-stage gradient ascent search is employed for inner maximisations. All function evaluations are recorded in a history set $H$. From within the uncertainty neighbourhood $N(\pmb{x})=\{\pmb{x}+\Delta\pmb{x} \mid \Delta\pmb{x}\in\U\}$ around a candidate point $\pmb{x}$, the points with the greatest objective function values are identified as high cost points. The high cost set $H_\sigma(\pmb{x})$ at any given point $\pmb{x}$ is defined as:
\[
	H_\sigma(\pmb{x}) := \{ \pmb{x}' \in H\cap N(\pmb{x}) \mid f(\pmb{x}')  \geq  \tilde{g}(\pmb{x}) -  \sigma \}
\]
where $\sigma$ is a threshold value for identifying hcps. 

The intention is then to identify the descent direction $\pmb{d}$ projecting from candidate point $\pmb{x}^k$ at iteration $k$, which optimally points away from the points in $H_\sigma(\pmb{x}^k)$. This is achieved by maximising the angle $\theta$ between $\pmb{d}$ and the vectors connecting the points in $H_\sigma(\pmb{x}^k)$ to $\pmb{x}^k$. This is a second order cone problem and can be tackled using mathematical programming:
\begin{align*}
\min_{\pmb{d},\beta} \ & \beta \tag{Soc1} \\
\text{s.t. } &\| \pmb{d}\| \le 1 \tag{Soc2} \\
& \pmb{d}^T \Bigg( \frac{\pmb{h} - \pmb{x}^k}{\| \pmb{h} - \pmb{x}^k  \|} \Bigg) \le \beta & \forall \pmb{h} \in H_\sigma(\pmb{x}^k) \tag{Soc3} \\
& \beta \le -\varepsilon \tag{Soc4}
\end{align*}
Setting $\varepsilon$ as a small positive scalar makes $\beta$ negative in (Soc4). The left hand side of (Soc3) equates to $\| \pmb{d}\| \cos \theta$, and is calculated for all points in $H_\sigma(\pmb{x}^k)$. (Soc3) therefore states that $\beta$ will correspond to the maximum value for $\cos \theta$ across all hcps. The objective (Soc1) is to minimise $\beta$. As $\beta$ is negative the angle $\theta$ will be greater than $90^o$ and maximised. Finally, minimising $\beta$ in combination with (Soc2) normalises $\pmb{d}$.

The final component of the algorithm is the calculation of the step size to be taken once a descent direction $\pmb{d}$ has been determined. At iteration $k$ in the local search a step size $\rho^k$ just large enough to ensure that all of the points in $H_\sigma(\pmb{x})$ are at least on the boundary of the $\Gamma$-uncertainty neighbourhood of the next candidate solution at step $k+1$ is used. We set $ \pmb{x}^{k+1} = \pmb{x}^k + \rho^k \pmb{d} $, where $\rho^k$ can be calculated using:
\[
	\rho^k = \min \left\{ \pmb{d}^T(\pmb{h} - \pmb{x}^k) + \sqrt{(\pmb{d}^T(\pmb{h} - \pmb{x}^k))^2 - {\| \pmb{h} - \pmb{x}^k  \|}^2 + \Gamma^2} \mid \pmb{h} \in H_\sigma(\pmb{x}^k) \right\}  \tag{Rho} 
\]
In the original formulation from \cite{BertsimasNohadaniTeo2010} several loops are potentially applied in the algorithm, in order to try to identify a valid direction and to ensure that that direction is reasonable. The parameter $\sigma$ is incrementally changed up to some limit, if (Soc) cannot be solved initially. Also given a valid direction vector $\pmb{d}$, a further check is used to ensure that the step to be taken does not immediately encounter additional hcps from $H$ beyond $N(\pmb{x}^k)$. 

This local stepping continues until no descent direction can be identified, and it is assumed that a robust local minimum has been reached. The approach can be extended to approximate a global search by randomly re-starting a new search each time the previous one completes. In \cite{HughesGoerigkWright2019} this is employed within the constraint of a fixed budget of function evaluations.

\subsection{Global robust search using largest empty hyperspheres}
\label{sec:globalRobustLEH}

The largest empty hypersphere metaheuristic \cite{HughesGoerigkWright2019} is a global method where the search progresses by moving to locations in the feasible region that are furthest away from all `bad' points previously visited. Using the d.d.\ idea of identifying high cost points in a global sense, LEH uses a history set $H$ of points evaluated and a high cost set $H_\tau$ which is a subset of $H$ containing all points with nominal objective function value $f(\pmb{x})$ greater than a threshold $\tau$, which is set to the current estimated robust global minimum value. Note that with this notation, $H_{\sigma}(\pmb{x}) = H_{\tilde{g}(\pmb{x})-\sigma}\cap N(\pmb{x})$.

Given $H_\tau$, LEH uses a genetic algorithm (GA) to estimate a point $\pmb{x}^k\in\X$ which is furthest from all hcps in $H_\tau$. The search then moves to this point. This is repeated until the budget of available function evaluations is exhausted or no point $\pmb{x}^k\in\X$ can be identified which is at least $\Gamma$ away from all hcps. In either case the current estimate for the global robust minimum is accepted.

At each candidate point $\pmb{x}^k$ an inner maximisation analysis is undertaken, however beyond the initial (random) start point each point in an inner maximisation analysis is compared to $\tau$ such that the $\Gamma$-radius uncertainty neighbourhood search can be stopped prematurely if any objective function value $f(\pmb{x}^k + \Delta\pmb{x}^k)$ is greater than $\tau$. This is a recognition of the fact that the current point $\pmb{x}^k$ cannot improve on the current estimate of the robust global optimum. This stopping condition potentially enables the LEH approach to explore $\X$ more efficiently.


\section{Particle swarm optimisation}
\label{sec:PSO}

\subsection{Motivation}

The overarching motivation for our work is the development of improved robust metaheuristics for black-box problems under implementation uncertainty. Of particular interest are approaches that can be applied to general problems of moderate dimension, where run-time issues limit the numbers of function evaluations or model runs that can be undertaken. Whilst some work has been undertaken in this area, compared to optimisation without any consideration of uncertainty, stochastic optimisation, or robust optimisation for mathematical programming problems, this remains a less developed field.

Of the techniques currently applicable in this setting, the local d.d.\ \cite{BertsimasNohadaniTeo2007, BertsimasNohadaniTeo2010, BertsimasNohadaniTeo2010nonconvex} and global LEH \cite{HughesGoerigkWright2019} approaches offer considerable insight into some key issues. This is not least because whilst d.d.\ is locally exploitation focussed, LEH is exploration focussed in its targeting of regions of the decision variable space devoid of `poor' points. That is these techniques span the exploitation versus exploration divide.

In particular these two approaches offer alternative stances on how to contend with what might be considered the additional 'burden' of a robust analysis under limitations on the number of function evaluations, i.e. the need to expend evaluations in the uncertainty neighbourhood analysis around individual points in the decision variable space in order to determine the robust value (the inner maximisation). Under budgetary restrictions we must therefore add the balancing of better estimating a candidate point's robust value versus the extent of the outer minimisation search, into the mix of exploration versus exploitation. This trade-off is complex, see for example \cite{MirjaliliLewisMostaghim2015, DiazHandlXu2017}.

The d.d.\ approach explicitly uses the additional information gained from an uncertainty neighbourhood inner maximisation search to direct a local search, by identifying the direction that optimally points away from the worst neighbourhood points. By contrast a key component of the LEH approach is what is termed a 'stopping condition', that is the ability to terminate an inner maximisation search prior to completion due to a recognition that the current global robust best cannot be improved upon at a given candidate point. This has the potential to introduce significant efficiencies when expending function evaluations, thereby enabling a more extensive outer minimisation search. In fact we recognise the contrast in exploitation (d.d.) versus exploration (LEH) as an echo of the contrast between enhancing a search through the exploitation of the inner maximisation information (d.d.) versus attempting to limit inner maximisation searches in order to expend function evaluations more efficiently (LEH).

Here we are interested in the potential benefits of both the better use of the information gained from previous function evaluations, and of efficiency savings in terms of numbers of function evaluations. In particular we are interested in addressing these elements within a single framework. Compared to the individual-based d.d.\ and LEH techniques, we consider a population-based approach more able to encompass these features under a single structure. In order to identify a suitable population-based framework here, consideration must be given to the features necessary to enable both the use of the stopping condition component of LEH, and the use of uncertainty neighbourhood directional information generated by the calculation of some form of descent direct at a given candidate point.

The inclusion of a stopping condition requires that an inner maximisation search can be terminated early.  Consider for example, how this might work in a fitness-based approach such as a genetic algorithm (GA). In order to effectively determine robust fitness at a candidate point, an inner maximisation must be complete. Early termination would not generate adequate fitness information: given a stopping threshold multiple members of the GA population could terminate their inner maximisation when an uncertainty neighbourhood point has been identified that exceeds that threshold, potentially leading to each being designated a similar fitness level close to the threshold objective function value. However by contrast if each inner maximisation were to complete, individual fitnesses could vary substantially. These two cases could lead to substantially different next generations due to the discrepancies in estimated fitnesses for members of the population. Therefore any such fitness-based approach does not suit a stopping condition of the kind under consideration here.

However considering a swarm-based approach such as PSO, each particle already has an in-built feature that can be exploited for stopping purposes, the best historic robust objective function value for each particle $j$. In the PSO case, stopping an inner maximisation prematurely if any $\Gamma$-radius uncertainty neighbourhood function evaluation exceeds that personal best threshold, $\tau^j$, has no negative impact as movement in a standard PSO formulation is based on some combination of personal and neighbourhood best information. Neither of these pieces of information are affected by particle-level stopping. For a particle $j$, a function evaluation exceeding $\tau^j$ establishes that neither the historic best particle level information nor the current neighbourhood (e.g. global) best can be improved upon by the particle's current location. In such a situation terminating an inner search and moving on is appropriate, and desirable.

Furthermore, the primary component of the d.d.\ approach is the determination of a direction vector, therefore an approach that already uses a vector-based approach is best placed to accommodate further vector information. Individual particle level movement in a swarm-based approach such as PSO is vector-based. In addition, in considering how a d.d.\ vector and particle level stopping features might be incorporated into a PSO formulation, it can be recognised that both features might be incorporated into the same algorithm independently, thus allowing them to be considered -- and their performance assessed -- individually or in combination. Therefore, whilst amongst the substantial number of population-based metaheuristics available, see e.g. \cite{Ghazali2009}, there may be other suitable frameworks, a PSO approach clearly meets our requirements.

Here, therefore, we seek to develop enhanced robust PSO-based approaches, based on adapting key elements of the d.d.\ and LEH approaches, combined with novel features. To that end we first consider PSO in more detail.

\subsection{Nominal PSO}
\label{sec:NomRobPSO}

There are many formulations of PSO, see \cite{KennedyEberhartShi2001, Kameyama2009, ZhangWangJi2015, SenguptaBasakPeters2018}. Here we describe one of the simplest, original formulations \cite{KennedyEberhart1995, KennedyEberhartShi2001, Ghazali2009}. This will form the basis of the robust approaches to be developed here. We will first consider a `standard', non-robust approach. A problem of the form (MM) can be considered in terms of its two constituent components: an inner maximisation and an outer minimisation. The PSO formulation described here should be appreciated as performing the outer minimisation component of (MM). We will return to the inner maximisation component when we discuss extending PSO to a baseline robust approach rPSO.

PSO starts at iteration $t=0$ with a population of $\Np$ particles at randomly selected points $\pmb{x}^j(0)$ in $\X$, where $j=1, \ldots, \Np$. The function is evaluated at these points. For each particle the best position it has visited is designated $\pb^j$, that is the position with the lowest objective function value $\tilde{g}(\pb^j)$.

Particles are interconnected for information sharing, so each particle has an associated neighbourhood of other particles within which information can be shared. Different PSO formulations employ different neighbourhood strategies. Here we use global PSO as described in \cite{ShiEberhart1998}. In global PSO the neighbourhood is the entire swarm and the information shared within the swarm is the global best position, that is the position $\gb$ in $\X$ with the lowest objective function value of all the points visited by all particles over all iterations.

From a particle's position $\pmb{x}^j(t)$ at iteration $t$ the particle's position is updated by the addition of its velocity vector $\pmb{v}^j$:
\[
	\pmb{x}^j(t)=\pmb{x}^j(t-1) \> + \> \pmb{v}^j(t)  \tag{Move}
	\]
Following the recommendation of \cite{Engelbrecht2012} to initialize particle velocities at zero or at random values close to zero, here we consider the approach where each $\pmb{v}^j(0)$ is separately initialised using uniform random sampling $\sim U(0\> , \>0.1)^n$. Beyond initialisation the following velocity formulation is used:
\[
	\pmb{v}^j(t)=\omega \cdot \pmb{v}^j(t-1) \> + \> C_1\cdot\pmb{r}_1\cdot(\pb^j-\pmb{x}^j(t-1)) \> + \> C_2\cdot\pmb{r}_2\cdot(\gb-\pmb{x}^j(t-1))  \tag{Vel1}
\]
Here $\pmb{r}_1 \> , \> \pmb{r}_2 \sim U(0\> , \>1)^n$, that is each component of the random vectors $\pmb{r}$ are randomly sampled individually, and multiplication between vectors is meant component wise. $C_1$, $C_2$ and $\omega$ are scalar terms. $C_1$ and $C_2$ represent `learning' factors that weight the priority that a particle puts on its own ($C_1$) versus the global ($C_2$) historic success (that is over all iterations, to date). $\omega$ is an inertia term which moderates the effect of the preceding velocity on the current velocity.

As the particles move through $\X$ their individual $\pb^j$ values and the global $\gb$ are updated as appropriate. If at any stage the next candidate position for any particle lies outside of the lower and upper bounds $l_i$ and $u_i$ of $\X$, here an invisible boundary condition is assumed, see \cite{RobinsonRahmatSamii2004}. Particles are allowed to leave the feasible region to naturally return to feasibility due to the pull of the $\pb^j$ and $\gb$ information. Note that when a candidate moves outside of the feasible region no function evaluations are undertaken. Rather the velocity equation is updated by the particle's new location, with the $\pb^j$ information remaining unchanged.

A standard PSO search can be extended to a baseline robust PSO search by adding an inner maximisation search component to the outer PSO minimisation search. Indeed it should be recognised that d.d.\ and LEH as described in Sections~\ref{sec:localRobustDD} and~\ref{sec:globalRobustLEH} focus on the outer minimisation component of the min max search (MM), so we will now give some consideration to inner maximisation.

\subsection{Inner maximisation}
\label{sec:InnerMax}

The key requirement of any inner maximisation approach is the ability to accurately identify the maximum objective function value within the $\Gamma$-radius uncertainty neighbourhood around any given candidate point. However when dealing with real-world problems we must additionally take account of practical considerations. For simulation problems a common limiting feature is the number of model runs that can be performed, primarily due to simulation run times. In such a situation it is common to be restricted to an upper budget $\budget$ on the number of model runs, which would in turn impact on the ability to accurately perform inner maximisation searches.

Where there are budgetary restrictions on the number of function evaluations, some trade-off must be achieved between the extent of each inner maximisation search (robustness) and the overall global search performance. However the trade-off between robustness and performance is not straightforward, see \cite{MirjaliliLewisMostaghim2015, DiazHandlXu2017}. In \cite{BertsimasNohadaniTeo2010} the inner maximisation involves a series of two-stage gradient ascent searches within the $\Gamma$-uncertainty neighbourhood of a given candidate point, and assumes the availability of gradient information. Such an approach to the inner maximisation is comprehensive, but is in practice likely to prove prohibitive with increasing number of dimensions, even assuming the availability of gradient information.

In \cite{HughesGoerigkWright2019} uniform random sampling is used in the $\Gamma$-radius hypersphere that forms the uncertainty neighbourhood around a candidate point, with the maximum value sampled taken as an approximation to the inner maximum. It is this approach that we adopt here for the inner maximisation in all heuristics considered. Pseudo code for inner maximisation by uniform random sampling in a $\Gamma$-radius hypersphere is shown in Algorithm~\ref{InnerAlgorithm}. In the following, we do not explicitly list $f$, $\X$, or $\Gamma$ as algorithm inputs, as they are always implied.


\begin{algorithm}[htbp] 
\caption{$\Gamma$-uncertainty neighbourhood inner maximisation inc. STOPPING option} \label{InnerAlgorithm}
\vspace{2mm} 
\hspace*{\algorithmicindent} \textbf{Input:} $\pmb{x}_c$, $\budget$, $\inner$, $stopping$, $\tau$
\vspace{2mm} 
 
\begin{algorithmic}[1]

	\State Calculate $f(\pmb{x}_c)$ and store in $F_H$ \label{Point1}	
	\State $H \gets H \cup \{\pmb{x}_c\}$
	\State $\tilde{g}(\pmb{x}_c) \gets f(\pmb{x_c})$ \label{UpdateRobust1}
	
	\State $\budget \gets \budget-1$ 
	
	\IIf{($stopping$) AND ($\tilde{g}(\pmb{x}_c) >\tau$)} goto line~\ref{ReturnOutput1} \EndIIf \label{Threshold1}	
		
	\IIf{($\budget=0$)} \textbf{break}: goto end of Outer Min algorithm \EndIIf \label{InnerFailedToComplete1}

	\ForAll{$i$ in ($1,\ldots,\inner - 1$)} \label{InnerStart}
	
		\State Choose $\Delta\pmb{x}^i_c \in \U$ uniformly at random \label{RandomInHyper1}
		\State$\pmb{x}_c^i \gets \pmb{x}_c + \Delta\pmb{x}^i_c$ \label{RandomInHyper2}
		\State Calculate $f(\pmb{x}^i_c)$ and store in $F_H$ \label{RandomInHyper3}
		\State $H \gets H \cup \{\pmb{x}^i_c\}$
		\State $\tilde{g}(\pmb{x}_c)  \gets \max\{\tilde{g}(\pmb{x}_c) ,f(\pmb{x}^i_c)\}$ \label{UpdateRobust2}
		
		\State $\budget \gets \budget-1$ 
		
		\IIf{($stopping$) AND ($\tilde{g}(\pmb{x}_c) >\tau$)} goto line~\ref{ReturnOutput1} \EndIIf \label{Threshold2}

		\IIf{($\budget=0$)} \textbf{break}: goto end of Outer Min algorithm \EndIIf \label{InnerFailedToComplete2}

	\EndFor  \label{InnerEnd} 

	\State \Return $\tilde{g}(\pmb{x}_c)$: estimated worst case cost at $\pmb{x}_c$  \label{ReturnOutput1}
	
\end{algorithmic}
\end{algorithm}

Within any robust heuristic, given a candidate point $\pmb{x}_c$ around which we want to perform a $\Gamma$-uncertainty neighbourhood inner maximisation, we call Algorithm~\ref{InnerAlgorithm}. As input this requires the point information $\pmb{x}_c$, a maximum number of function evaluations that can be undertaken in the entire search budget $\budget$, the defined number of points to be evaluated within the inner maximisation analysis $\inner$ (in the case of stopping this is the maximum number of points that could be evaluated), a boolean specifying whether or not the stopping condition is to be invoked $stopping$, and if required the stopping threshold $\tau$.

It should be noted that $\inner$ is a parameter that is tuned for all heuristics, within the experimental testing here, see Section~\ref{sec:parameterTuning}.

The algorithm starts by evaluating the function at the candidate point $\pmb{x}_c$ (line~\ref{Point1}), prior to moving on to uniformly randomly sample points in the $\Gamma$-uncertainty neighbourhood of $\pmb{x}_c$ and evaluating the function at each of these points (lines~\ref{RandomInHyper1} to~\ref{RandomInHyper3}). Function evaluations are recorded in the set $F_H$ associated with $H$. When a function evaluation is performed the budget counter is reduced by 1 and a check is performed to ensure that $\budget$ has not been exceeded. Note, however, that when the inner maximisation analysis has been prematurely ended due to $\budget$  being exhausted, we do not want to return an estimate for $\tilde{g}(\pmb{x}_c)$, but instead return to the end of the outer minimisation algorithm where the extant estimate for the robust global minimum is accepted (lines~\ref{InnerFailedToComplete1} and~\ref{InnerFailedToComplete2}). 

As the inner sampling proceeds the estimate for $\tilde{g}(\pmb{x}_c)$ is updated as appropriate (lines~\ref{UpdateRobust1} and~\ref{UpdateRobust2}). If the input value for $stopping$ is TRUE and a function evaluation is detected which exceeds $\tau$, the inner maximisation is terminated with the current estimate for $\tilde{g}(\pmb{x}_c)$ returned (lines~\ref{Threshold1} and~\ref{Threshold2}). Otherwise the full inner maximisation is completed, at which point the estimate for $\tilde{g}(\pmb{x}_c)$ is returned to the outer minimisation algorithm (line~\ref{ReturnOutput1}).

\subsection{Baseline robust PSO}
\label{sec:bruteForcePSO}

The easiest way to extend a PSO approach to an rPSO version and tackle the problem (MM) is to perform an inner maximisation search around any point in $\X$ visited by a particle, in order to replace the nominal objective function value $f(\pmb{x})$ at any given point with the corresponding worst case cost value $\tilde{g}(\pmb{x})$. With this approach the PSO formulation would remain unchanged. This is a baseline rPSO, which is used as the starting point for developing enhanced rPSO metaheuristics here.

Pseudo code for this baseline rPSO heuristic is given in Algorithm~\ref{BaselineRPSO} for defined input parameter values $\Np$, $C_1$, $C_2$, $\omega$, and $\inner$. In the experimental testing, and for all rPSO based heuristics, these five parameters are tuned, see Section~\ref{sec:parameterTuning}. Note that we do not need to define the number of iterations over which the swarm is progressed, as this will be controlled by $\budget$ within the inner maximisation Algorithm~\ref{InnerAlgorithm}, which is called from line~\ref{InnerSub} of Algorithm~\ref{BaselineRPSO}. In addition we must input information for $\budget$. 

The PSO algorithm loops over iterations $i$ until the budget $\budget$ is exceeded (line~\ref{Iterations}), and over the $\Np$ particles in the swarm (line~\ref{Swarm}). At the first iteration the particles are randomly initialised (line~\ref{RandomInitialise}), but for subsequent iterations the particle positions, velocities, and personal best information are updated according to equations (Vel1) and (Move) (lines~\ref{UpdateVelocity} and ~\ref{MoveParticle}). Prior to performing any inner maximisation function evaluations the feasibility of the candidate point $\pmb{x}^j(i)$ is confirmed. Here we use the boolean $ToBeEvaluated$ (lines~\ref{FeasibilityFlag} and~\ref{Feasibility1}) to flag feasibility. Note that the use of the flag $ToBeEvaluated$ is exploited further subsequently in Algorithm~\ref{partLEHAlgorithm}.

If $\pmb{x}^j(i)$ is not in $\X$ the inner maximisation, and associated function evaluations, are skipped (line~\ref{Feasibility2}). This means that the particles personal best information $\pb^j$ will not be updated, but otherwise the subsequent movement of particle $\pmb{x}^j$ will continue according to equations (Vel1) and (Move).

Particle's $\pb^j$ and the estimate of the robust global optimum $\gb$ are updated as appropriate (lines~\ref{PBestInit} and ~\ref{RobGlobInit}). At the end of the swarm-iterations loops the extant estimate of the robust global optimum $\gb$ is returned (line~\ref{ReturnOutput2}).


\begin{algorithm}[htbp] 
\caption{A baseline robust particle swarm optimisation algorithm} \label{BaselineRPSO}
\vspace{2mm} 
\hspace*{\algorithmicindent} \textbf{Input:} $\budget$ \\
\hspace*{\algorithmicindent} \textbf{Parameters:} $\Np$, $C_1$, $C_2$, $\omega$, $\inner$
\vspace{2mm} 
 
\begin{algorithmic}[1]

	\While{($\budget>0$)} \label{Iterations}

		\State $i\gets 0$ 
	
		\ForAll{($j$ in $1,\ldots,\Np$)} \label{Swarm}

			\State $ToBeEvaluated \gets$ TRUE \label{FeasibilityFlag}
			\If{($i = 0$)}  \label{FirstIter}
				
				\State Choose $\pmb{x}^j(i) \in \X$ uniformly at random \label{RandomInitialise}
			
			\Else  \label{NotFirstIter} 	
			
				\State Update particle velocity $\pmb{v}^j(i)$ according to (Vel1) \label{UpdateVelocity}
				\State Update particle position $\pmb{x}^j(i)$ according to (Move)  \label{MoveParticle}
				\IIf{($\pmb{x}^j(i) \notin \X$)} $ToBeEvaluated \gets$ FALSE \EndIIf \label{Feasibility1}				
			
			\EndIf  \label{FirstIterEnd} 	
							
			\If{($ToBeEvaluated$)}  \label{Feasibility2}

				\State $\tilde{g}(\pmb{x}^j(i)) \gets$ \textbf{CALL} Algorithm~\ref{InnerAlgorithm}($\pmb{x}^j(i)$, $\budget$, $\inner$, FALSE, 0) \label{InnerSub}	
									
				\IIf{($i = 0$) OR ($\tilde{g}(\pmb{x}^j(i)) < \tilde{g}(\pb^j)$)}  
					 $\pb^j \gets \pmb{x}^j(i)$ \EndIIf			\label{PBestInit}

				\IIf{($i = 0$ AND $j = 0$) OR ($\tilde{g}(\pmb{x}^j(i)) < \tilde{g}(\gb)$)}  
					 $\gb \gets \pmb{x}^j(i)$	 \EndIIf \label{RobGlobInit}									
			
			\EndIf  \label{FeasibilityEnd} 								

		\EndFor  \label{SwarmEnd} 		

		\State $i\gets i+1$ 
		
	\EndWhile  \label{IterationsEnd} 

	\State \Return A robust solution $\gb$ \label{ReturnOutput2}
	
\end{algorithmic}
\end{algorithm}

\medskip


\section{Comparison of baseline heuristics}
\label{sec:baselineComparators}

When it comes to the testing of new heuristics in Section~\ref{sec:experimentsResults} we require comparator robust heuristics against which to assess performance. Here we use the three baseline approaches already discussed: re-starting d.d.\, LEH and the baseline rPSO. In order to give some indication of the different natures of the searches due to each comparator robust metaheuristic considered we introduce a two dimensional problem and plot exemplar searches due to each heuristic.

Consider one of the test problems to be used in our experimental test suite Section~\ref{sec:setUp}, the multi-dimensional Pickelhaube problem. A full description of this function is given in Appendix~\ref{sec:testFunctionFormulae}, and plots of the nominal and worst case ($\Gamma$=1) 2D Pickelhaube problem are shown in Figures~\ref{fig:PickelhaubeNom} and~\ref{fig:PickelhaubeWorst}. In our formulation for the 2D problem the nominal global optimum  is at $(-35, -35)$, whilst the robust global optimum is at $(-25, -25)$. 

\begin{figure}[htbp]
	\centering
	
	\begin{subfigure}{.38\textwidth}
		\centering
		\includegraphics[width=2.5in, height=2.6in]{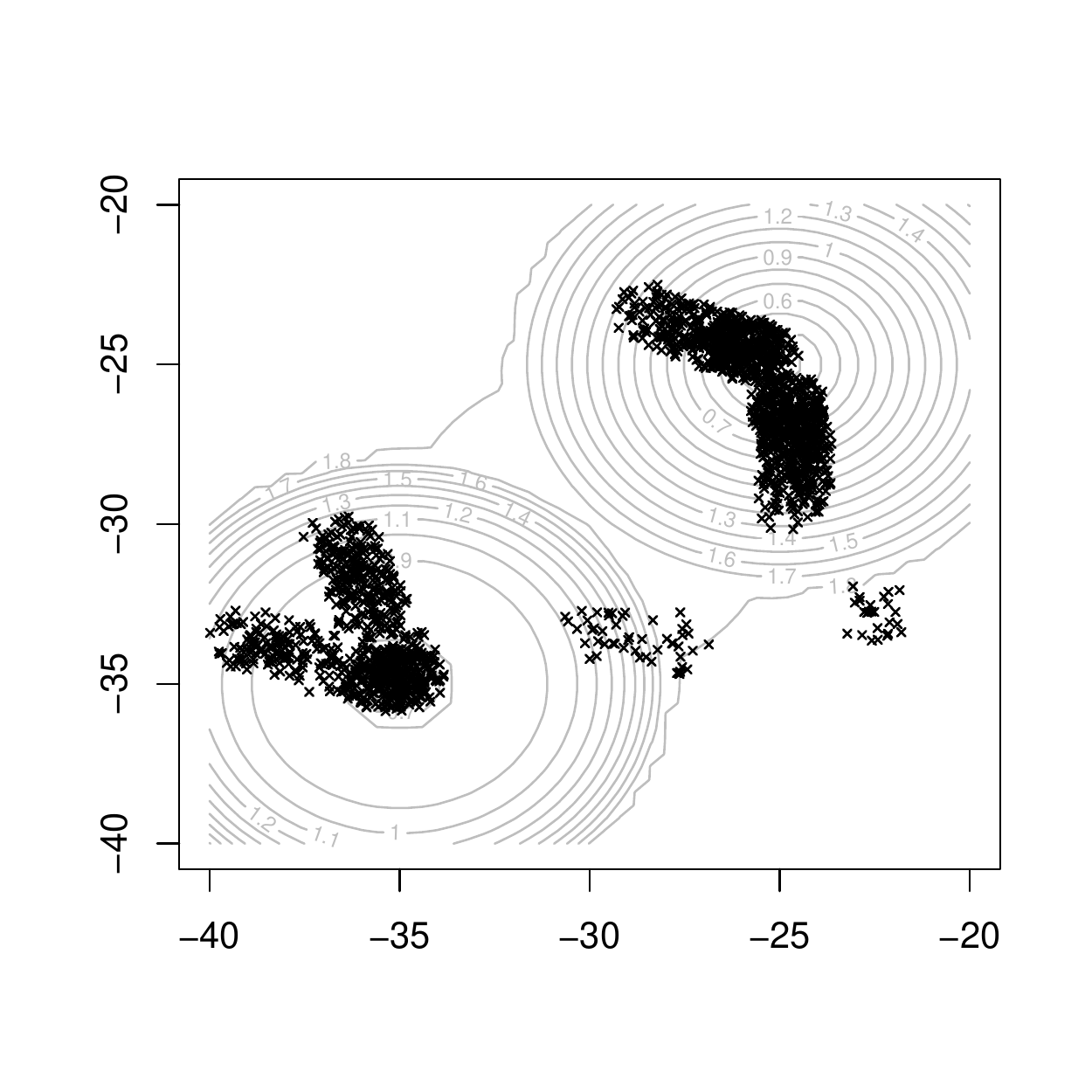}	
		\vspace{-11mm} 
		\caption{d.d.\ points}
		\label{fig:DDPoints}
	\end{subfigure}%
	\hspace{7mm} 
	\begin{subfigure}{.38\textwidth}
		\centering
		\includegraphics[width=2.5in, height=2.6in]{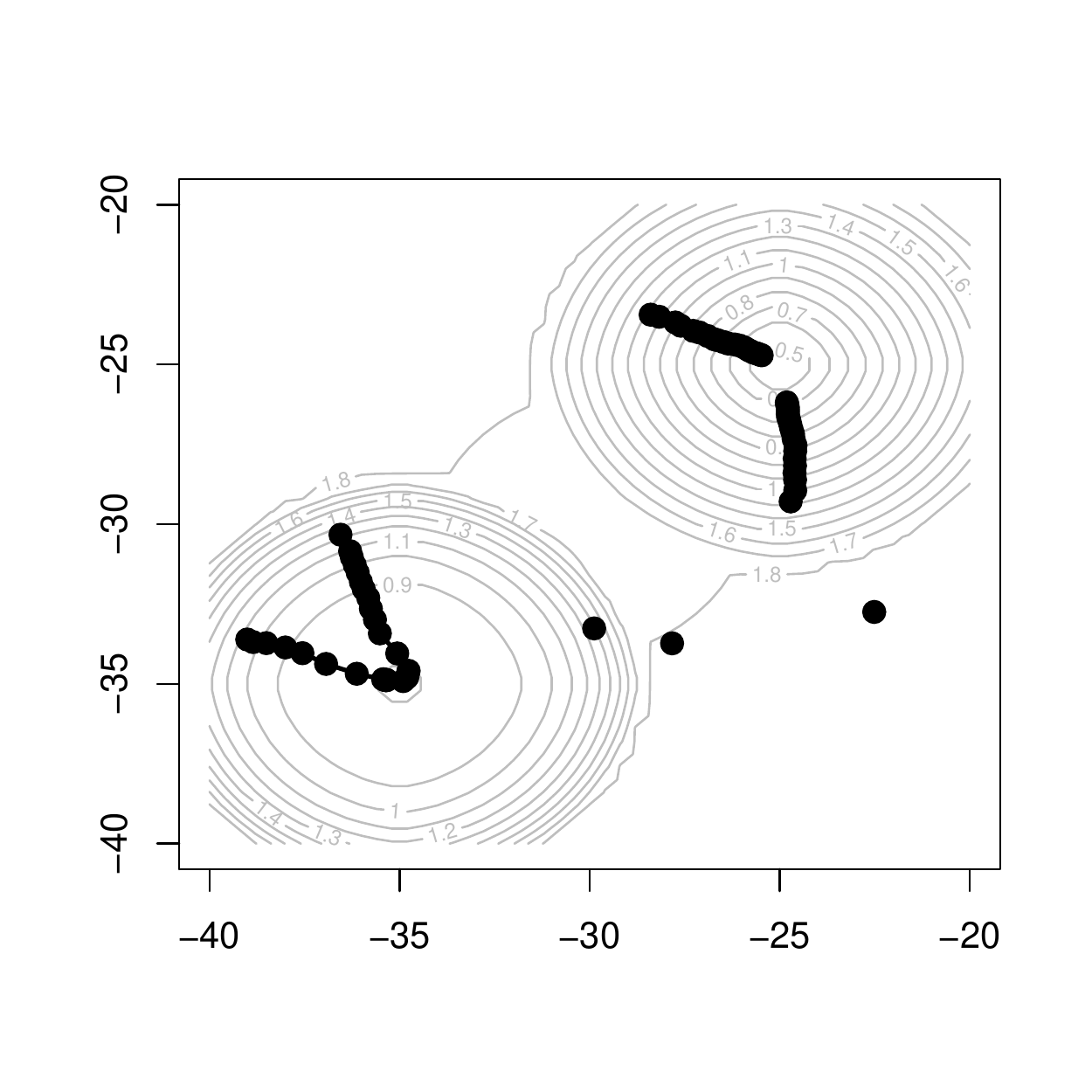}	
		\vspace{-11mm} 
  	\caption{d.d.\ current best}
		\label{fig:DDSearch}
	\end{subfigure}	
		
	\vspace{-5mm} 
			
	\begin{subfigure}{.38\textwidth}
		\centering
		\includegraphics[width=2.5in, height=2.6in]{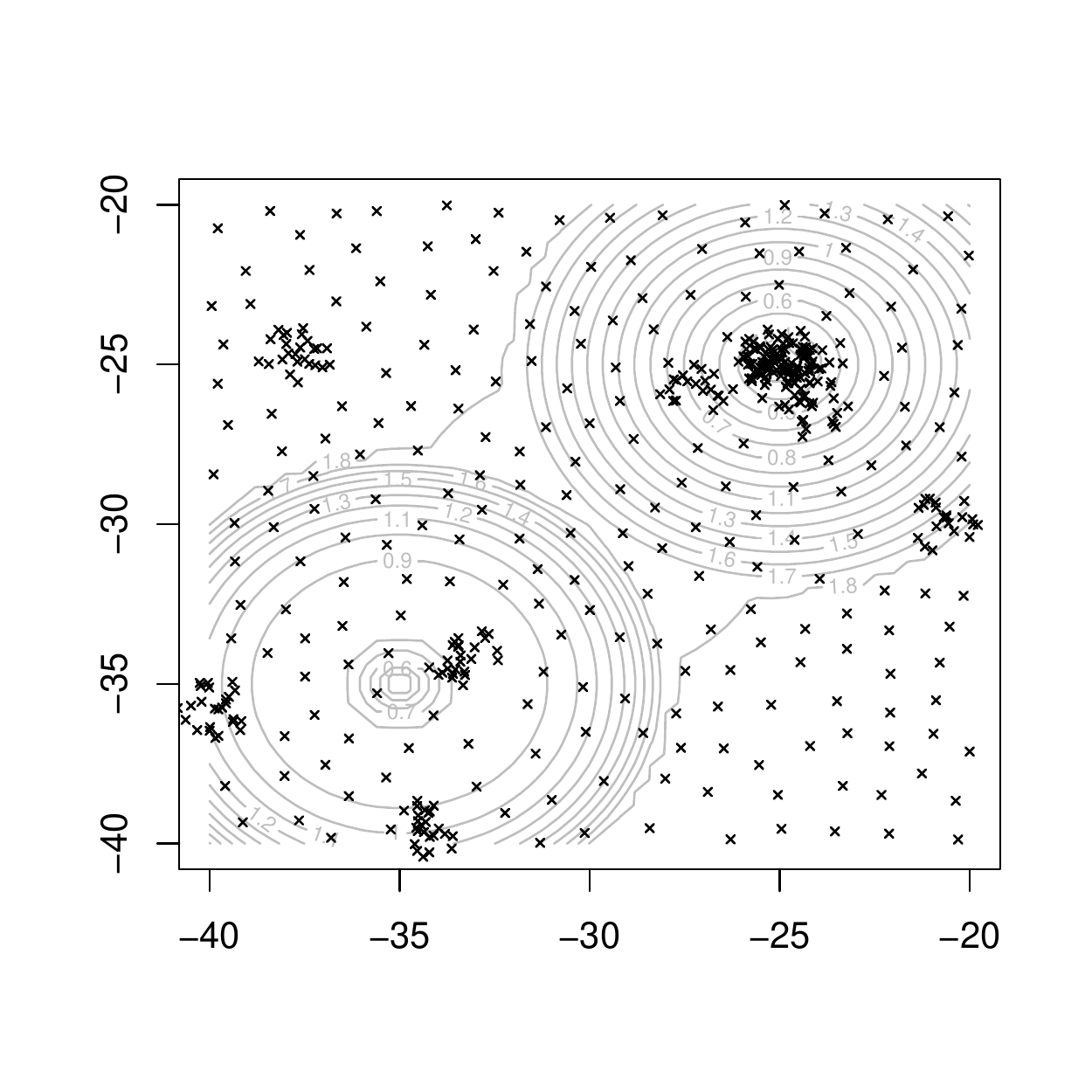}		
		\vspace{-11mm} 
	  \caption{LEH points}	
		\label{fig:LEHPoints}
	\end{subfigure}%
	\hspace{7mm} 
	\begin{subfigure}{.38\textwidth}
		\centering
		\includegraphics[width=2.5in, height=2.6in]{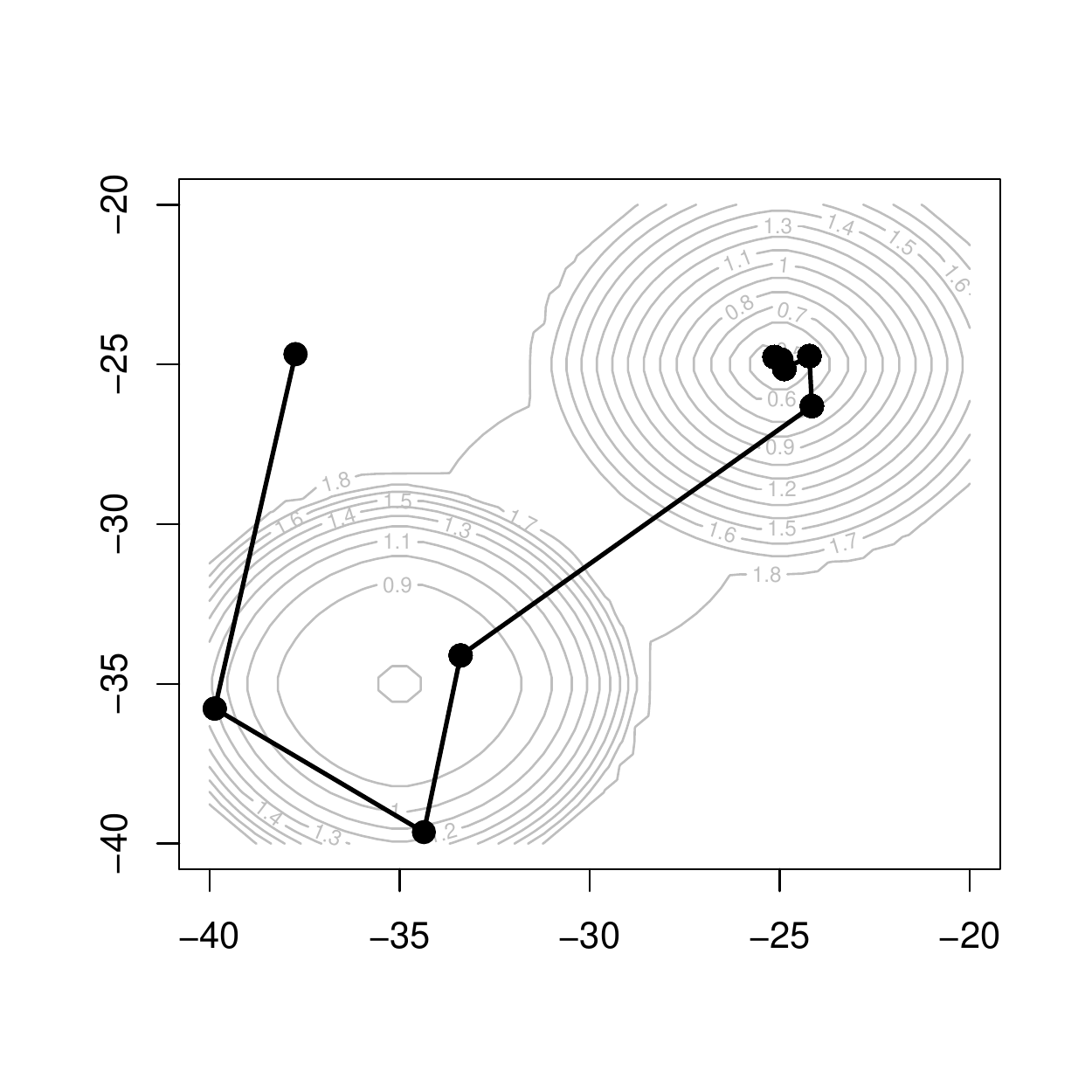}		
		\vspace{-11mm} 
  	\caption{LEH current best}
		\label{fig:LEHSearch}
	\end{subfigure}
		
	\vspace{-5mm} 
	
	\begin{subfigure}{.38\textwidth} 
		\centering
		\includegraphics[width=2.5in, height=2.6in]{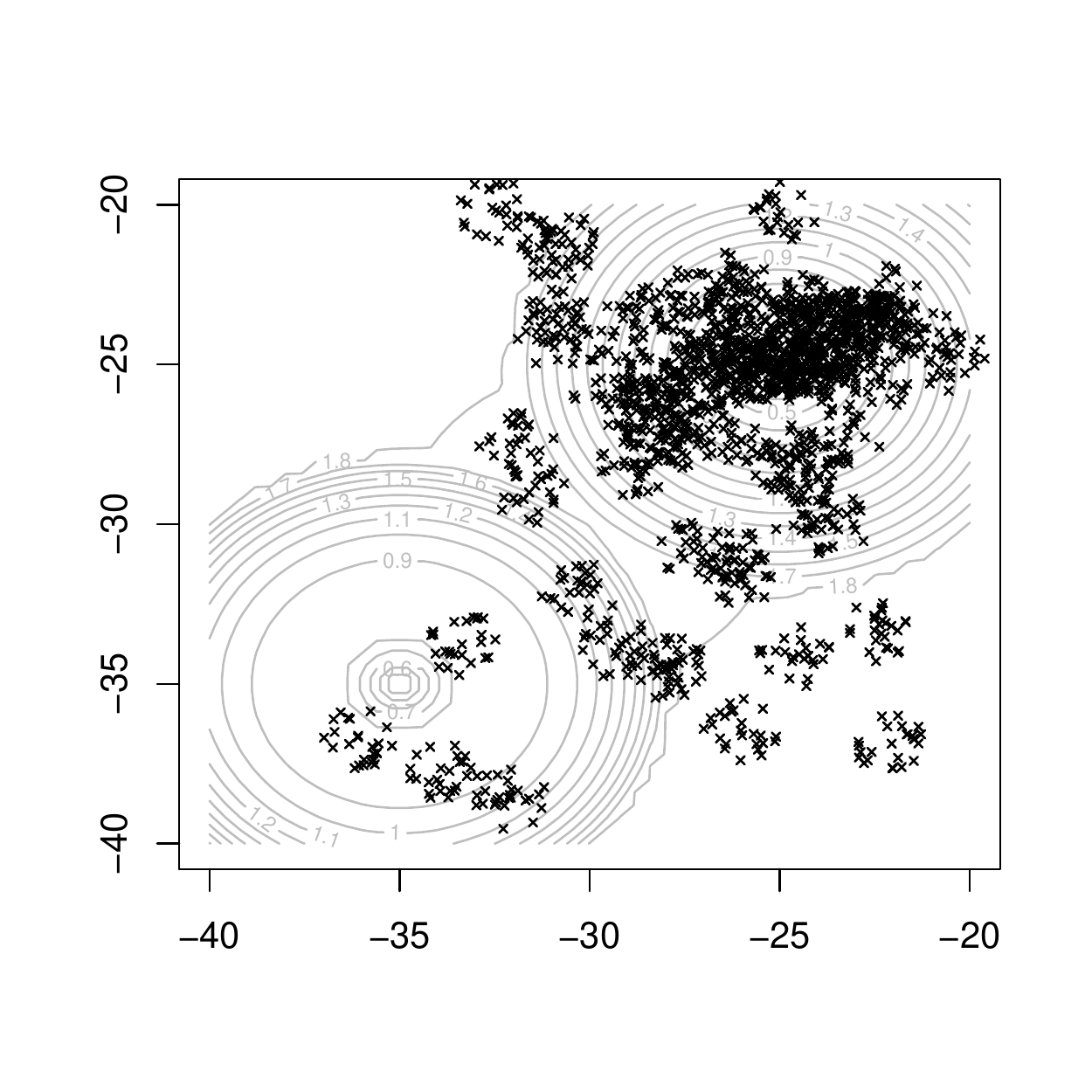}
		\vspace{-11mm} 
  	\caption{rPSO points}
		\label{fig:PSOPoints}
	\end{subfigure}%
	\hspace{7mm} 
	\begin{subfigure}{.38\textwidth}
		\centering
		\includegraphics[width=2.5in, height=2.6in]{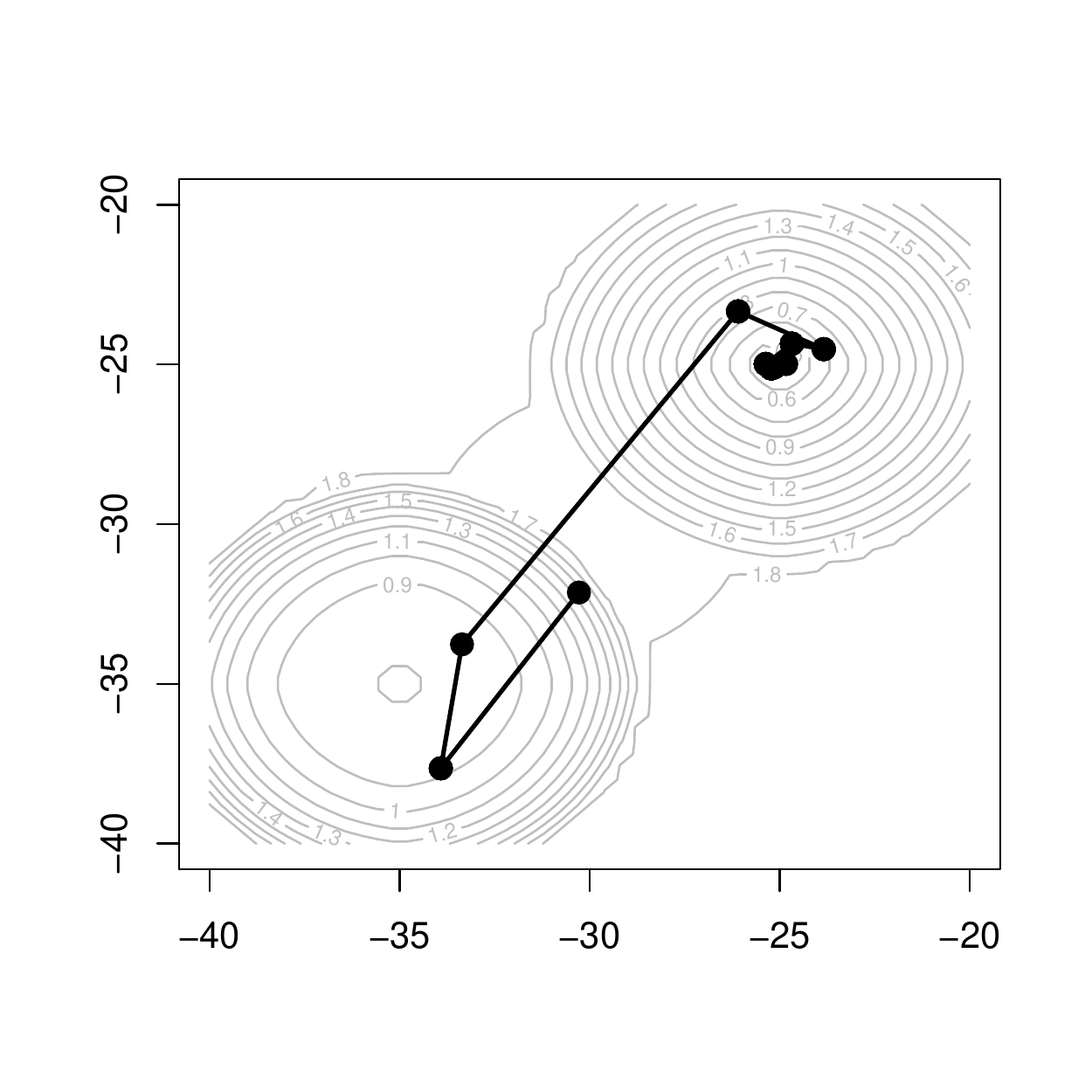}	
		\vspace{-11mm} 
  	\caption{rPSO current best}
		\label{fig:PSOSearch}
	\end{subfigure}
	
	\caption{Example searches of the 2D version of the Pickelhaube problem (see Figures~\ref{fig:PickelhaubeNom} and ~\ref{fig:PickelhaubeWorst}) with $\Gamma=1$, for the baseline metaheuristics. Evaluated points are shown on the left, with the path of the improving current robust best on the right. The outer heuristics are (top to bottom): d.d.\ with re-start, LEH, and robust PSO. The nominal global optimum is at $(-35, -35)$, and the robust global optimum is at $\sim (-25, -25)$.}
	\label{fig:BaselineRobust2DSearches}
\end{figure}

Contour plots of individual example searches for each of the three baseline approaches applied to this 2D problem are shown in Figure~\ref{fig:BaselineRobust2DSearches}. These plot should be seen as indicative exemplars. Plots on the left indicate all points evaluated and the underlying contour is the nominal plot. Plots on the right show the improving search path of the currently identified global robust optima, over the underlying worst case contour.

For the d.d.\ search in Figures~\ref{fig:DDPoints} and~\ref{fig:DDSearch} the inner maximisation groupings of evaluated points can be seen to follow a series of paths, each towards a robust local optimum and based on a series random re-starts at the completion of the previous local search. It can be seen that two such local searches successfully hone in on the robust global optimum.

The nature of the LEH search in Figures~\ref{fig:LEHPoints} and~\ref{fig:LEHSearch} is very different to the d.d.\ search, and two obvious features are apparent. First is the extensive exploration of the solution space as the search jumps to centres of regions devoid of poor points. Second is that in many cases the inner maximisation groupings of evaluated points are sparse, and in fact are often just single point evaluations as the stopping condition recognises that the current location cannot improve on the existing estimate of a global robust optima -- and so the algorithm immediately moves on without undertaking a full $\Gamma$-uncertainty neighbourhood search. The location of the robust global optima is successfully identified.

Finally for the baseline rPSO search Figures~\ref{fig:PSOPoints} and~\ref{fig:PSOSearch} the combined exploratory-exploitation nature of a PSO search can be observed. The approximate location  of each particle within each iteration is identifiable by the cluster of points evaluated within that particle's $\Gamma$-uncertainty neighbourhood. Again the location of the robust global optima is successfully identified.

\section{Enhanced robust particle swarm optimisation}
\label{sec:EnhancedPSO}

In both d.d.\ and LEH approaches the additional function evaluations required to calculate robust, as opposed to nominal, values are used to direct the search. In the case of d.d.\ this points the search towards the optimal local direction to avoid hcps, whilst for LEH this points the search towards the location which is globally furthest from all existing hcps. Here we seek to exploit both of these local and global search directions, separately and combined, in the form of three new enhanced rPSO metaheuristics:

\begin{itemize}[leftmargin=*]	
	\item rPSO + d.d.: Baseline rPSO velocity formulation with the addition of a d.d.\ component (rPSOdd).
	
	The basic PSO movement formulation involves a weighted combination of two direction vectors, added to a particle's current position, generating the movement to a new location. The two weighted direction vectors are based on the differences between a particle's current location and its historic personal best and, here, the current estimated global best. An obvious approach, given a third piece of vector information -- a local uncertainty neighbourhood d.d.\ vector --  is to simply add the weighted d.d.\ vector to the original vector calculation. Here each d.d.\ vector can be calculated at the particle level at each candidate point, given an inner maximisation at that point. Furthermore parameter tuning including all three weighting parameters will enable identification of the best combination of vector information.

	\item rPSO + LEH:  Using historic information, the stopping condition and LEH in the baseline rPSO algorithm (rPSOleh).
	
  In the LEH formulation at each new candidate point an inner maximisation begins but can terminate early if an uncertainty neighbourhood point exceeds the currently estimated robust global minimum. At its most efficient this may frequently mean only a single function evaluation is undertaken at a candidate, if it is immediately determined that that evaluation exceeds the global best. Incorporating this into PSO at a particle level, using particle best information as a stopping threshold, will introduce the desired inner maximisation efficiency savings. However this can be taken further. For LEH in a restricted budget on function evaluations setting, the historic record of function evaluations is unlikely to include any points in the uncertainty neighbourhood of a candidate location, as the very nature of LEH is to move to previously unexplored regions. In a population-based approach however this may not be the case, and in fact it may be desirable that there is some convergence of members of the population in the decision variable space. This introduces the potential that no function evaluations may need to be undertaken at an individual particle's candidate location if the historic record identifies a previously evaluated point in the candidate uncertainty neighbourhood with objective function value greater than the particle's threshold. This is an additional efficiency, but also opens up the possibility of a further feature.

	The non-requirement to perform even a single function evaluation at an individual particle's current candidate location allows for the consideration of a particle becoming 'dormant'. Not requiring any function evaluations at a location could be due to moving into an already visited region of the decision variable space and identifying neighbouring points in the historic record exceeding the particle's threshold, as discussed above. However alternatively this could be due to the particle moving outside of the feasible region. That is the invisible boundary condition \cite{RobinsonRahmatSamii2004} employed here also means that individual particles may become dormant in the sense that they move outside of the feasible region, necessitating no function evaluations and under the expectation of subsequent movements ultimately returning the particle to feasibility, driven by personal and neighbourhood best information.

	If either of the situations described above were to repeat over several swarm iterations it would be reasonable to consider an approach that interrupts individual particles that have become dormant: 'stuck' either in previously visited or infeasible areas. By introducing a dormancy threshold representing the number of iterations before a particle is deemed truly dormant and requiring further action, that dormancy threshold can be parameter-tuned.

	Given a particle's dormancy it seems reasonable to consider some action, so here we introduce an exploratory component. Taking the other main element of the LEH approach, the calculation of the largest empty region devoid of all previously evaluated points with objective value greater than some defined threshold, we relocate-reinitialise the particle at the centre of that LEH. Here we set the 'high cost' threshold, identifying which points from the historical record of function evaluations to avoid in the calculation of the largest empty region, equal to the current robust global minimum.
	
	\item rPSO + LEH + d.d.: This is a combination of all new features in our PSO framework: using historic information, the stopping condition and LEH in the baseline rPSO algorithm, as well as an additional d.d.\ component in the rPSO velocity formulation (rPSOlehdd).
\end{itemize}

Our enhanced rPSO metaheuristics use a framework of the baseline rPSO Algorithm~\ref{BaselineRPSO} in conjunction with new algorithms to execute the features described above. 

\subsection{rPSOdd}
\label{sec:rPSOdd}

We start by considering the use of d.d.\ information to enhance the PSO velocity equation (Vel1) through the addition of a further velocity component of the form:
\[
	C_3 \cdot \pmb{r}_3 \cdot \pmb{dd}^j(t-1) \tag{ddVel}
\]
For each particle $j$ at each iteration $t$ we perform a single-step d.d.\ calculation using the uncertainty neighbourhood points around the particle's position. Given a valid unit length direction vector $\pmb{d}^j(t-1)$ as a result of solving (Soc), a final vector $\pmb{dd}^j(t-1)$ is calculated by scaling using the term $\rho^j(t-1)$, which is calculated as in equation (Rho) in Section~\ref{sec:localRobustDD}:
\[
	\pmb{dd}^j(t-1)= \rho^j(t-1) \cdot \pmb{d}^j(t-1)
\]
If no such direction can be calculated then $\pmb{dd}^j(t-1)$ is set equal to the zero vector. Incorporating this additional velocity component into the original velocity formulation (Vel1), beyond the initialisation of the particles at iteration 0 the following velocity formulation is used:
\[
	\pmb{v}^j(t)=\omega \cdot \pmb{v}^j(t-1) \> + \> C_1 \cdot \pmb{r}_1 \cdot (\pb^j-\pmb{x}^j(t-1)) \> + \> C_2 \cdot \pmb{r}_2 \cdot (\gb-\pmb{x}^j(t-1) \> + \> C_3 \cdot \pmb{r}_3 \cdot \pmb{dd}^j(t-1) \tag{Vel2}
\]
Pseudo code for the calculation of the additional d.d.\ velocity component (ddVel) is given in Algorithm~\ref{ddVelComponent}. In terms of the overarching enhanced metaheuristic rPSOdd, Algorithm~\ref{BaselineRPSO} is still valid but with a single change: in line~\ref{UpdateVelocity} the enhanced velocity equation (Vel2) now replaces equation (Vel1).


\begin{algorithm}[htb] 
\caption{Calculating the additional d.d.\ velocity component (ddVel)} \label{ddVelComponent}
\vspace{2mm} 
\hspace*{\algorithmicindent} \textbf{Input:} $\pmb{x}_c$, $H$, $F_H$ \\
\hspace*{\algorithmicindent} \textbf{Parameters:} $\sigma$, $\sigma_{limit}$, $\sigma_{no}$, $C_3$
\vspace{2mm} 
 
\begin{algorithmic}[1]

	\State $\pmb{dd}_c \gets \pmb{0}$  \label{zeroInitialiseDD} 

	\If{($\pmb{x}_c \in \X$)}  \label{ddFeasibility}

		\State Initialise $H_\sigma$ and $F_{H_\sigma}$
		\label{initialHCPs}  
	
		\State $SolvedSOCP \gets \text{FALSE}$
	
		\While{(!$SolvedSOCP$) AND ($\sigma \geq \sigma_{limit}$)} \label{TrySOCP}	
	
			\State Try: $\pmb{d}_c \gets$ Solve (Soc) for $\pmb{x}_c$ and $H_\sigma$
			\IIf{(Solve (Soc) is successful)} $SolvedSOCP \gets \text{TRUE}$ \EndIIf	

			\If{($SolvedSOCP$)}  \label{SolvedSOCP}

				\State Calculate $\rho$ according to (Rho) \label{RoCalculation}
				\State $\pmb{d}_c \gets \rho\cdot\pmb{d}_c$

			\Else  \label{FailedSOCP} 	
			
				\State $\sigma \gets \sigma-(\sigma-\sigma_{limit})/\sigma_{no}$	 \label{sigmaReduction}
				\State Update $H_\sigma$ and the associated $F_{H_\sigma}$   \label{updateHCPs2}
			
			\EndIf  \label{SolvedSOCPEnd} 	

		\EndWhile  \label{TrySOCPEnd} 
	
		\IIf{($SolvedSOCP$)}  
			 $\pmb{dd}_c \gets C_3 \cdot \pmb{r_3} \cdot \pmb{d}_c$
    \EndIIf  \label{FinalSolvedSOCPEnd} 		

	\Else  \label{FailedddFeasibility}
		\ForAll{($i$ in $1,\ldots,n$)} \label{AdjustForInfeasible}
		
			\If{($(\pmb{x}_c)_i \le l_i$)}
				\State $ (\pmb{dd}_c)_i \gets \Gamma$
			\ElsIf{($(\pmb{x}_c)_i \ge u_i$)}
				\State $ (\pmb{dd}_c)_i \gets -\Gamma$			
			\EndIf 
		
		\EndFor  \label{AdjustForInfeasibleEnd} 	
		
		\State $\pmb{dd}_c \gets C_3 \cdot \pmb{r}_3 \cdot \pmb{dd}_c$   \label{ddVelCalcInFeasible}	

	\EndIf  \label{ddFeasibilityEnd} 		

	\State \Return Additional d.d.\ velocity component $\pmb{dd}_c$ 

\label{ReturnOutput3}
	
\end{algorithmic}
\end{algorithm}

Algorithm~\ref{ddVelComponent} requires input information for the candidate point of interest $\pmb{x}_c$, the history set $H$ and the associated set of function evaluation values $F_H$. In addition d.d.\ calculation parameter values are required: the initial $\sigma$ value, the lowest value this can take $\sigma_{limit}$, and the number of reduction-steps which can be applied in reducing $\sigma$ from its initial value down to $\sigma_{limit}$, in repeated attempts to solve (Soc) when the previous attempts have failed, $\sigma_{no}$ -- see Section~\ref{sec:localRobustDD}. Also the new input scalar parameter value $C_3$ is required for the final calculation of the $\pmb{dd}_c$ vector.

In the experimental testing the $\sigma$, $\sigma_{limit}$ and $C_3$ parameters are tuned, see Section~\ref{sec:parameterTuning}. This is in addition to the tuned baseline rPSO parameters, see Section~\ref{sec:bruteForcePSO}. $\sigma_{no}$ could be tuned but is pre-set in the experimental testing here in order to better control heuristic run times.

The algorithm begins with a feasibility check for the input candidate point $\pmb{x}_c$ (line~\ref{ddFeasibility}). If no function evaluations are undertaken at a given candidate point due to it lying outside of $\X$ (line~\ref{FailedddFeasibility}) then we  instead set $\pmb{dd}^j(t-1)$ to include only values in the dimensions which are infeasible (in all feasible dimensions the vector component retains the initialisation setting of 0, see line~\ref{zeroInitialiseDD}). The non-zero vector components are all set to magnitude $\Gamma$, with the sign for each dimension determined in order to point back into the feasible region (lines~\ref{AdjustForInfeasible} to~\ref{AdjustForInfeasibleEnd}). As with all d.d.\ vectors this is then multiplied by a scalar $C_3$ value and a random vector $\pmb{r_3}$ (line~\ref{ddVelCalcInFeasible}). The intention here is to promote a return to $\X$, beyond the existing draw of a particle's $\pb^j$ and the global $\gb$ information (which are both in $\X$).

The high cost set $H_\sigma$ and associated set of function evaluation values $F_{H_\sigma}$ are initialised (line~\ref{initialHCPs}). Next is the attempt to solve the second order cone problem (Soc) and identify a valid descent direction based on the high cost set $H_\sigma$, see lines~\ref{TrySOCP} to~\ref{TrySOCPEnd}. As this is a mathematical programming problem in practice this is achieved by a call to an optimisation software package. 

As described in Section~\ref{sec:localRobustDD}, if (Soc) cannot be solved immediately it may be re-tried multiple times with reducing values of $\sigma$ and hence with $H_\sigma$ containing fewer points (lines~\ref{sigmaReduction} and~\ref{updateHCPs2}). Flagging of the solution to (Soc) is controlled by the boolean $SolvedSOCP$. If (Soc) is solved the original normalised direction vector is re-scaled according to equation (Rho) (line~\ref{RoCalculation}), prior to the rPSO velocity equation update due to equation (ddVel) (line~\ref{FinalSolvedSOCPEnd}). If ultimately (Soc) cannot be solved the initialisation setting of $\pmb{dd}_c$ to $\pmb{0}$ is retained (line~\ref{zeroInitialiseDD}).

Contour plots of an example search using the rPSOdd heuristic applied to the 2D Pickelhaube problem are shown in Figures~\ref{fig:PSOddPoints} and~\ref{fig:PSOddSearch} on page~\pageref{fig:PSOddPoints}, to give an indication of the nature of a rPSOdd search. Unsurprisingly the nature of these plots is somewhat similar to the baseline rPSO search seen in Figures~\ref{fig:PSOPoints} and~\ref{fig:PSOSearch}. However in addition the new $\Gamma$-uncertainty neighbourhood descent directions component in the velocity function, for each particle at each location, does appear to add elements of robust local search at the particle level. The extent of any d.d.\ component will be heavily influenced by the rPSOdd parameter value settings.

\subsection{rPSOleh}
\label{sec:rPSOleh}

Our second enhanced rPSO formulation involves augmenting the baseline rPSO Algorithm~\ref{BaselineRPSO} with an additional algorithm to perform elements of the LEH approach due to \cite{HughesGoerigkWright2019}. Here we incorporate two components  into our enhanced heuristic rPSOleh: the stopping condition, and the calculation of the largest empty hypersphere devoid of high cost points and placement of candidates at the centre of the calculated LEH. This further leads to an increased role for the use of the historic function evaluation information from the history set $H$.

We have already included the stopping condition when developing the $\Gamma$-uncertainty neighbourhood inner maximisation Algorithm~\ref{InnerAlgorithm}. That algorithm is set up to accept boolean input information $stopping$ to flag whether or not the inner maximisation stopping condition should be invoked, plus the associated stopping threshold value $\tau$ if $stopping$ is TRUE. Whereas in the baseline rPSO Algorithm~\ref{BaselineRPSO} $stopping$ was set to FALSE (Algorithm~\ref{BaselineRPSO} line~\ref{InnerSub}), here the stopping condition is invoked at the particle level by setting $stopping$ to TRUE and using the candidate particle $j$ personal best value $\tilde{g}(\pb^j)$ for the stopping threshold $\tau$. So for rPSOleh, in line~\ref{InnerSub} of Algorithm~\ref{BaselineRPSO} the call to Algorithm~\ref{InnerAlgorithm} becomes:
\[
	\tilde{g}(\pmb{x}^j(i)) \gets \textbf{CALL} \> Algorithm~\ref{InnerAlgorithm} \> (\pmb{x}^j(i), \> \budget, \> \inner, \> TRUE, \> \tilde{g}(\pb^j)) \tag{call~\ref{InnerAlgorithm}}
\]
In addition, however, within our new Algorithm~\ref{partLEHAlgorithm}, we introduce a pre-inner maximisation check of the history set $H$. Again this algorithm requires input information for the candidate point of interest $\pmb{x}_c$, and the history set $H$ and associated set $F_H$. Further, Algorithm~\ref{partLEHAlgorithm} needs to access the full particle information for particle $j$($i$) associated with the current candidate point $\pmb{x}_c$, and the global best value $\tilde{g}(\gb)$. Also Algorithm~\ref{partLEHAlgorithm} uses the counter $dormancy_{count}^j$ and the parameters $ToBeEvaluated$ and $dormancy_{limit}$; the latter is tuned in our experiments along with the baseline rPSO parameters, Section~\ref{sec:parameterTuning}. A further parameter $placement_{limit}$ is also introduced here but is not tuned, and instead pre-set to control heuristic run times. $dormancy_{count}^j$, $dormancy_{limit}$ and $placement_{limit}$ are described below.

The pre-inner maximisation check of $H$ is to identify if there are existing points in the uncertainty neighbourhood of the candidate location $\pmb{x}_c$. If there are, we further check whether any such points already have nominal objective function value $f(\pmb{x})$ greater than the particle threshold $\tilde{g}(\pb^j)$. If so we determine that we do not need to perform any inner maximisation search for candidate location $\pmb{x}_c$. We have already seen the use of boolean $ToBeEvaluated$ to flag the need to undertake inner maximisation function evaluations, due to feasibility issues, in Algorithm~\ref{BaselineRPSO}. Here we extend the use of $ToBeEvaluated$ to additionally flag when there is no need to undertake inner maximisation function evaluations due to the the pre-inner maximisation check (lines~\ref{lehFeasibility} to~\ref{preStopEnd}). In a similar fashion to the case of $\pmb{x}_c$ being infeasible the particle $j$ velocity information is updated by the particle $j$ location $\pmb{x}_c$, but $\pb^j$ remains unchanged.


\begin{algorithm}[htbp] 
\caption{Re-locating particles using elements of the LEH heuristic}  \label{partLEHAlgorithm}
\vspace{2mm} 
\hspace*{\algorithmicindent} \textbf{Input:} Particle $j(i)$, $\pmb{x}_c$, $H$, $F_H$, $\gb$  \\
\hspace*{\algorithmicindent} \textbf{Parameters:} $ToBeEvaluated$, $dormancy_{count}^j$, $dormancy_{limit}$, $placement_{limit}$
\vspace{2mm} 
 
\begin{algorithmic}[1]

	\If{($\pmb{x}_c \notin \X$) OR ($\max\{\tilde{g}(\pmb{h}) \mid \pmb{h}\in N(\pmb{x}_c)\} > \tilde{g}(\pb^j)$)}  \label{lehFeasibility}
		\State $ToBeEvaluated \gets$ FALSE			
		\State $dormancy_{count}^j \gets dormancy_{count}^{j}+1$  \label{dormancyInc1}

	\EndIf  \label{preStopEnd} 	

	\If{($dormancy_{count}^j \> > \> dormancy_{limit}$)}  \label{dormancyCheck}
	
		\State $lehComplete \gets$ FALSE  \label{initialiseLEHflag}
		\State $countLEHtry\gets 0$  \label{initialiseCountLEH}						
		
		\While{(!$lehComplete$) AND ($countLEHtry \> < \> placement_{limit}$)} \label{TryLEH}	
		
			\State $\pmb{p} \gets$ solution to (lehMM) \label{solvelehMM}
			\State Calculate $f(\pmb{p})$ and store in $F_H$  \label{singlePointEval}
			\State $H \gets H \cup \{\pmb{p}\}$	
			\State $countLEHtry \gets countLEHtry+1$  \label{updateCountLEH}
				
			\IIf{($f(\pmb{p}) \> < \> \tilde{g}(\gb)$)} $lehComplete \gets \text{TRUE}$ \EndIIf  \label{updateLEHflag}	

		\EndWhile  \label{TryLEHEnd} 		
		
		\State Re-set all particle $j(i)$ details to initialisation values  \label{resetParticle}
		\State $dormancy_{count}^j \gets$ 0  \label{resetdormancyCount}		
		\State $\pmb{x}^j(i) \gets \pmb{p}$  \label{relocateParticle}	
		
	\Else
	
		\State Update particle position $\pmb{x}^j(i)$ according to (Move)  \label{MoveParticleLEH}	
	
	\EndIf  \label{dormancyCheckEnd} 			

	\State \Return $ToBeEvaluated$ and input particle $j(i)$ updated if appropriate \label{ReturnOutput4}
	
\end{algorithmic}
\end{algorithm}


The second feature of LEH that we exploit here is the exploration-based locating of new candidates at points furthest away from all previously visited `bad' high cost points -- which equates to placing candidates at the centre of the LEH (empty of hcps). Here we apply such an approach to relocate individual particles $j$ based on a determination that particle $j$ is `dormant'. Dormancy is based on a count of the number of swarm iterations over which no function evaluations have been performed for particle $j$, which may either be due to an infeasible candidate location, or due to the pre-inner maximisation check described  (lines~\ref{lehFeasibility} to~\ref{preStopEnd}). To this end we introduce a particle level count of the number of dormant iterations $dormancy_{count}^j$, and the dormancy count level which triggers the relocation of a particle $dormancy_{limit}$. The counter is incremented as appropriate (line~\ref{dormancyInc1}). 

For a particle $j$, given the exceeding of the dormancy count level (line~\ref{dormancyCheck}), an LEH calculation is undertaken in order to re-position particle $j$ to the centre of the identified LEH. If, however, particle $j$ is not identified as being dormant, the original update rule for the particle is used (line~\ref{MoveParticleLEH}).

In the original LEH algorithm, see Section~\ref{sec:globalRobustLEH}, the high cost set is defined as $H_\tau$, a subset of the history set $H$ containing all points with nominal objective function value $f(\pmb{x})$ greater than a threshold $\tau$. Within that algorithm the same threshold $\tau$ is employed for the stopping condition and for the identification of an LEH. Here we employ $\tau=\tilde{g}(\gb)$ with the intention of trying to re-locate particle $j$ away from points that we already know cannot improve on our current estimate for the robust global minimum $\gb$. 

We seek to estimate the point $\pmb{p}\in\X$ furthest from all designated high cost points $\pmb{h} \in H_{\tilde{g}(\gb)}$. This is the max min problem: 
\[
	\max_{\pmb{p}\in\X} \min_{\pmb{h}\in H_{\tilde{g}(\gb)}} \|\pmb{p}- \pmb{h}\|, \tag{lehMM}
\] 
where $\|\cdot\|$ is the Euclidean norm.

We use the approach due to \cite{HughesGoerigkWright2019} to estimate the solution to (lehMM), employing a genetic algorithm (line~\ref{solvelehMM}). As is the case for the baseline LEH comparator heuristic, within our experimental testing the parameters controlling the GA applied to solving problem (lehMM) are tuned, see Section~\ref{sec:parameterTuning}.

There is a final element of this particle relocation. In a further attempt to enhance the exploratory nature of the rPSOleh heuristic, at the potential new candidate point $\pmb{p}$ we perform a single point function evaluation (line~\ref{singlePointEval}). If this value is less than $\tilde{g}(\gb)$ we accept the new candidate point (line~\ref{updateLEHflag}), however if not we perform further LEH calculations up to some input number of times $placement_{limit}$ that this retry can occur (lines~\ref{TryLEH} to~\ref{TryLEHEnd}). This process is controlled by a counter, lines~\ref{initialiseCountLEH} and~\ref{updateCountLEH}, and success flag in line~\ref{initialiseLEHflag}. If the number of retries is exhausted the final potential candidate $\pmb{p}$ is accepted. Note that with each LEH calculation an additional point is added to $H$, impacting subsequent LEH calculations. 

As this LEH calculation effectively re-initialises particle $j$, the previous particle $j$ information including initial velocity and $\pb^j$ need to be re-set to the particle initialisation settings, overwriting the existing information (lines~\ref{resetParticle} and~\ref{resetdormancyCount}). Subsequently new particle information will be established as Algorithm~\ref{BaselineRPSO} progresses. Particle $j$ is then re-located to $\pmb{p}$ (line~\ref{relocateParticle}). This relocation could happen to the same particle $j$ more than once over the course of the heuristic search.

Algorithm~\ref{partLEHAlgorithm} can then be accessed from the baseline rPSO Algorithm~\ref{BaselineRPSO} by replacing lines~\ref{MoveParticle} and~\ref{Feasibility1} there with the single line reference to Algorithm~\ref{partLEHAlgorithm}:
\[
	\text{Update particle position } \pmb{x}^j(i) \text{ according to Algorithm~\ref{partLEHAlgorithm}} \tag{relocateLEH}
\]
In addition, to use Algorithm~\ref{BaselineRPSO} $dormancy_{count}^j$ values need to be introduced and initialised (set to zero) in Algorithm~\ref{BaselineRPSO}. The addition of a line within the  IF statement, lines~\ref{FirstIter} to~\ref{NotFirstIter}, achieves this:
\[
	dormancy_{count}^j\gets 0 \tag{addCounter}
\]
Also $dormancy_{limit}$ and $placement_{limit}$ would need to be defined as additional input parameters in Algorithm~\ref{BaselineRPSO}. 

\begin{figure}[htbp]
	\centering
	\begin{subfigure}{.38\textwidth} 
		\centering
		\includegraphics[width=2.5in, height=2.6in]{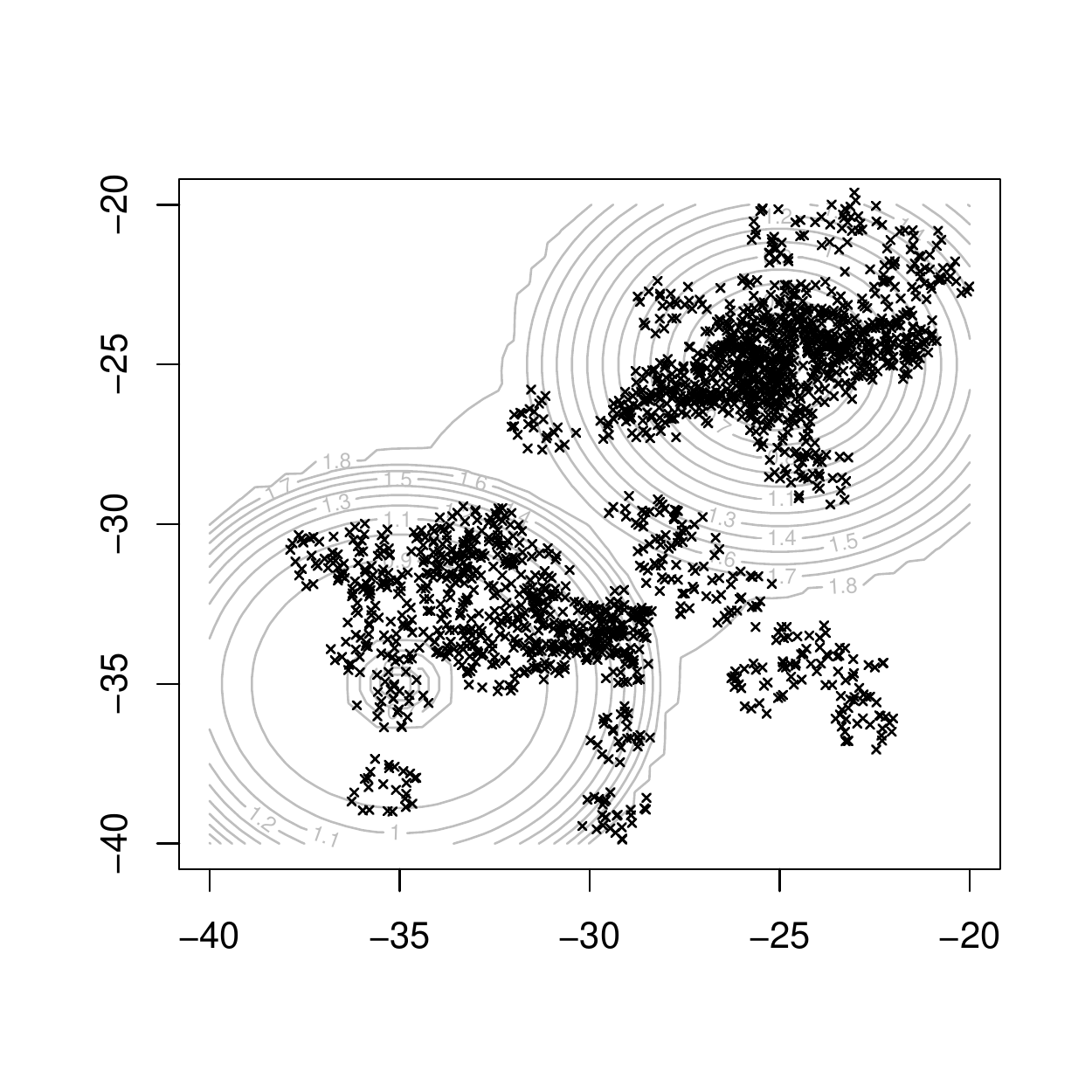}
		\vspace{-11mm} 
  	\caption{rPSOdd points}
		\label{fig:PSOddPoints}
	\end{subfigure}%
	\hspace{7mm} 
	\begin{subfigure}{.38\textwidth}
		\centering
		\includegraphics[width=2.5in, height=2.6in]{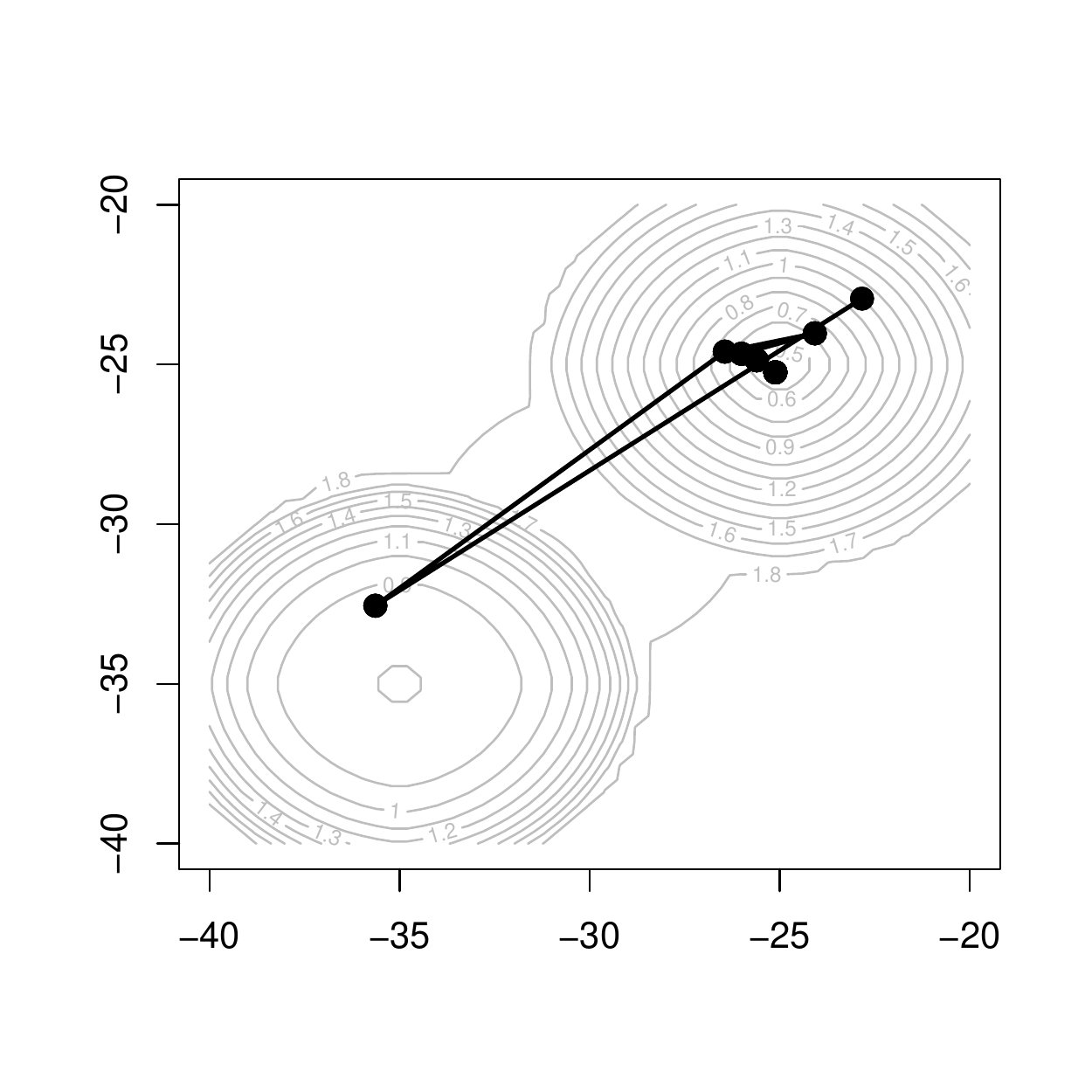}	
		\vspace{-11mm} 
  	\caption{rPSOdd current best}
		\label{fig:PSOddSearch}
	\end{subfigure}
		
	\vspace{-5mm} 
	
	\begin{subfigure}{.38\textwidth}
		\centering
		\includegraphics[width=2.5in, height=2.6in]{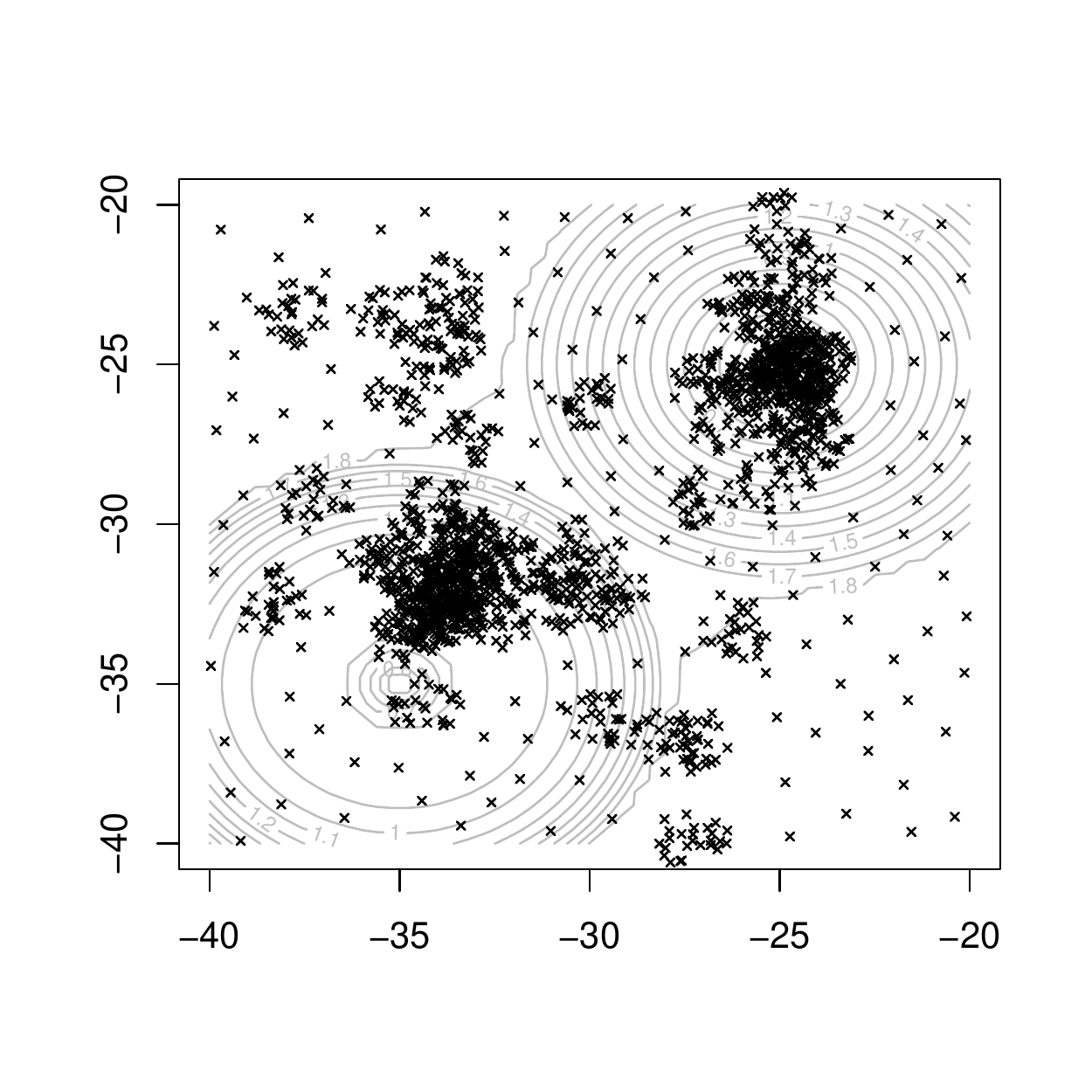}	
		\vspace{-11mm} 
		\caption{rPSOleh points}
		\label{fig:PSOlehPoints}
	\end{subfigure}%
	\hspace{7mm} 
	\begin{subfigure}{.38\textwidth}
		\centering
		\includegraphics[width=2.5in, height=2.6in]{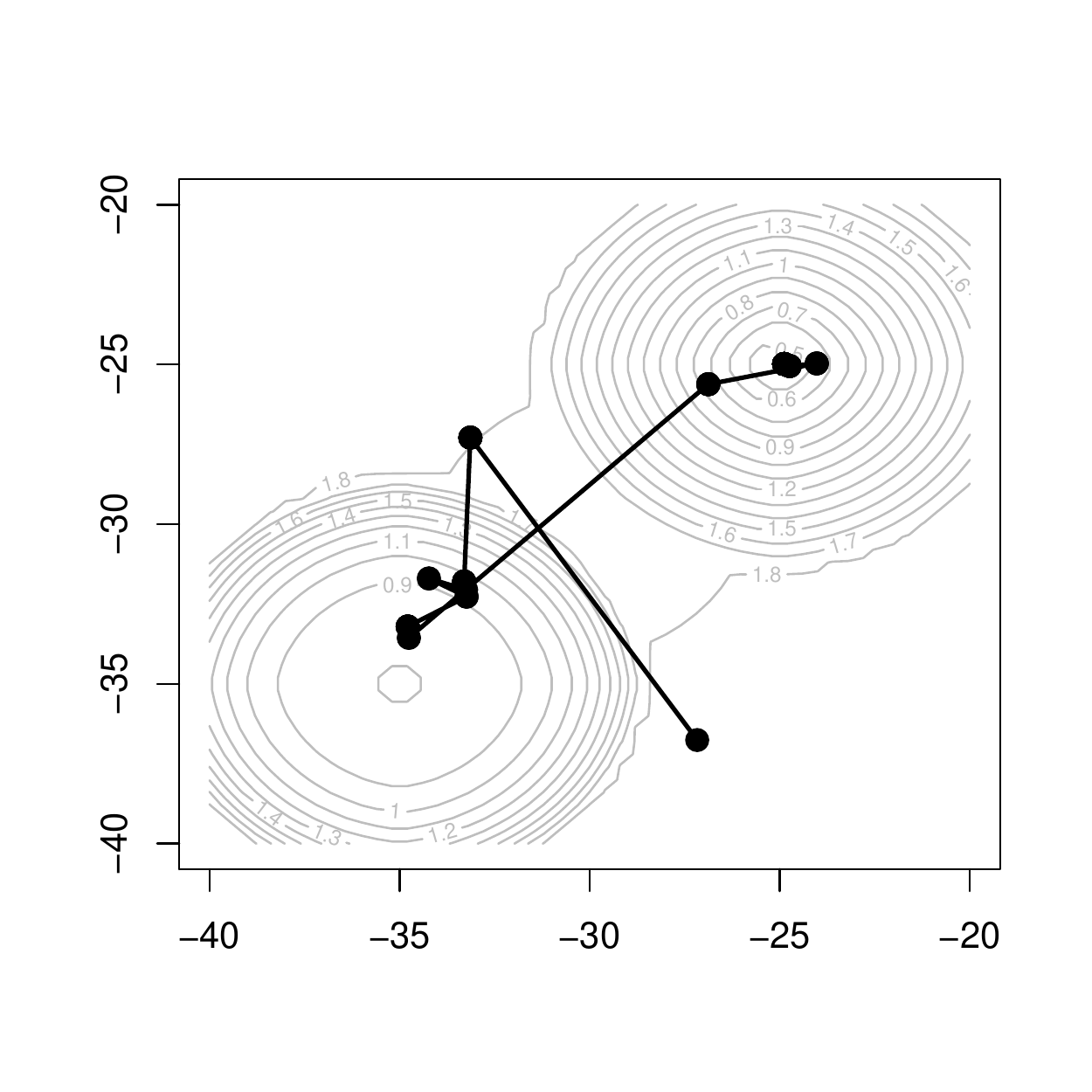}	
		\vspace{-11mm} 
  	\caption{rPSOleh current best}
		\label{fig:PSOlehSearch}
	\end{subfigure}	
		
	\vspace{-5mm} 
			
	\begin{subfigure}{.38\textwidth}
		\centering
		\includegraphics[width=2.5in, height=2.6in]{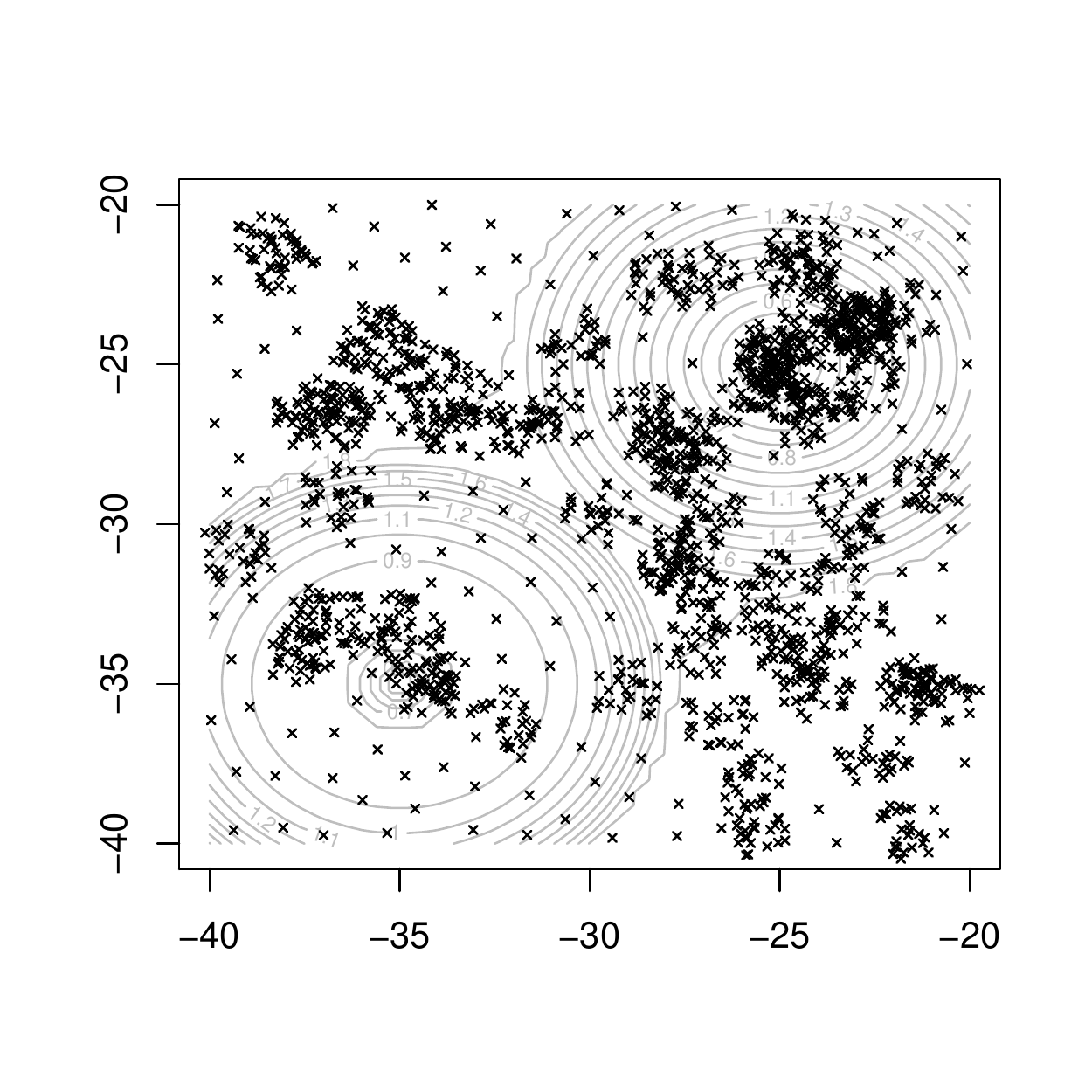}		
		\vspace{-11mm} 
	  \caption{rPSOlehdd points}	
		\label{fig:PSOlehddPoints}
	\end{subfigure}%
	\hspace{7mm} 
	\begin{subfigure}{.38\textwidth}
		\centering
		\includegraphics[width=2.5in, height=2.6in]{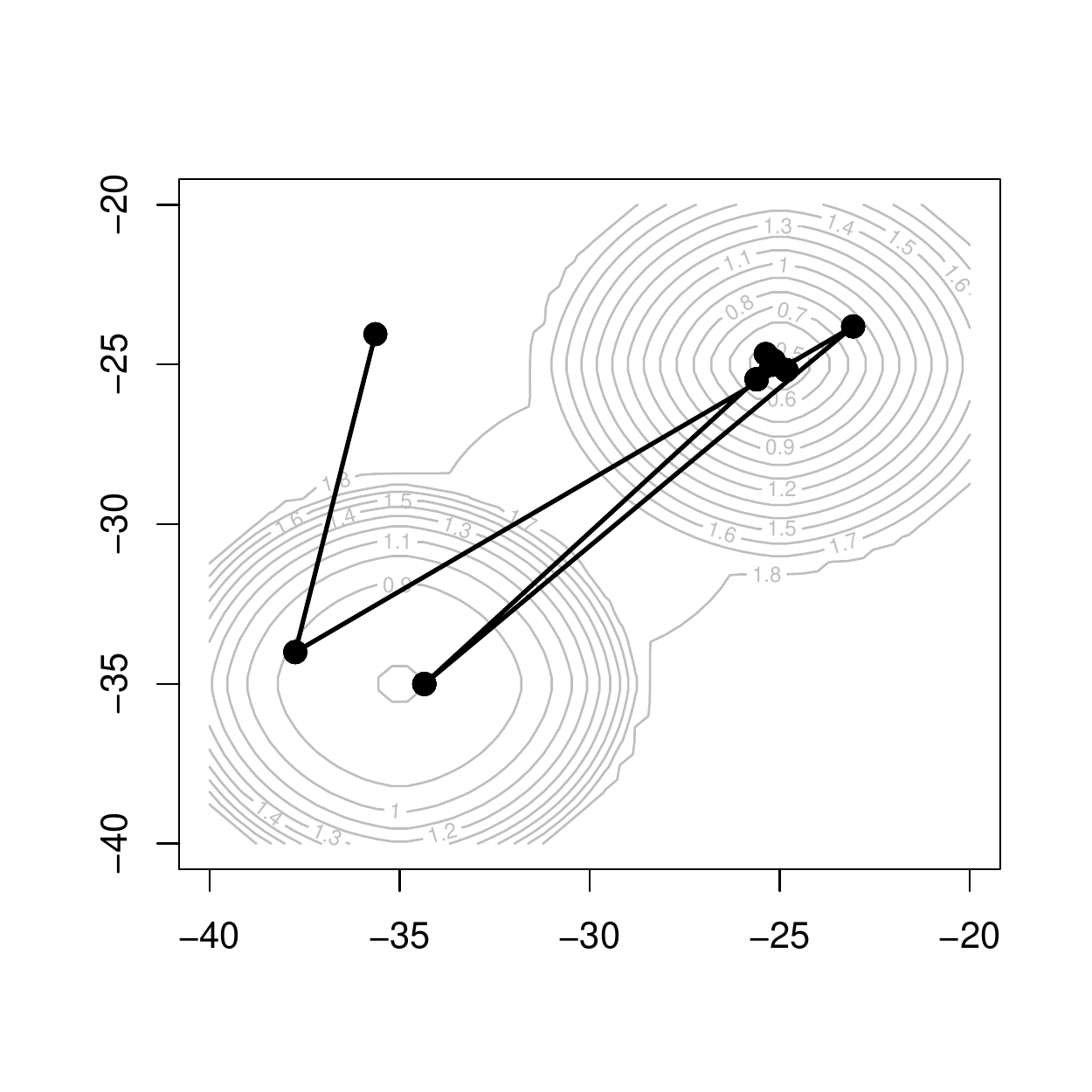}		
		\vspace{-11mm} 
  	\caption{rPSOlehdd current best}
		\label{fig:PSOlehddSearch}
	\end{subfigure}	
	
	\caption{Example searches of the 2D version of the Pickelhaube problem Figures~\ref{fig:PickelhaubeNom} and ~\ref{fig:PickelhaubeWorst} with $\Gamma$=1, for the new enhanced rPSO metaheuristics. Evaluated points are shown on the left, with the path of the improving current robust best on the right. The outer heuristics are: rPSOdd, rPSOleh, and rPSOlehdd. The nominal global optimum is at $(-35, -35)$, and the robust global optimum is at $\sim(-25, -25)$.}
	\label{fig:EnhancedPSO Robust2DSearches}
\end{figure}

Contour plots of an example search using the rPSOleh heuristic applied to the 2D Pickelhaube problem are shown in Figures~\ref{fig:PSOlehPoints} and~\ref{fig:PSOlehSearch}. The nature of the rPSOleh plots is primarily a combination of the example rPSO and LEH searches seen in Figure~\ref{fig:BaselineRobust2DSearches}. There are some inner maximisation groupings of evaluated points for each particle, as the particles iterate, plus some limiting of the extent of these inner searches -- moving on without undertaking a full $\Gamma$-uncertainty neighbourhood search, plus the extensive LEH exploration of the solution space as the search moves to regions devoid of hcps.

\subsection{rPSOlehdd}
\label{sec:rPSOlehdd}

Our final enhanced rPSO heuristic is in fact just the combination of the two new heuristics rPSOdd and rPSOleh: we augment the baseline rPSO Algorithm~\ref{BaselineRPSO} with Algorithms~\ref{ddVelComponent} and~\ref{partLEHAlgorithm}. In combination the new algorithms require input information for $\pmb{x}_c$, $H$, $F_H$, full particle information for particle $j$($i$), $\gb$ and $dormancy_{count}^j$. Parameter values for $\sigma$, $\sigma_{limit}$, $\sigma_{no}$, $C_3$, $ToBeEvaluated$, $dormancy_{limit}$ and $placement_{limit}$ are also required.

In the experimental testing the $\sigma$, $\sigma_{limit}$, $C_3$ and $dormancy_{limit}$ parameters are tuned, see Section~\ref{sec:parameterTuning}, in addition to the tuned baseline rPSO parameters (Section~\ref{sec:bruteForcePSO}) and the parameters controlling the GA applied to solving problem (lehMM) (Section~\ref{sec:rPSOleh}).

Algorithm~\ref{BaselineRPSO} still requires a single change to accommodate the d.d.\ component: in line~\ref{UpdateVelocity} the enhanced velocity equation (Vel2) replaces equation (Vel1). To accommodate the LEH components Algorithm~\ref{BaselineRPSO} requires several minor updates. The call to the inner maximisation Algorithm~\ref{InnerAlgorithm} must be updated in accordance with (call~\ref{InnerAlgorithm}). Plus the updates to include (relocateLEH) and (addCounter), and the addition of input parameters $dormancy_{limit}$ and $placement_{limit}$ to Algorithm~\ref{BaselineRPSO}, as described in Section~\ref{sec:rPSOleh}. 

These updates and calls to Algorithms~\ref{InnerAlgorithm},~\ref{ddVelComponent} and~\ref{partLEHAlgorithm} from the baseline rPSO Algorithm~\ref{BaselineRPSO} complete the rPSOlehdd heuristic.

Contour plots of an example search using the rPSOlehdd heuristic applied to the 2D Pickelhaube problem are shown in Figures~\ref{fig:PSOlehddPoints} and~\ref{fig:PSOlehddSearch}. The nature of these plots is somewhat similar to the rPSOleh plots also shown in Figure~\ref{fig:EnhancedPSO Robust2DSearches}, combining PSO and LEH search elements. Also, however, as with the rPSOdd heuristic the additional particle level d.d.\ component in the velocity function appears to introduce some elements of robust local search. As with rPSOdd the extent of any such d.d.\ component will be influenced by the heuristic parameter value settings.


\section{Computational experiments}
\label{sec:experimentsResults}

\subsection{Experimental set up}
\label{sec:setUp}

We discussed three global robust metaheuristics as a baseline against which to test our three new approaches, based on enhancements of a baseline rPSO approach. This gives six heuristics:

\begin{enumerate}[leftmargin=*]		
	\item Baseline: largest empty hypersphere (LEH) \cite{HughesGoerigkWright2019}.
	\item Baseline: repeating descent direction (d.d.) based on \cite{BertsimasNohadaniTeo2010}. 
	\item Baseline: robust particle swarm optimisation (rPSO) based on `standard' PSO formulation \cite{KennedyEberhart1995, KennedyEberhartShi2001, Ghazali2009}.
	\item Enhanced rPSO: original PSO velocity formulation with the inclusion of a d.d.\ component (rPSOdd).
	\item Enhanced rPSO: using historic information, the stopping condition and LEH (rPSOleh).
	\item Enhanced rPSO: using historic information, the stopping condition, the inclusion of a d.d.\ component and LEH (rPSOlehdd). 
\end{enumerate}	
Recall that all inner maximisation analysis is undertaken exclusively using uniform random sampling in the $\Gamma$-radius hypersphere that forms the uncertainty neighbourhood of any given point, see Algorithm~\ref{InnerAlgorithm}. The level of $\Gamma$-radius sampling is a tuned value for all heuristics, and is also a maximum level of sampling when a stopping condition is employed by the heuristic.

Each run of each heuristic identifies an estimate of the location and value of a robust global optimum. The robust global value used by the heuristic is likely to be an inaccurate estimate of the actual worst case value at the identified location. Therefore we post process all estimated values by randomly sampling 1,000,000 points in the $\Gamma$-uncertainty neighbourhood of the identified robust location, and taking the maximum sampled value as the estimated worst case cost. This robust value is taken as the output of a single heuristic-test problem run. This post processing does not impact on the heuristic search. 

We employ ten multi-dimensional test functions over six dimensions: 2D, 5D, 10D, 30D, 60D and 100D, plus the 2D polynomial test problem from \cite{BertsimasNohadaniTeo2010}. This gives 61 test problem instances for each of the six comparator robust metaheuristics to be applied to. Each heuristic is applied to each test instance 200 times to give reasonable sample sets of results for comparison, in order to identify the best performing heuristics. In line with assumed restrictions on numbers of function evaluations when handling real-world problems, each test problem run is limited to a budget of 5,000 function evaluations. Prior to undertaking the sample runs parameter tuning has been applied, as described in Section~\ref{sec:parameterTuning}.

\begin{figure}[htbp]
	\centering

	\vspace{-6mm} 

	\begin{subfigure}{.38\textwidth}
		\centering
		\includegraphics[width=2.5in, height=2.6in]{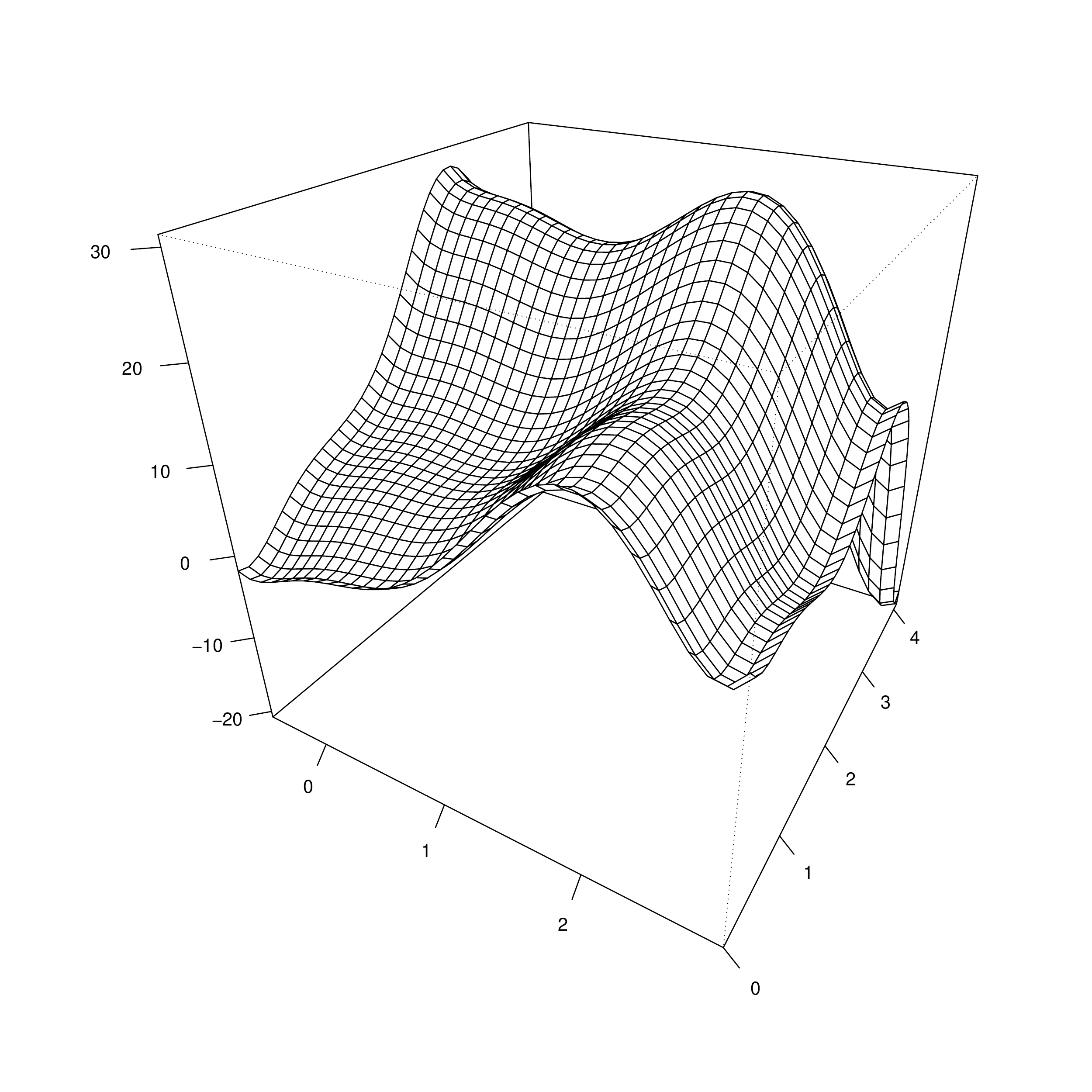}		
		\vspace{-11mm} 
	  \caption{2D Polynomial Nominal}	
		\label{fig:Bertsimas2DNom}
	\end{subfigure}%
	\hspace{7mm} 
	\begin{subfigure}{.38\textwidth}
		\centering
		\includegraphics[width=2.5in, height=2.6in]{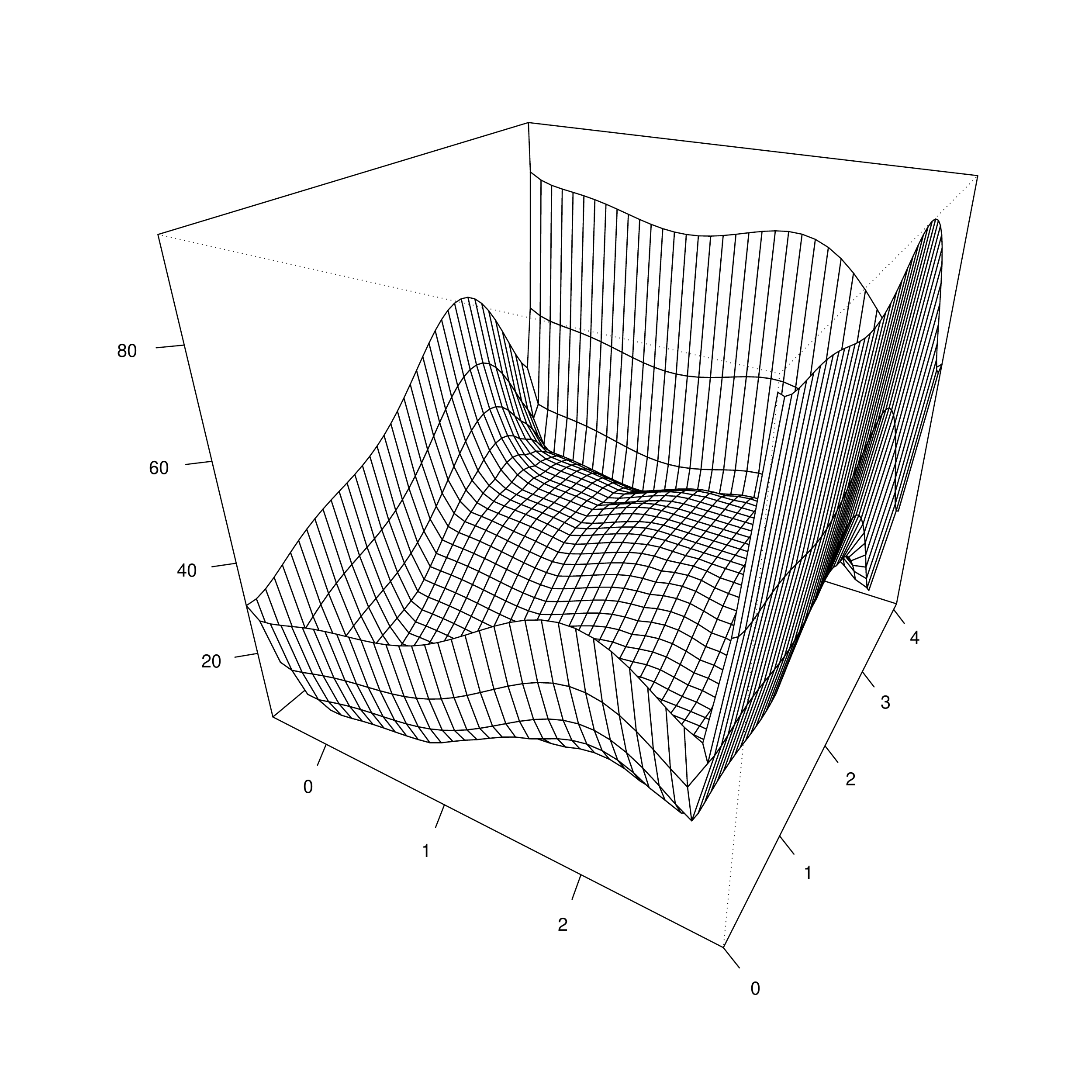}	
		\vspace{-11mm} 
  	\caption{2D Polynomial Worst}
		\label{fig:Bertsimas2DWorst}
	\end{subfigure}			

	\caption{Plots of 2D Polynomial test function \cite{BertsimasNohadaniTeo2010} in the rPSO test suite.}
	\label{fig:poly2D}
\end{figure}

\begin{figure}[htbp]
	\centering
		
	\begin{subfigure}[t]{.24\textwidth}
		\includegraphics[width=\textwidth]{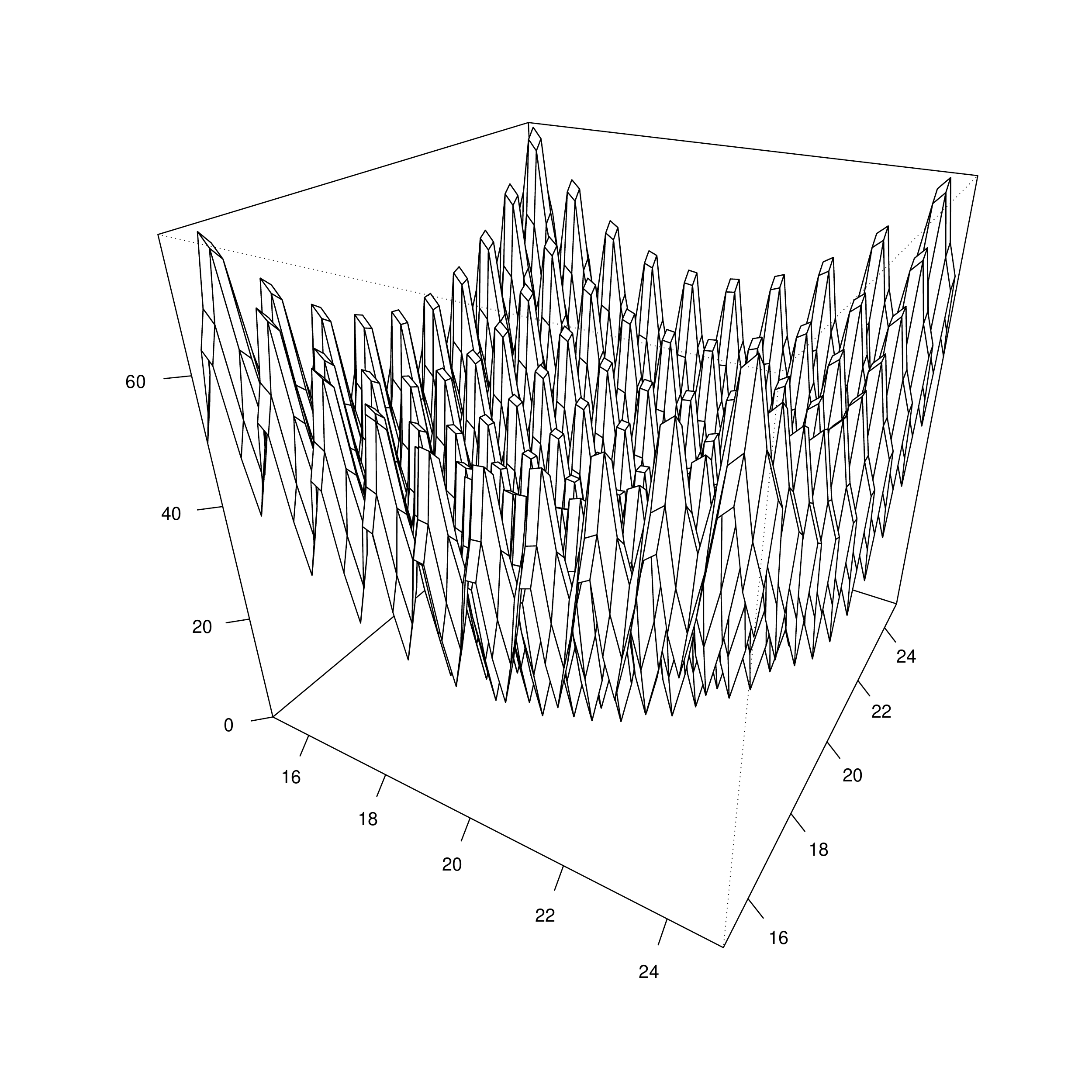}		
  	\caption{Rastrigin Nominal} \label{fig:RastriginNom}
	\end{subfigure}%
	\begin{subfigure}[t]{.24\textwidth}
		\includegraphics[width=\textwidth]{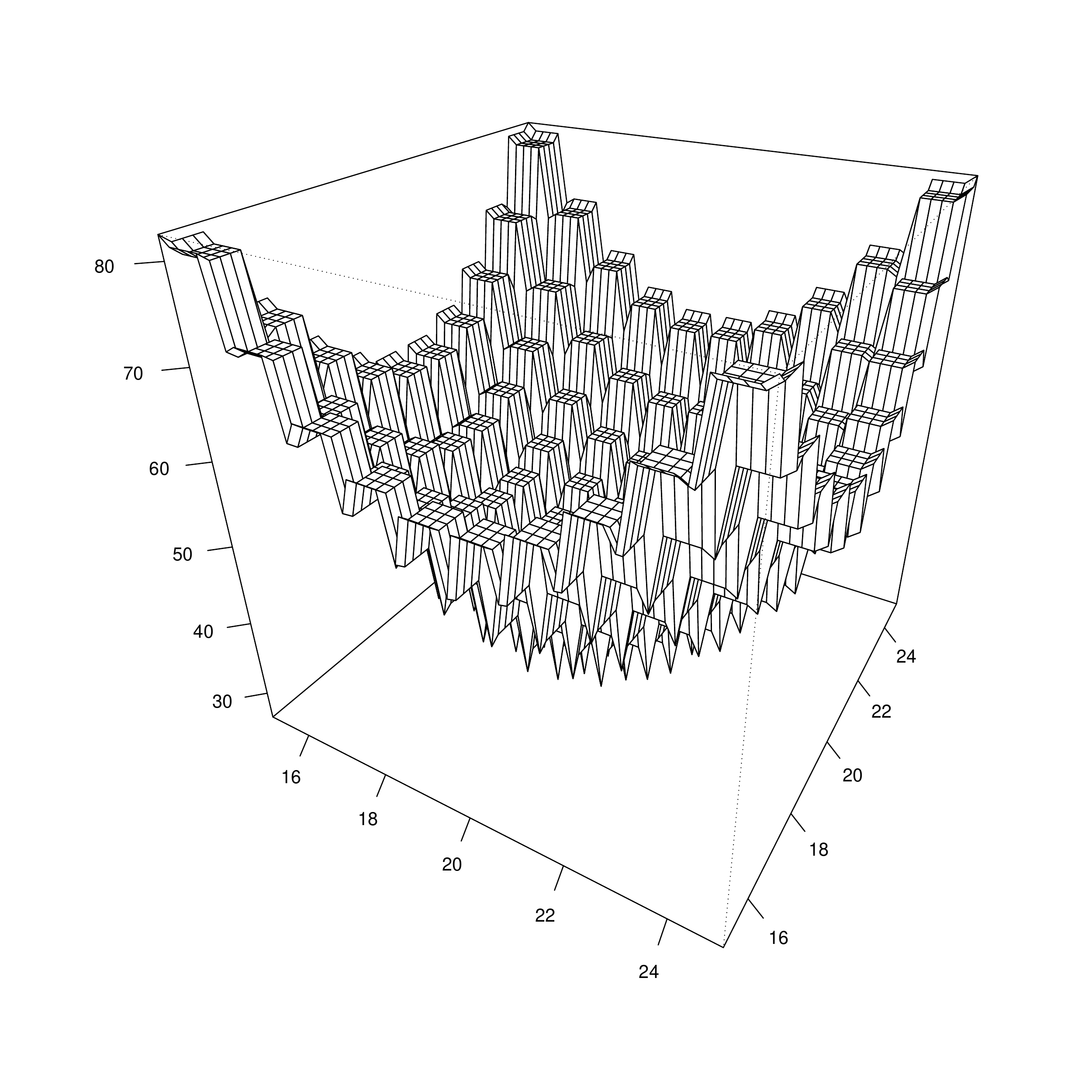}	
  	\caption{Rastrigin Worst} \label{fig:RastriginWorst}
	\end{subfigure}		
	\begin{subfigure}[t]{.24\textwidth}
		\includegraphics[width=\textwidth]{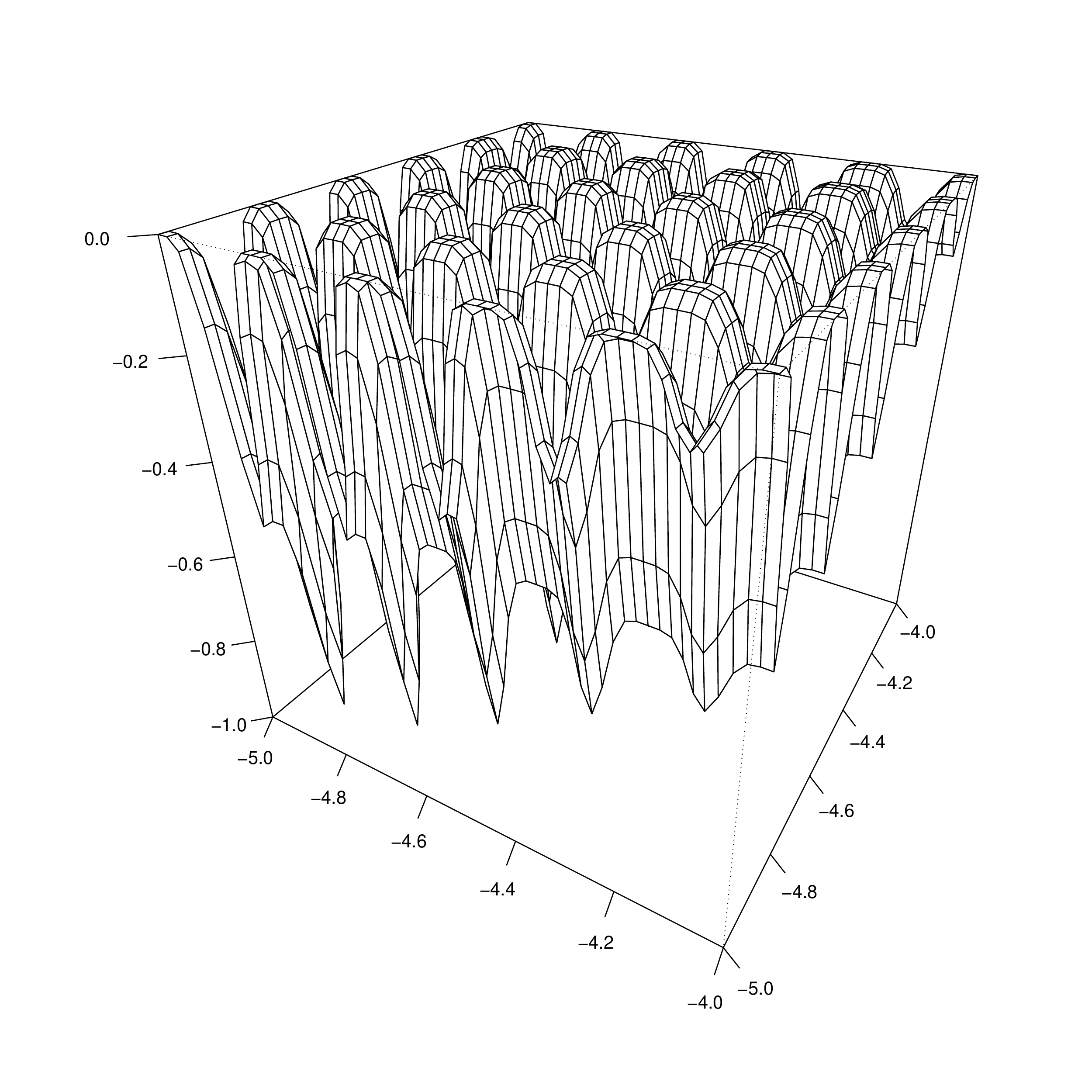}		
  	\caption{Multipeak F1 Nominal} \label{fig:MultipeakF1Nom}
	\end{subfigure}%
	\begin{subfigure}[t]{.24\textwidth}
		\includegraphics[width=\textwidth]{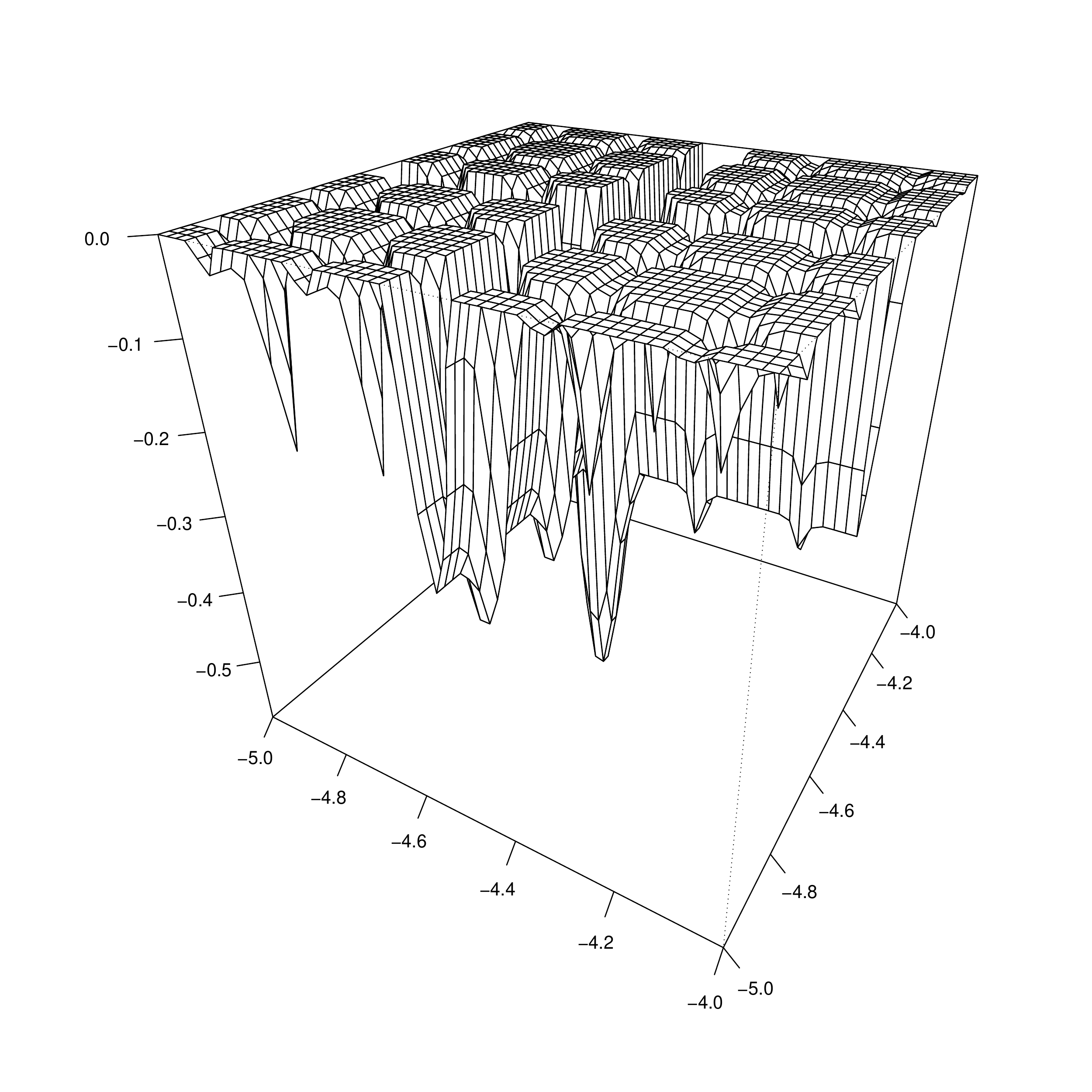}	
  	\caption{Multipeak F1 Worst} \label{fig:MultipeakF1Worst}
	\end{subfigure}
	\begin{subfigure}[t]{.24\textwidth}
		\includegraphics[width=\textwidth]{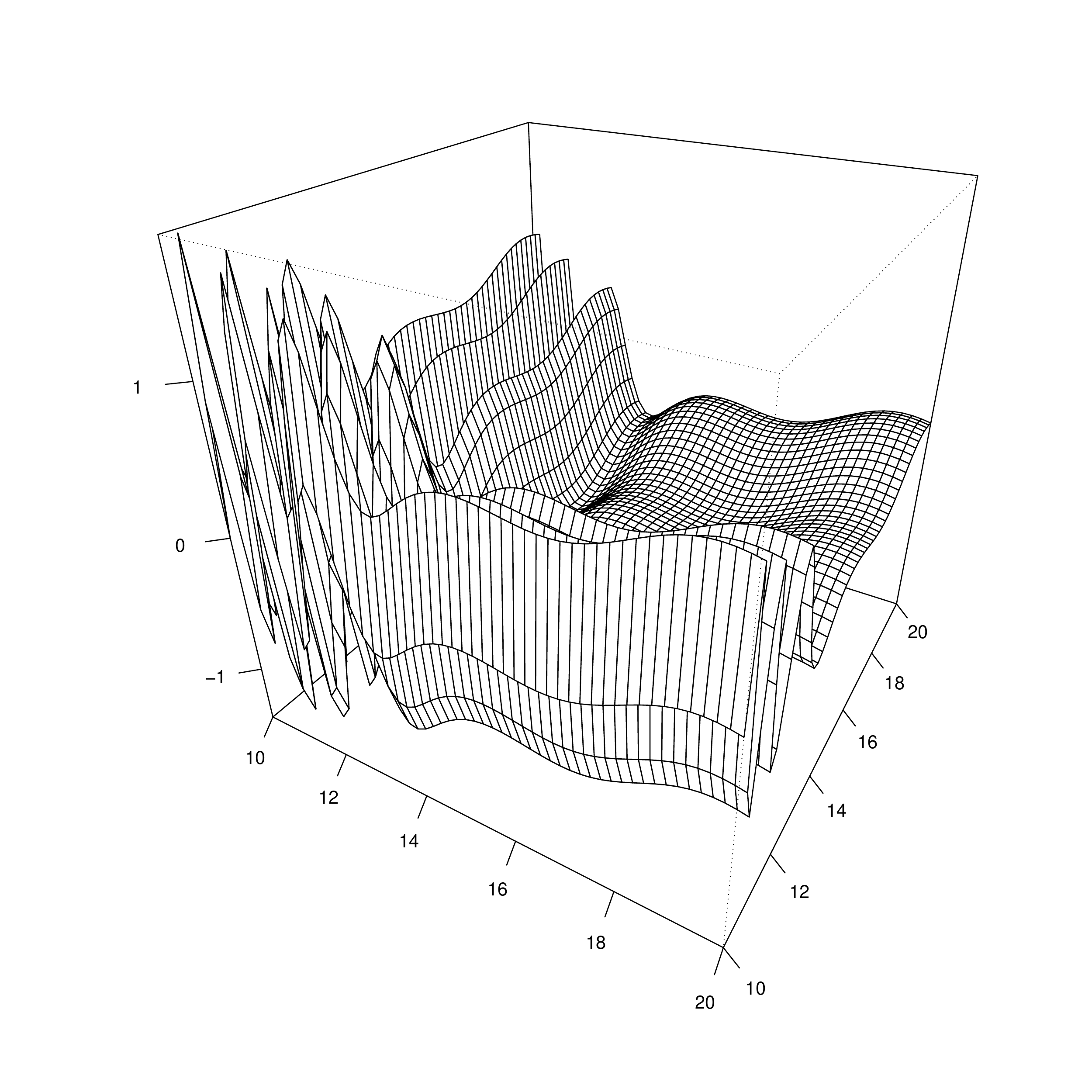}	
  	\caption{Multipeak F2 Nominal} \label{fig:MultipeakF2Nom}
	\end{subfigure}%
	\begin{subfigure}[t]{.24\textwidth}
		\includegraphics[width=\textwidth]{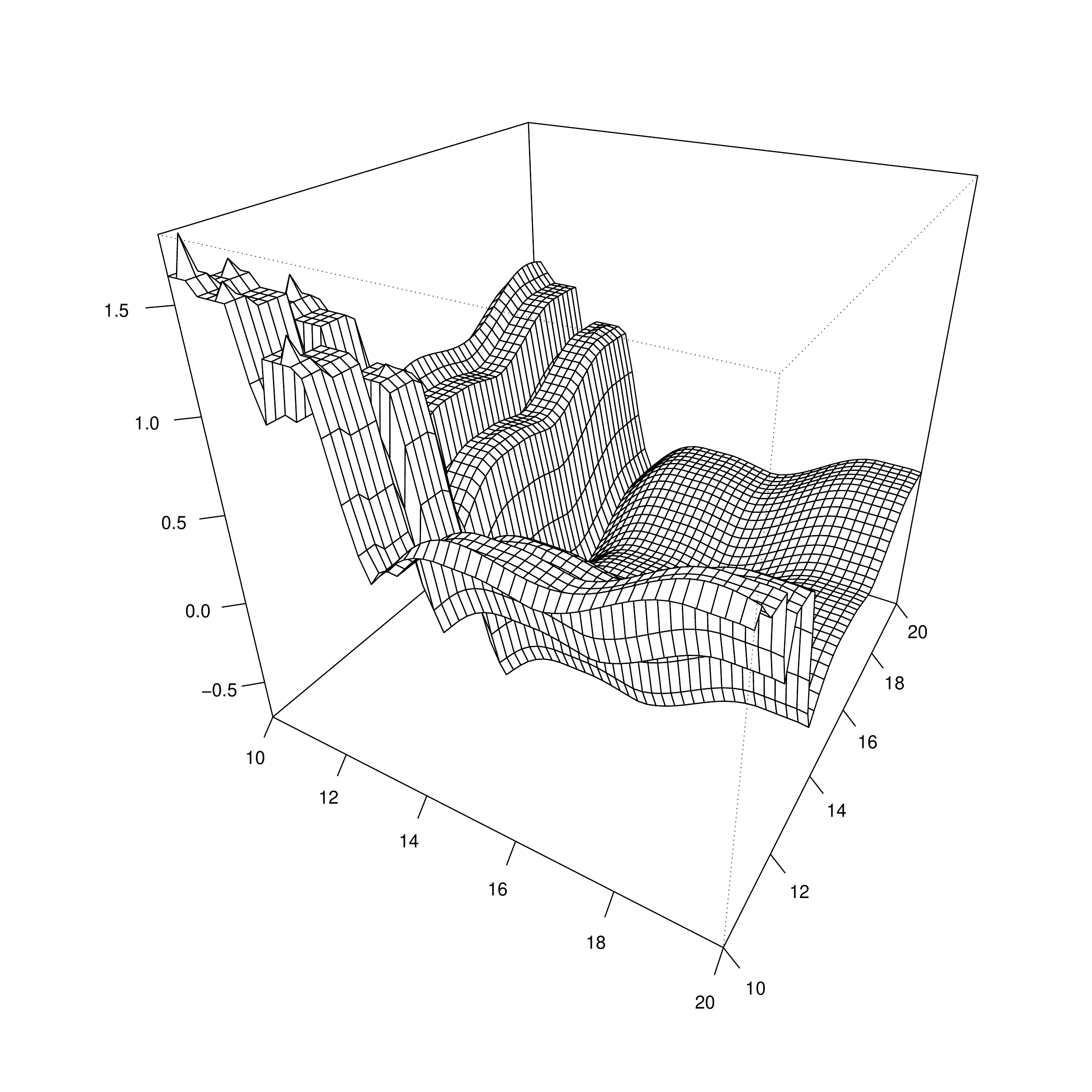}	
  	\caption{Multipeak F2 Worst} \label{fig:MultipeakF2Worst}
	\end{subfigure}
        \begin{subfigure}[t]{.24\textwidth}
		\includegraphics[width=\textwidth]{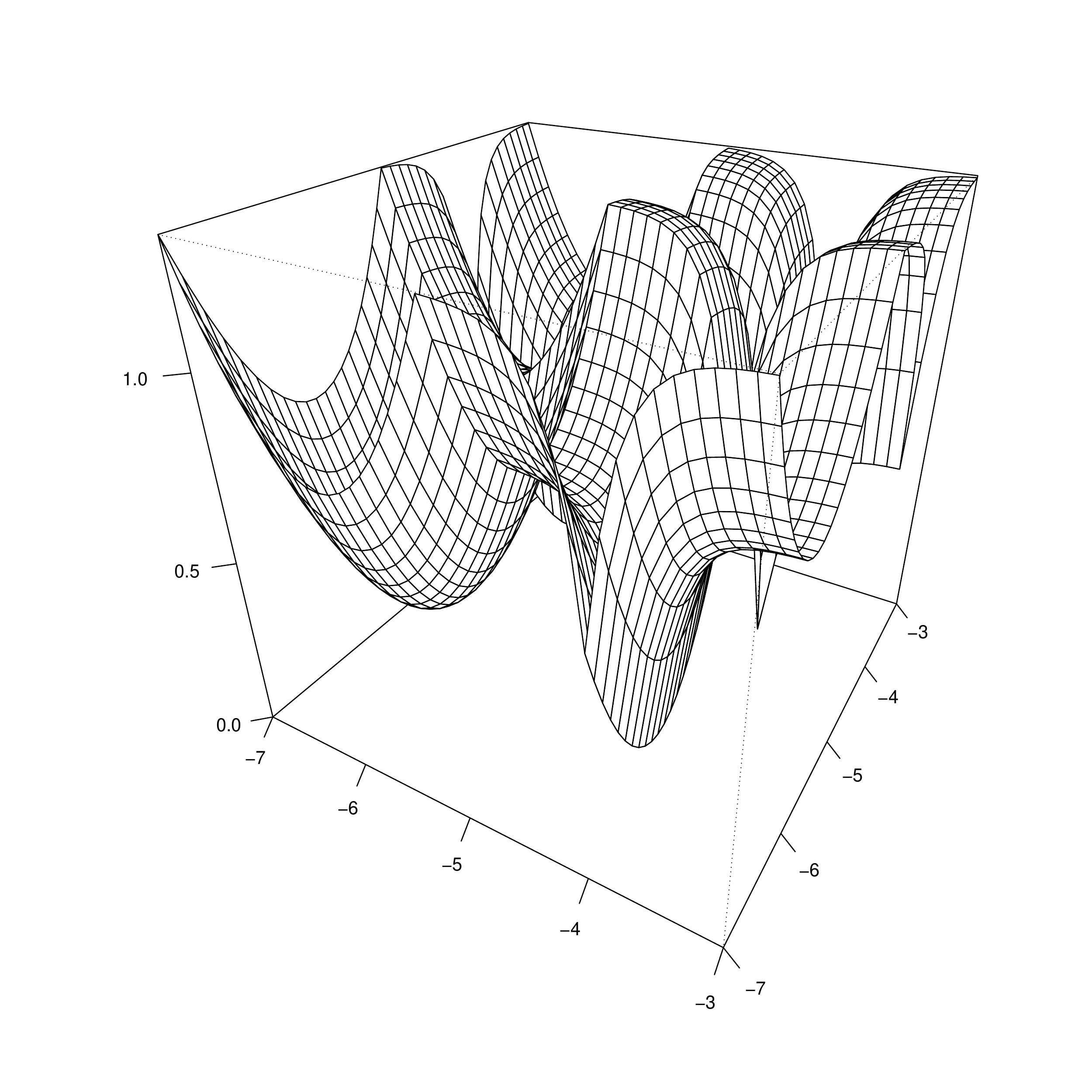}
  	\caption{Branke's Multipeak Nominal} \label{fig:BrankesMultipeakNom}
	\end{subfigure}%
	\begin{subfigure}[t]{.24\textwidth}
		\includegraphics[width=\textwidth]{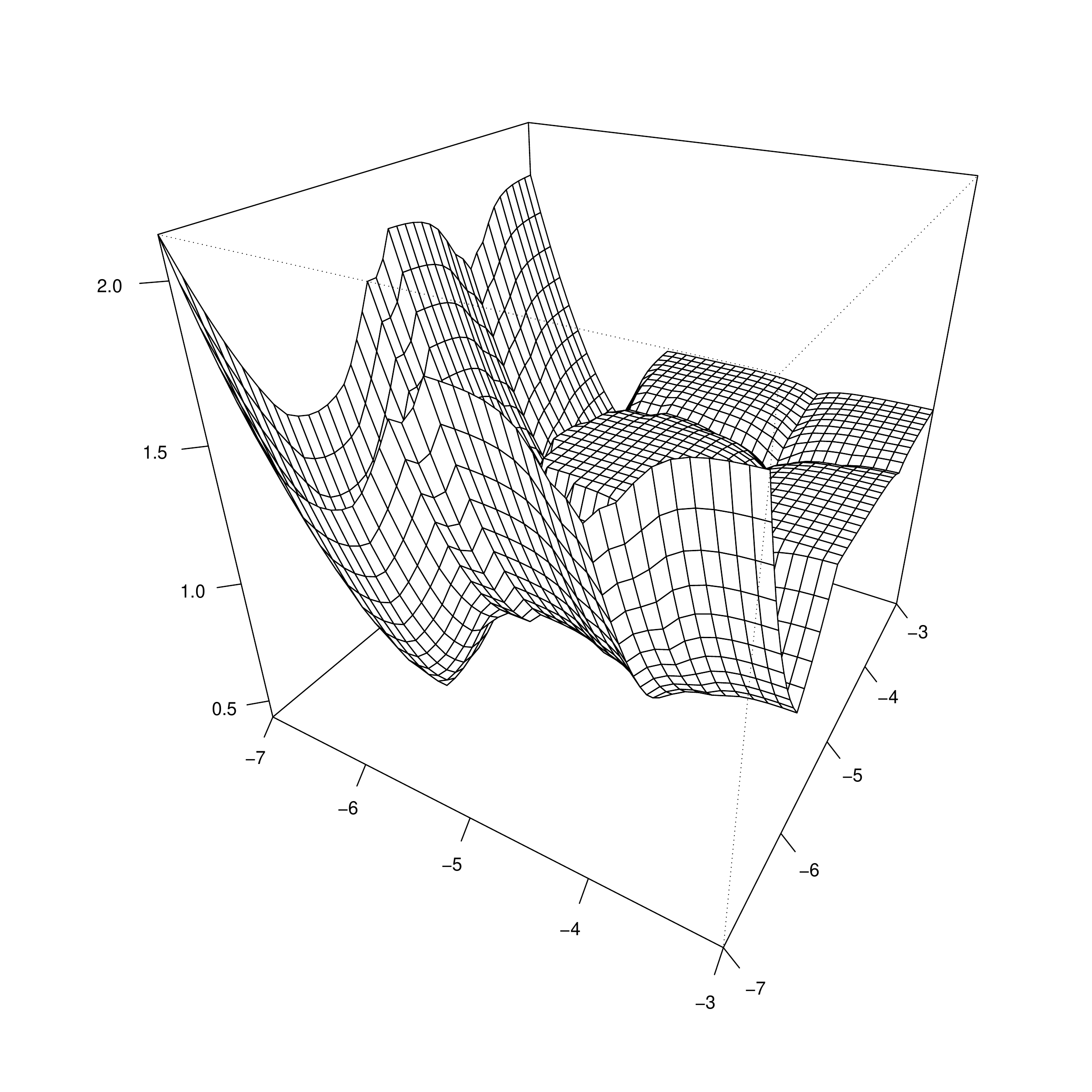}	
  	\caption{Branke's Multipeak Worst} \label{fig:BrankesMultipeakWorst}
	\end{subfigure}
	\begin{subfigure}[t]{.24\textwidth}
		\includegraphics[width=\textwidth]{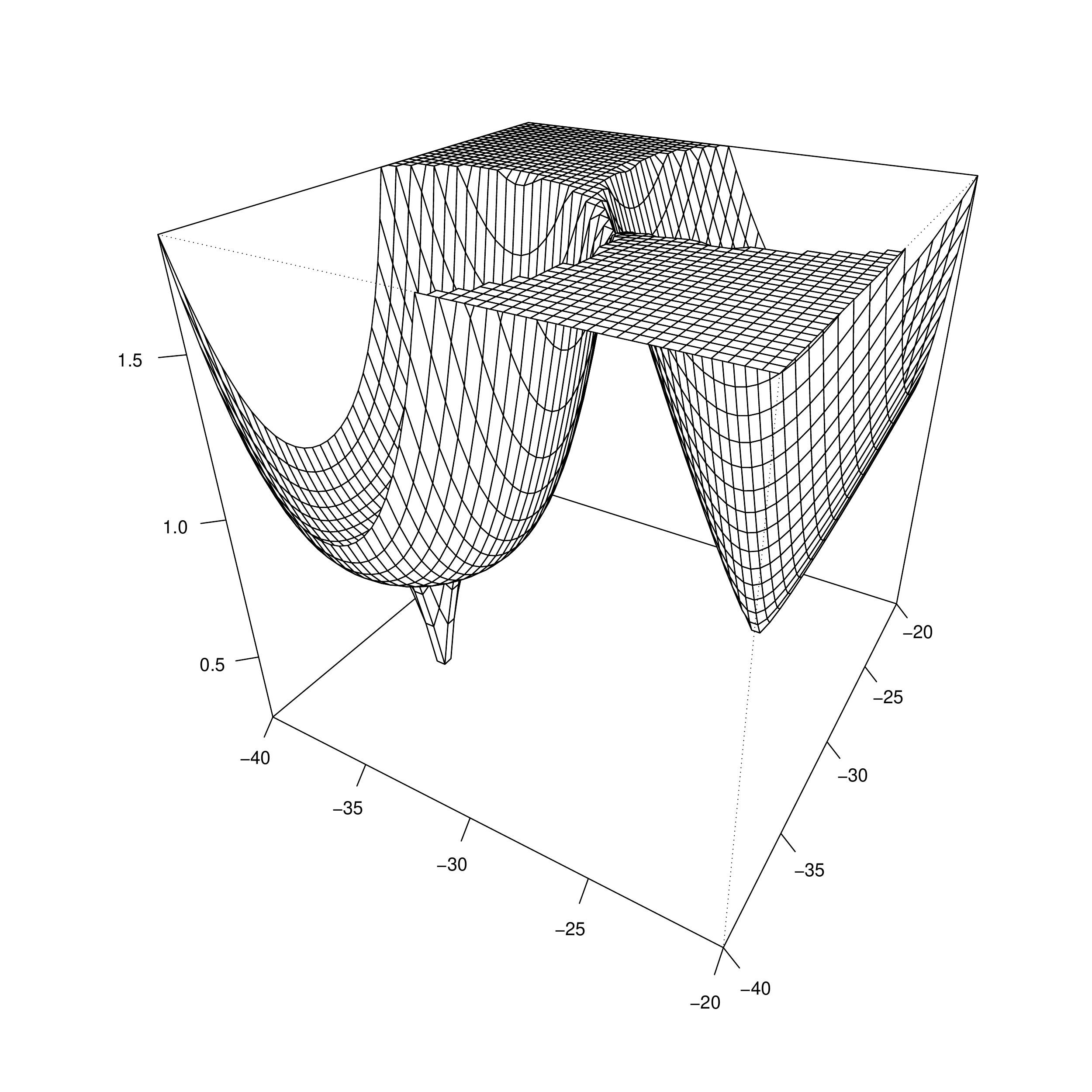}		
	  \caption{Pickelhaube Nominal}	\label{fig:PickelhaubeNom}
	\end{subfigure}%
	\begin{subfigure}[t]{.24\textwidth}
		\includegraphics[width=\textwidth]{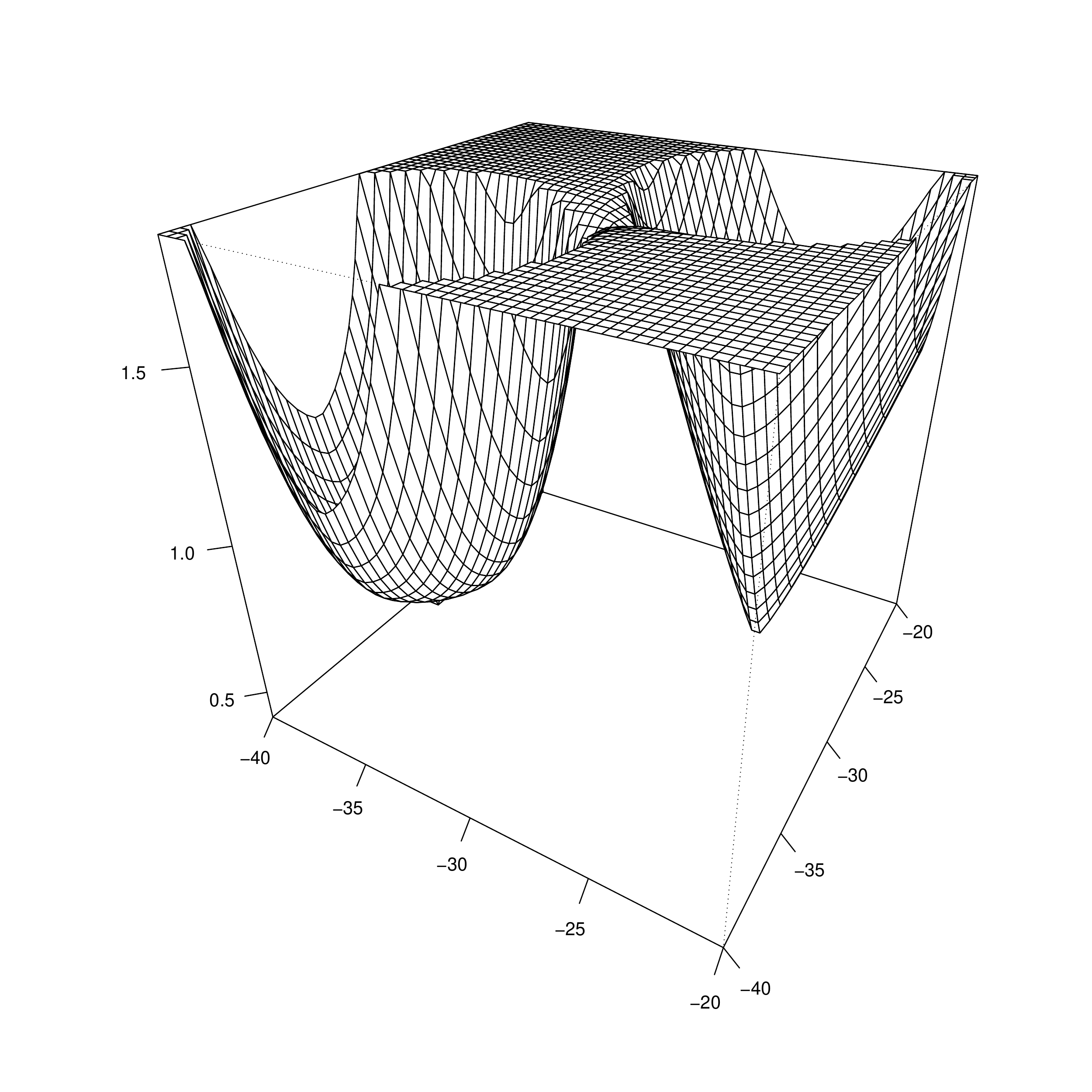}	
  	\caption{Pickelhaube Worst}\label{fig:PickelhaubeWorst}
	\end{subfigure}	
	\begin{subfigure}[t]{.24\textwidth}
		\includegraphics[width=\textwidth]{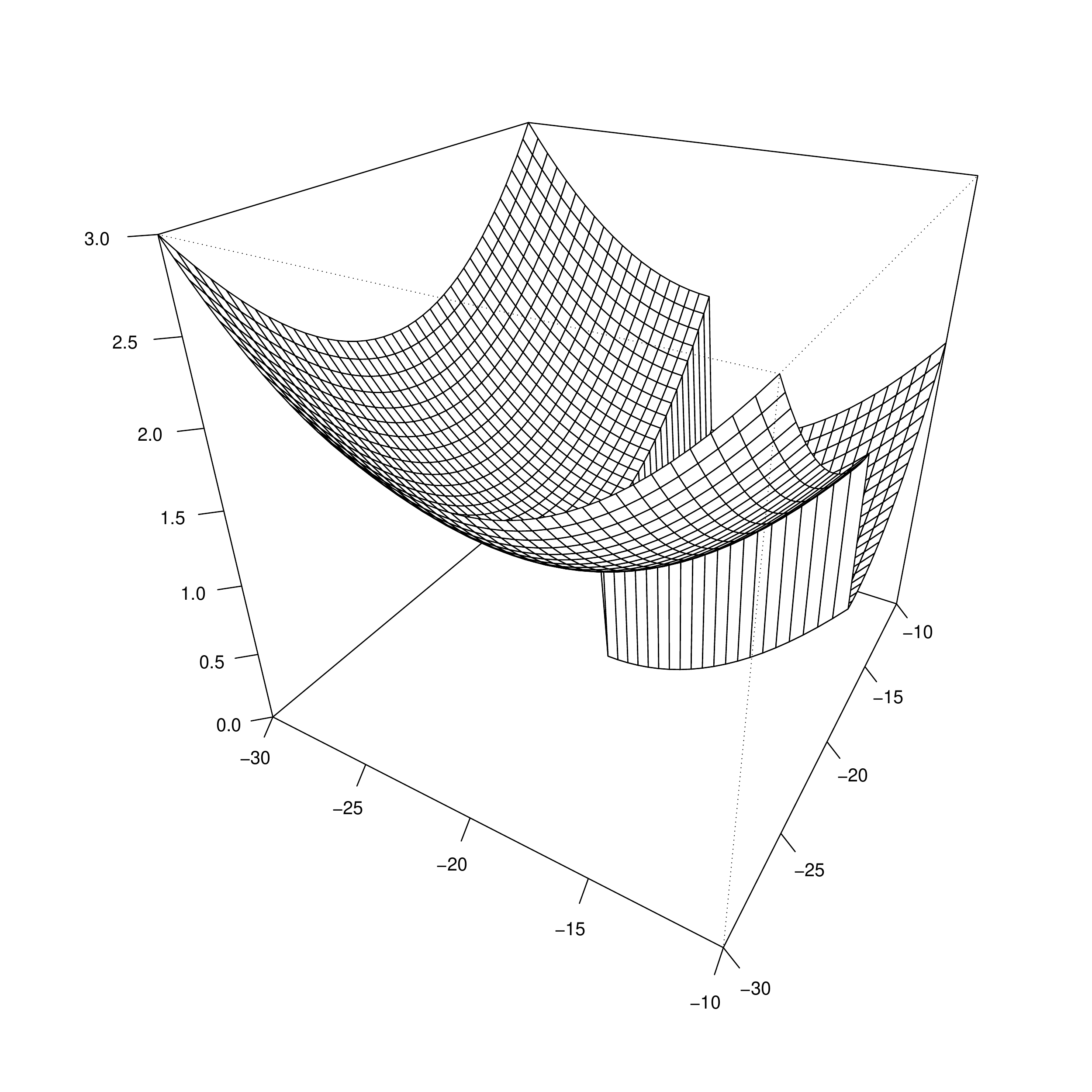}	
		\caption{Heaviside Sphere Nominal}\label{fig:HeavisideSphereNom}
	\end{subfigure}%
	\begin{subfigure}[t]{.24\textwidth}	
		\includegraphics[width=\textwidth]{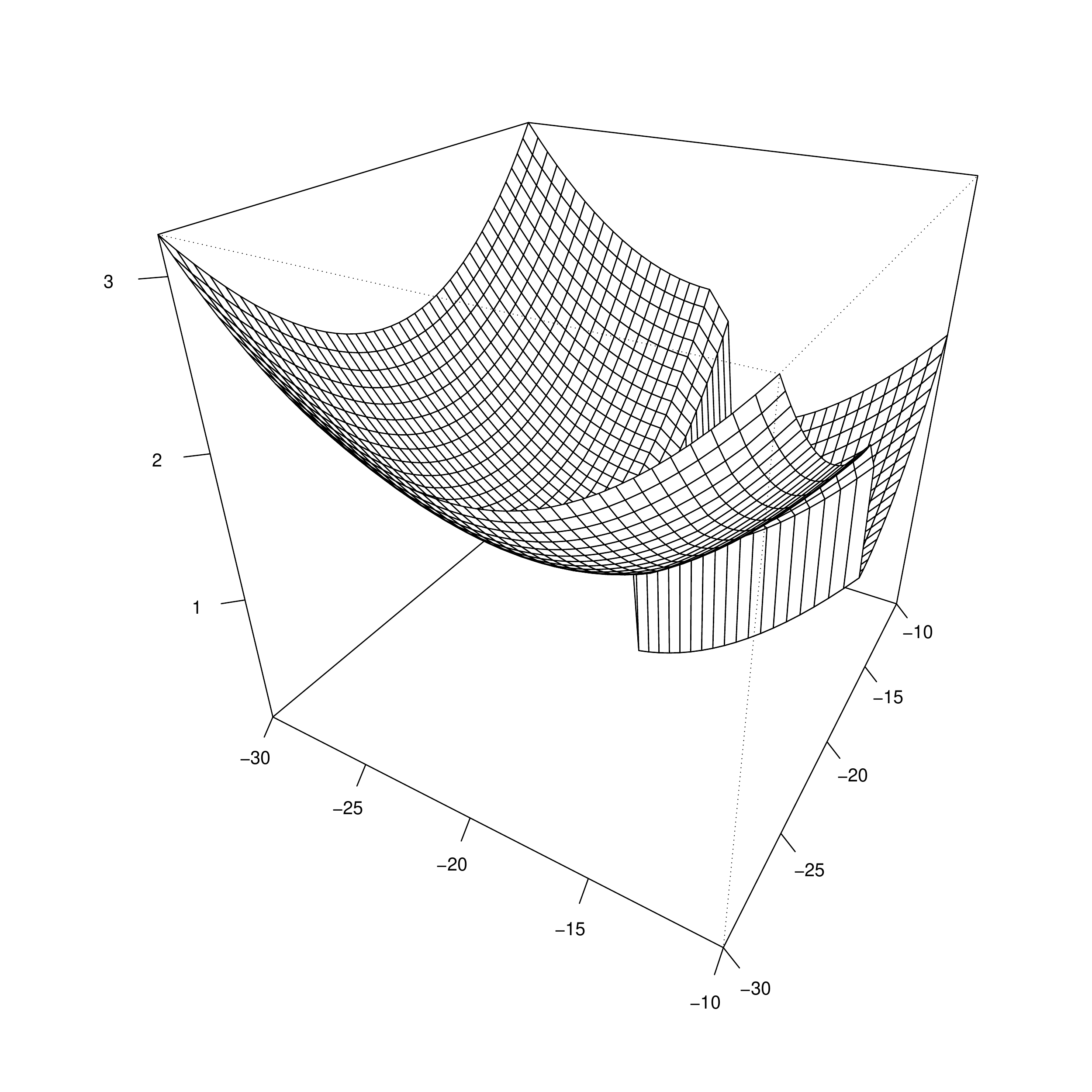}	
  	\caption{Heaviside Sphere Worst} \label{fig:HeavisideSphereWorst}
	\end{subfigure}	
        \begin{subfigure}[t]{.24\textwidth}
		\includegraphics[width=\textwidth]{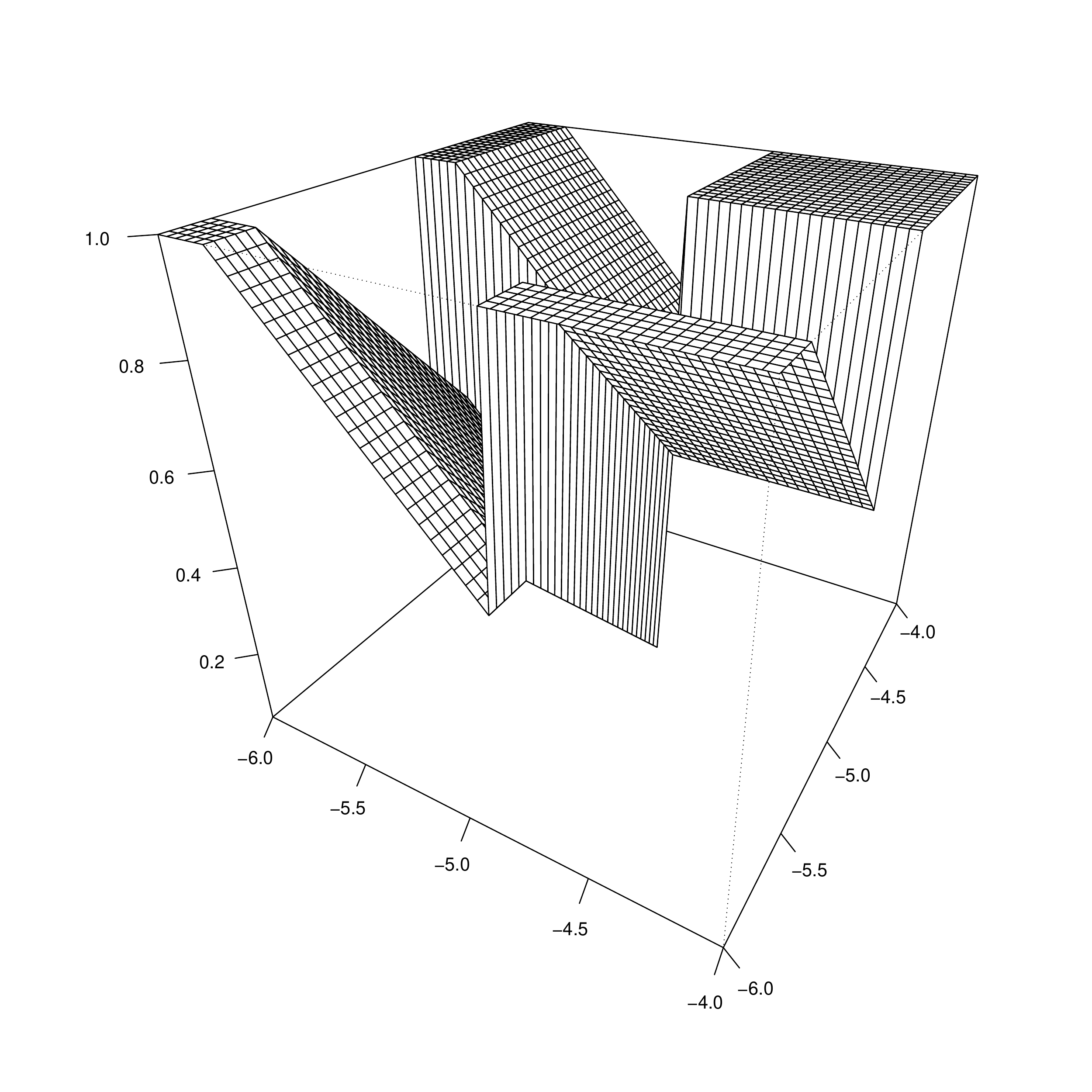}	
	  \caption{Sawtooth Nominal} \label{fig:SawtoothNom}
	\end{subfigure}%
	\begin{subfigure}[t]{.24\textwidth}
		\includegraphics[width=\textwidth]{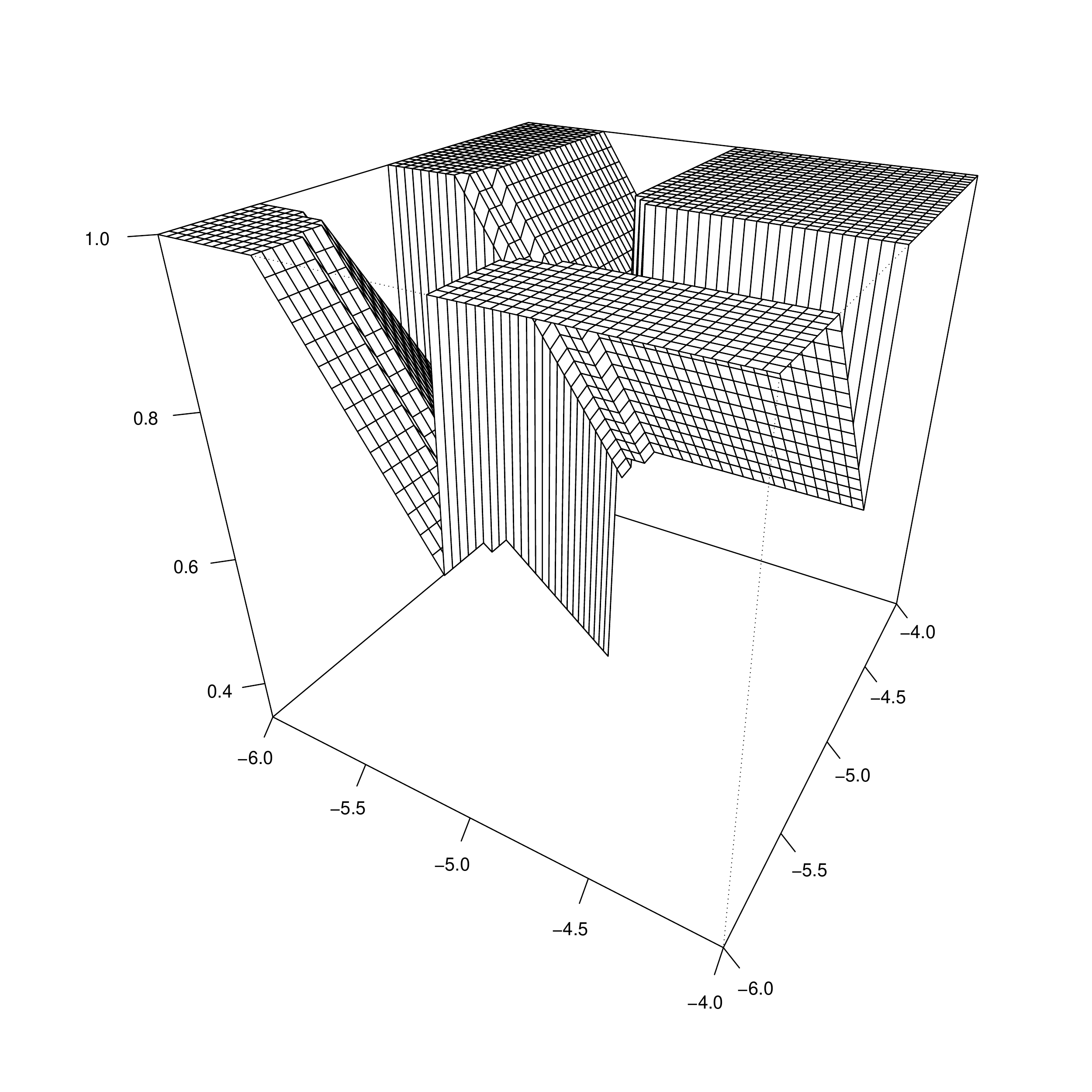}	
  	\caption{Sawtooth Worst} \label{fig:SawtoothWorst}
	\end{subfigure}	
	\begin{subfigure}[t]{.24\textwidth}
		\includegraphics[width=\textwidth]{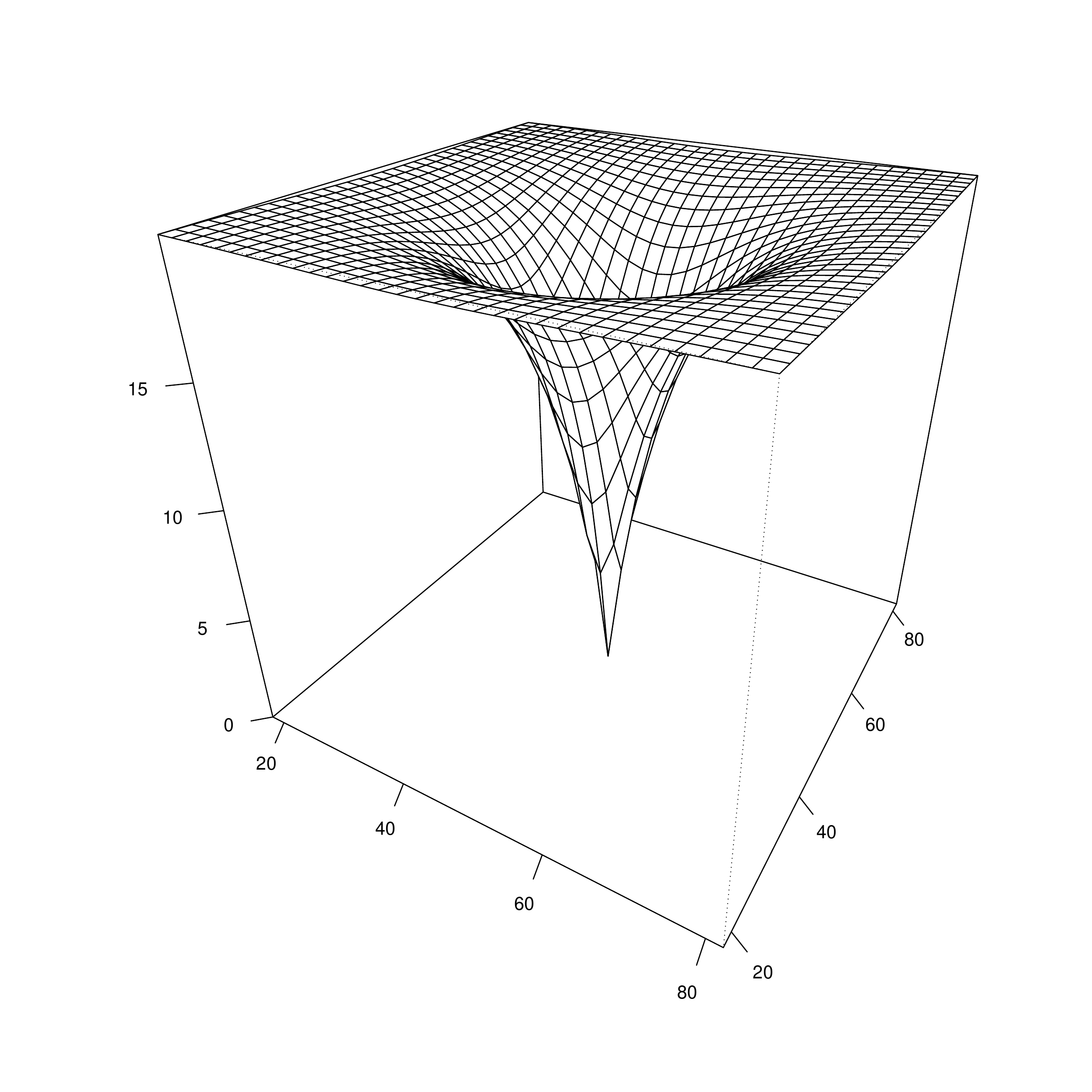}	
  	\caption{Ackley Nominal}\label{fig:AckleyNom}
	\end{subfigure}%
	\begin{subfigure}[t]{.24\textwidth}
		\includegraphics[width=\textwidth]{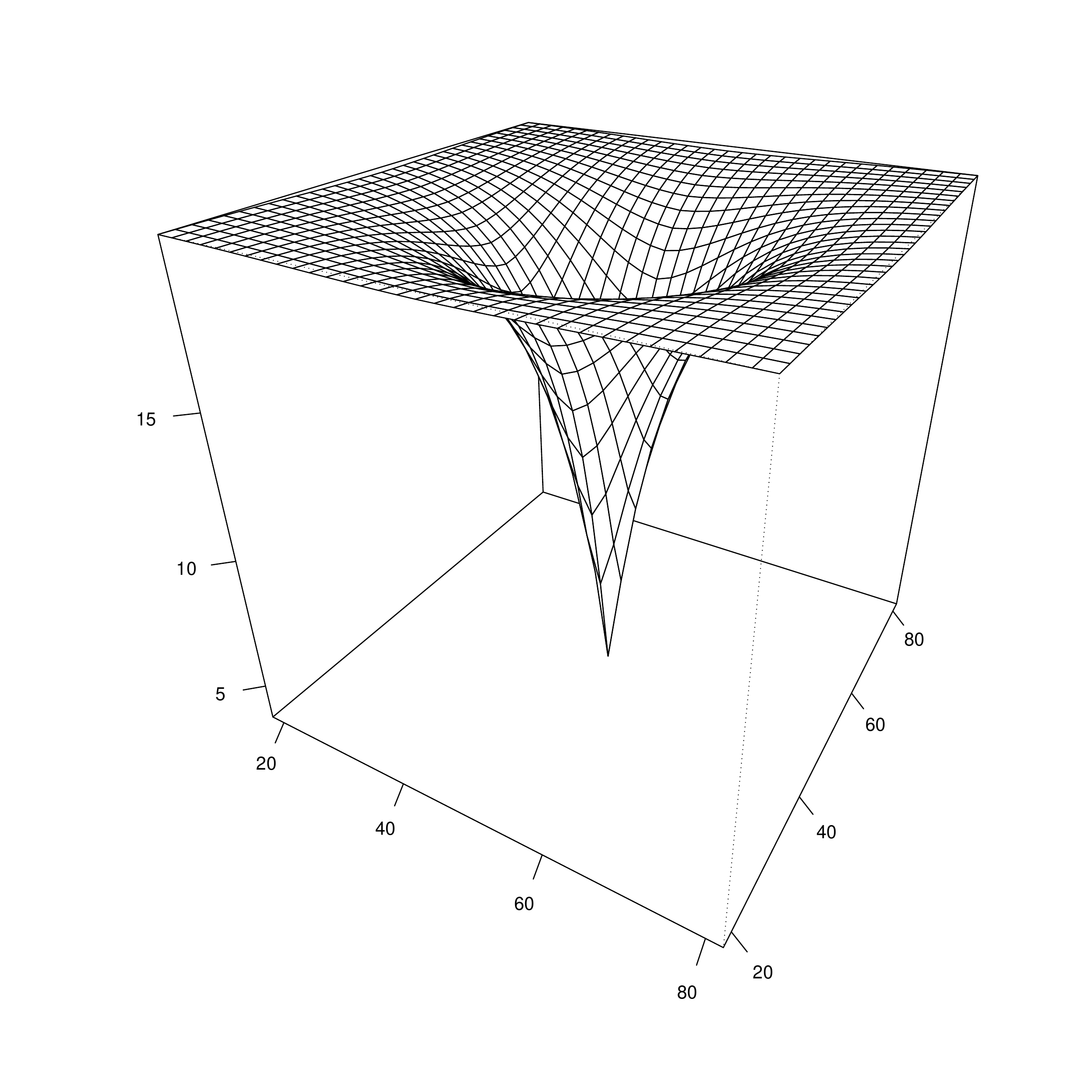}	
  	\caption{Ackley Worst} \label{fig:AckleyWorst}
	\end{subfigure}
	\begin{subfigure}[t]{.24\textwidth}
		\includegraphics[width=\textwidth]{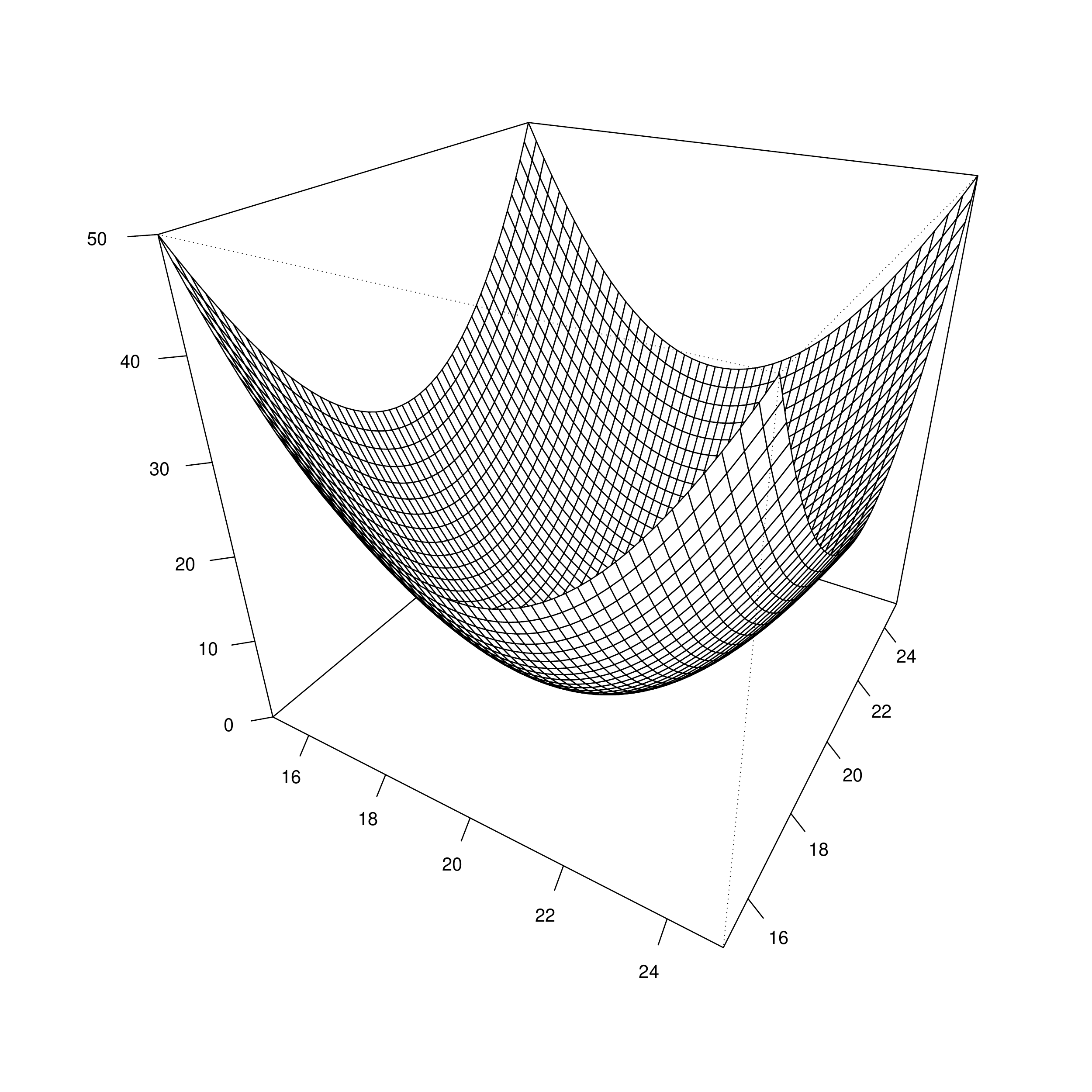}
  	\caption{Sphere Nominal} \label{fig:SphereNom}
	\end{subfigure}%
	\begin{subfigure}[t]{.24\textwidth} 
		\includegraphics[width=\textwidth]{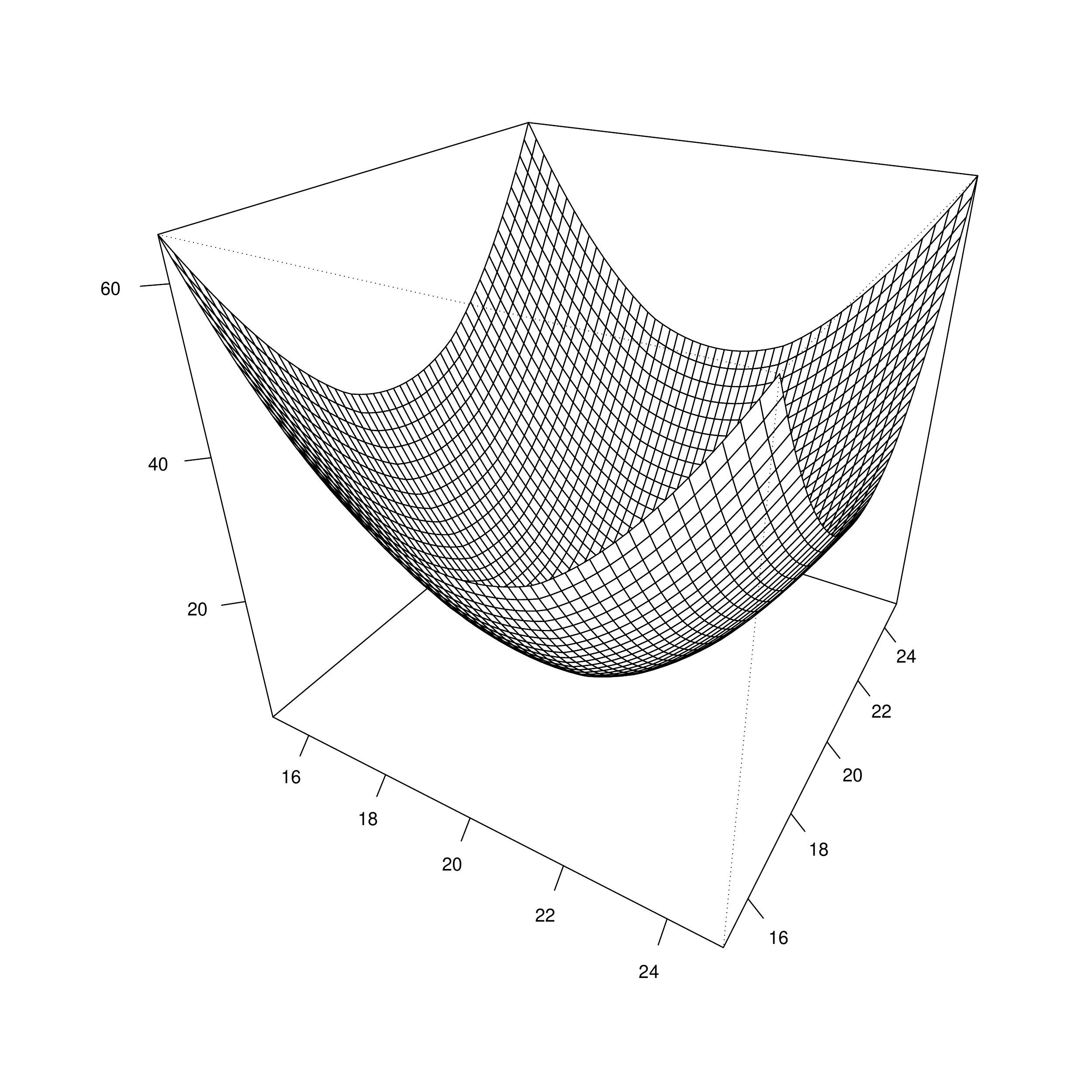}	
  	\caption{Sphere Worst} \label{fig:SphereWorst}
	\end{subfigure}		
	\begin{subfigure}[t]{.24\textwidth}
		\includegraphics[width=\textwidth]{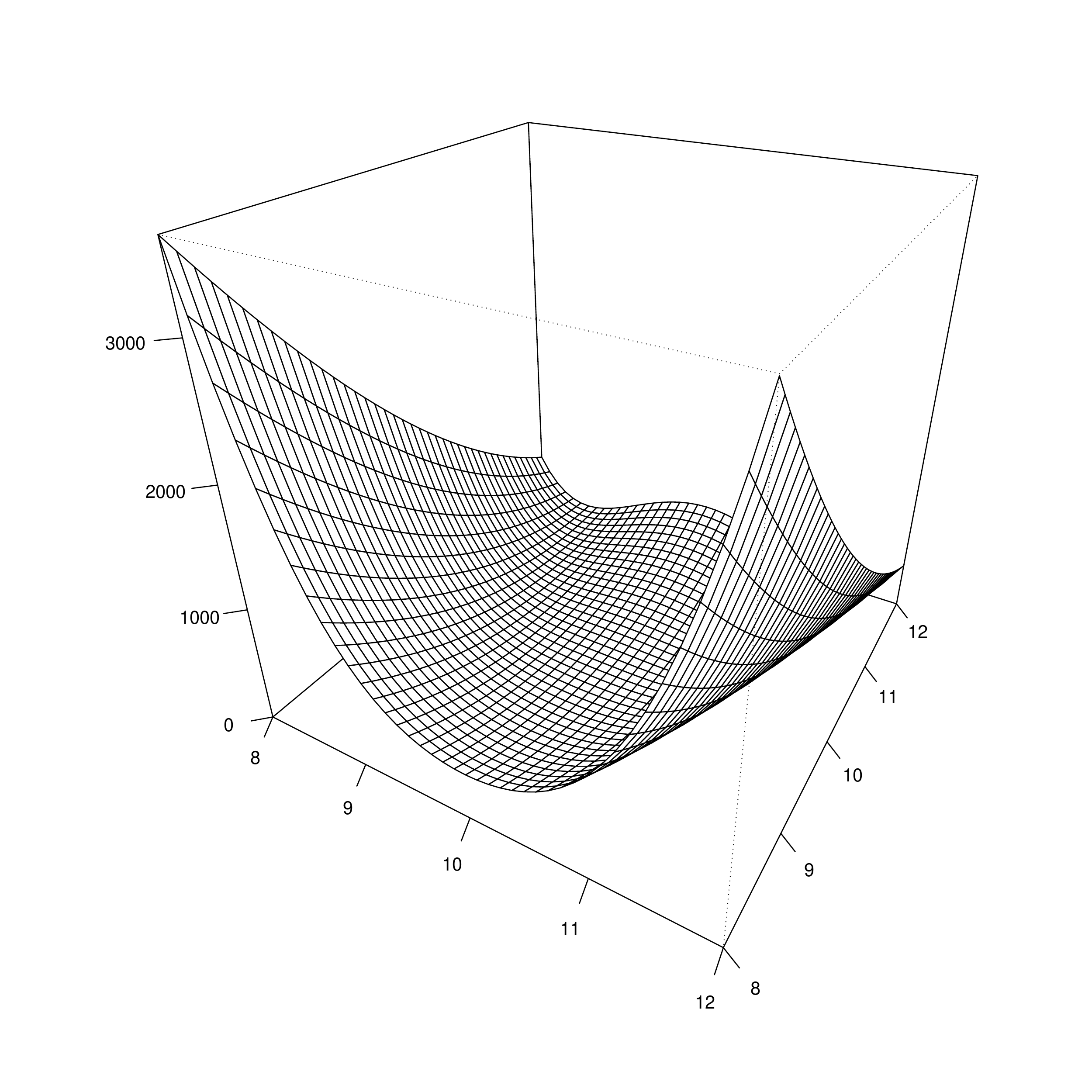}
	  \caption{Rosenbrock Nominal} \label{fig:RosenbrockNom}
	\end{subfigure}%
	\begin{subfigure}[t]{.24\textwidth}
		\includegraphics[width=\textwidth]{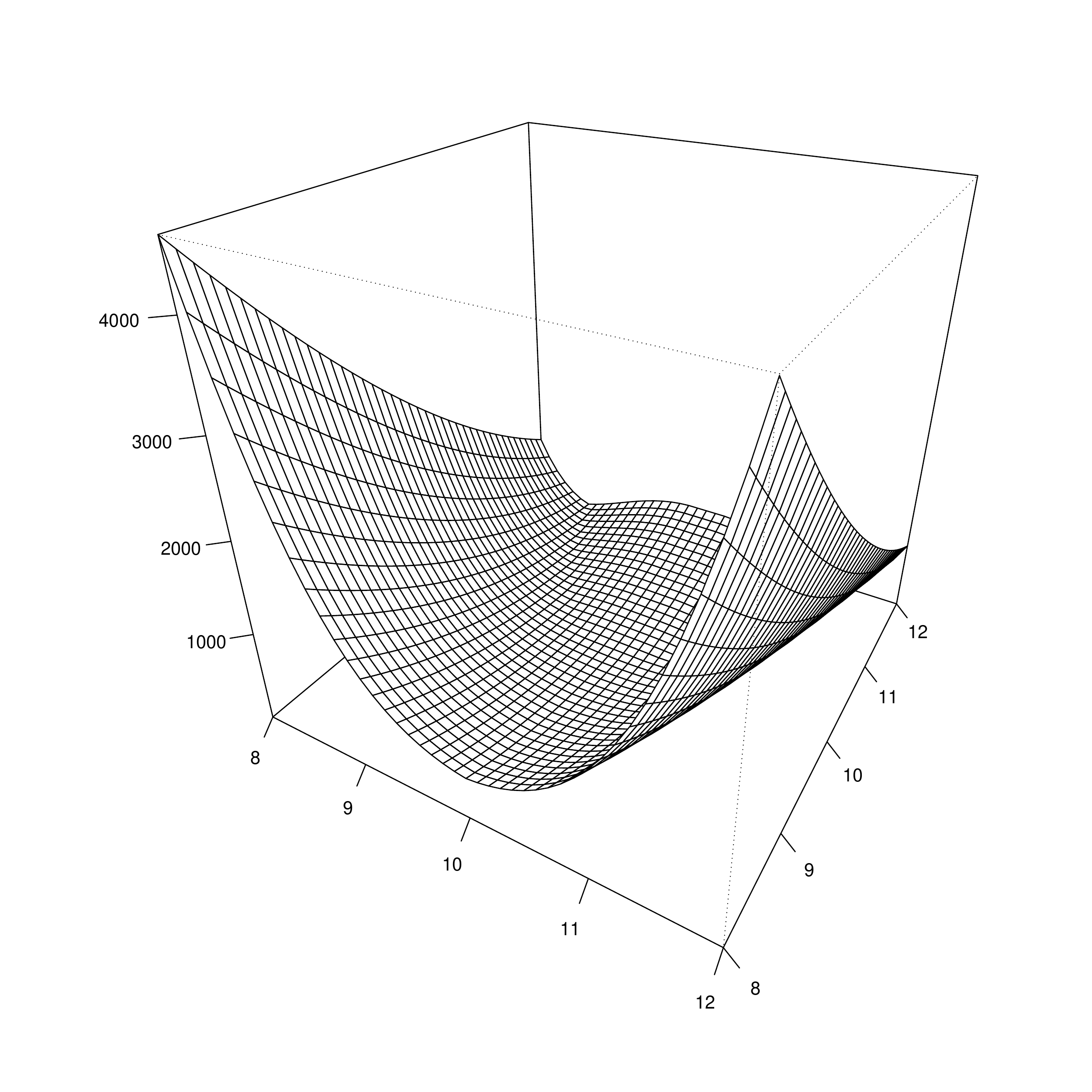}		
  	\caption{Rosenbrock Worst} \label{fig:RosenbrockWorst}
	\end{subfigure}	
	\caption{Plots of 2D versions of the 10 multi-dimensional problems in the rPSO test suite.}
	\label{fig:rPSOtestSuite}
\end{figure}

Table~\ref{tab:funcs} presents the eleven test functions used within our experimental testing. The functions are based on the literature:  \cite{Branke1998, KruisselbrinkEmmerichBack2010, KruisselbrinkReehuisDeutzBackEmmerich2011, Kruisselbrink2012, JamilYang2013,BertsimasNohadaniTeo2010}, and their mathematical description is given in Appendix~\ref{sec:testFunctionFormulae}. 3D plots of the 2D versions of these multi-dimensional functions are shown in Figures~\ref{fig:poly2D} and~\ref{fig:rPSOtestSuite}.

\begin{table}[h]
\begin{center}
\begin{tabular}{r|ll}
Name & $\X$ & $\Gamma$ \\
\hline
Rastrigin & $[14.88, 25.12]^n$ & $0.5$\\
MultipeakF1 & $[-5, -4]^n$ & $0.0625$\\
MultipeakF2 & $[10, 20]^n$ & $0.5$\\
Branke's Multipeak & $[-7, -3]^n$ & $0.5$\\
Pickelhaube & $[-40, -20]^n$ & $1$\\
Heaviside Sphere & $[-30, -10]^n$ & $1$\\
Sawtooth & $[-6, -4]^n$ & $0.2$\\
Ackley & $[17.232, 82.768]^n$ & $3$\\
Sphere & $[15, 25]^n$ & $1$\\
Rosenbrock & $[7.952, 12.048]^n$ & $0.25$\\
2D polynomial & $[-1, 4]^2$ & $0.5$
\end{tabular}
\caption{Test functions.}\label{tab:funcs}
\end{center}
\end{table}

The common features of the objective function surfaces of individual test problems means that such problems can be associated with high level categorisations such as multi-modality, basins or valleys, see e.g. \cite{JamilYang2013}. Here we have arranged the ten multi-dimensional test problems in an approximate order based on modality and common features, which should be apparent in the 3D plots shown in Figure~\ref{fig:rPSOtestSuite}.

All of our algorithms have been coded in Java, with calls to the mathematical programming software IBM ILOG CPLEX Optimization Studio V12.6.3 to solve the second order cone problem element of a d.d.\ calculation.

\subsection{Parameter tuning}
\label{sec:parameterTuning}

Parameter tuning has been undertaken at the dimensional level, separately for each comparator metaheuristic, and for each of the six dimensions. This seems reasonable in a practical setting, where a decision maker is likely to have a good advance understanding of the dimension of the problem at hand. This generates a set of parameters for each heuristic separately for each dimension in the testing process.

Our tuning uses an evolutionary tuning approach applied to a subset of the test suite (four instances), including only problems where the nominal and robust global optima are differently located. This is primarily to ensure that the tuned level of inner maximisation is not biased by test problems where the nominal and robust optima are coincident, and therefore where no $\Gamma$-uncertainty neighbourhood analysis might be desirable.

Within the tuning GA each member of the evolving population represents a set of parameter values for a given heuristic operating at a specific dimension. For all comparator heuristics the same level of tuning is employed for a given dimension, i.e., each problem in the tuning subset of the test suite is run for the same number of samples, and the same population size and number of generations is used. For a given member of the evolutionary population the mean of the sample runs on a single tuning test problem is calculated, and ranked in comparison to the other individuals in the population. The average ranking of these means over the tuning test suite is used as the measure of utility within an evolutionary tournament selection, see e.g. \cite{Eiben2012}. 

For all heuristics a key tuned parameter is the extent of the inner maximisation analysis. For d.d.\ the parameters that control the parameter $\sigma$ which impacts the determination of high cost points, and the minimum step size required for any descent direction step, are tuned. This also applies to the enhanced rPSO heuristics employing a d.d.\ calculation. For LEH the parameters that control the genetic algorithm employed in the identification of the largest empty hypersphere devoid of all high cost points, see equation (lehMM), are tuned. Again this applies to the enhanced rPSO heuristics employing an LEH calculation.  For all four rPSO based heuristics (baseline rPSO and three enhanced heuristics developed here) the new parameters that are tuned are stated in Sections~\ref{sec:bruteForcePSO},~\ref{sec:rPSOdd},~\ref{sec:rPSOleh} and~\ref{sec:rPSOlehdd}.

\subsection{Results}
\label{sec:results}

\begin{figure}[htbp]
\begin{center}
	\includegraphics[width=0.6\textwidth, height=0.16\textheight]{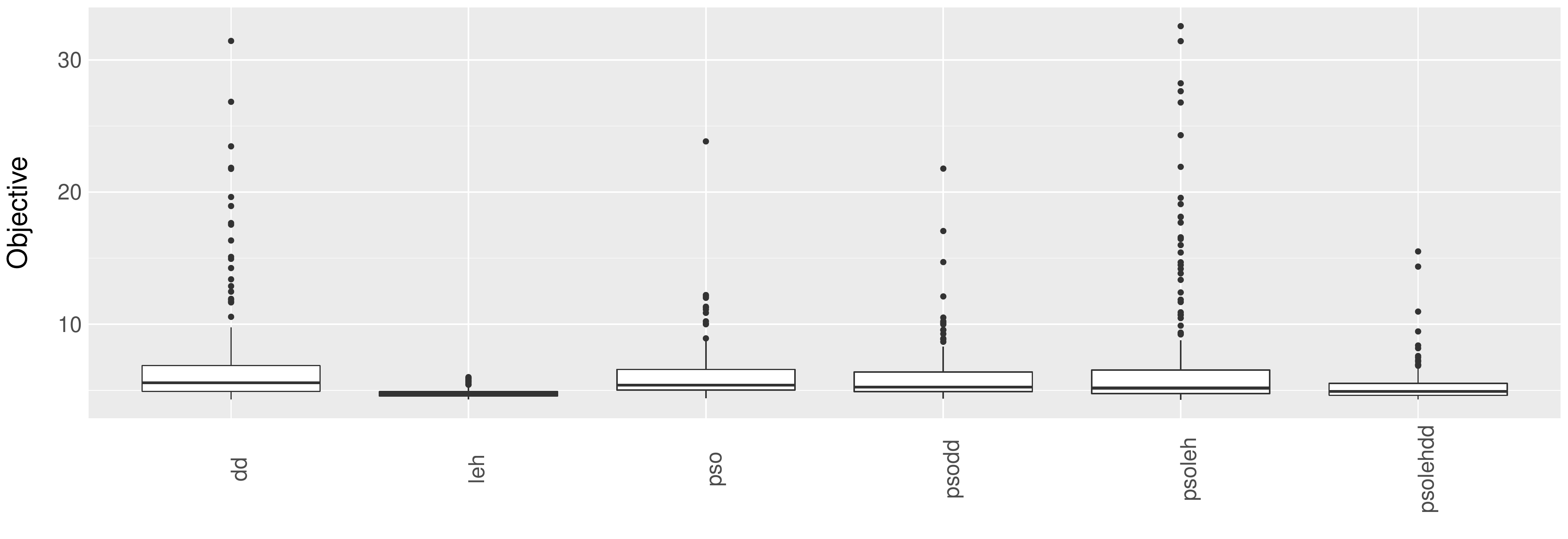}	
	\caption{Box plots of 2D Polynomial test function \cite{BertsimasNohadaniTeo2010} robust objective values for 200 sample runs.}
	\label{fig:poly2Dbox}
\end{center}
\end{figure}

Here we present results of the 200 sample runs for each of the six comparator robust metaheuristics when applied to each of the 61 test problem-dimension instances. The mean results are shown in Tables~\ref{fig:BertMeanResults} and~\ref{fig:meanResults}. The best or statistically equivalent to the best results (best-equivalent) due to Wilcoxon rank-sum testing with 95\% confidence and employing a Bonferonni correction (see e.g. \cite{HastieTibshiraniFriedman2009}) are highlighted. That is a highlight on a method means that no other method is statistically better. For Table~\ref{fig:meanResults} this applies at the cell level. Note that being best with respect to the mean objective value does not always correspond to being statistically best-equivalent. Statistical analysis was undertaken in R, see \cite{RCore2019, Dinno2017}. Details of the distributions of the 200 samples runs for the multi-dimensional problems are provided by box plots in Appendix~\ref{sec:boxPlots}. The box plots for the 2D Bertsimas polynomial are shown in Figure~\ref{fig:poly2Dbox}.

\begin{table}[htbp]

\vspace{-8mm} 

\begin{center}
\begin{tabular}{r|l}
Heuristic & Mean \\
\hline
d.d.  & 6.91 \\
LEH & \textbf{4.80} \\
rPSO (\hspace{-0.8mm}~\ref{BaselineRPSO}) & 6.10 \\
rPSOdd (\hspace{-0.8mm}~\ref{BaselineRPSO},\hspace{-0.8mm}~\ref{ddVelComponent})  & 5.97 \\
rPSOleh (\hspace{-0.8mm}~\ref{BaselineRPSO},\hspace{-0.8mm}~\ref{partLEHAlgorithm})  & 7.13 \\
rPSOlehdd (\hspace{-0.8mm}~\ref{BaselineRPSO},\hspace{-0.8mm}~\ref{ddVelComponent},\hspace{-0.8mm}~\ref{partLEHAlgorithm})  & 5.29
\end{tabular}
\caption{Mean results for 200  sample runs for the 2D polynomial function due to \cite{BertsimasNohadaniTeo2010}. Statistically equivalent best heuristics are highlighted. Bracketed numbering on rPSO based heuristics refers to the outer minimisation algorithms used.}
\label{fig:BertMeanResults}
\end{center}
\end{table}

\begin{table}[htbp]

\vspace{-5mm} 

\begin{center}
	\includegraphics[width=1.0\textwidth]{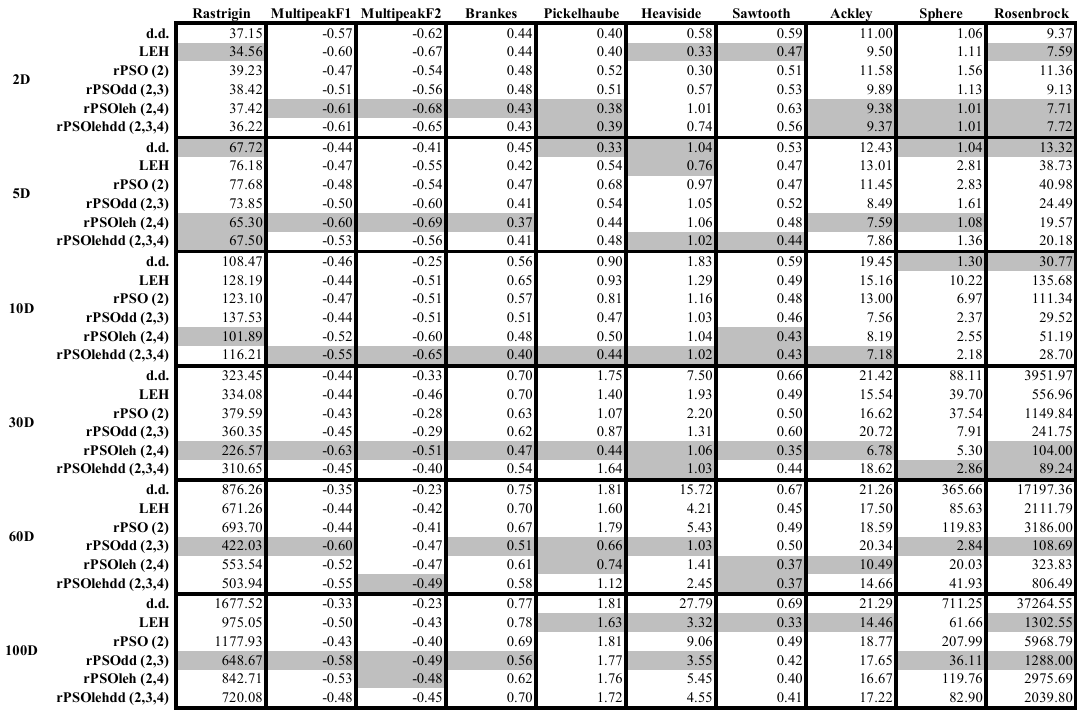} 		
	\caption{Mean results for 200 sample runs for the 10 multi-dimensional problems. Statistically equivalent best heuristics are highlighted. Bracketed numbering on rPSO based heuristics refers to the outer minimisation algorithms used.}
	\label{fig:meanResults}
\end{center}
\end{table}

We see that the 2D instances are dominated by the rPSOleh, baseline LEH and rPSOlehdd heuristics, which are best-equivalent in 64\%, 45\% and 36\% of instances respectively; here, e.g. 64\% refers to rPSOleh being best or statistically equivalent to the best in seven out of the eleven 2D test problem instances. For 5D, rPSOleh is best-equivalent in 60\% of cases, followed by the baseline d.d.\ with 50\%, rPSOlehdd 30\%, and the baseline LEH 10\%. However for 10D the picture changes, with rPSOlehdd dominating as best-equivalent in 70\% of cases, followed by rPSOleh and the baseline d.d.\ with 20\% each.

The 30D and 60D instances are completely dominated by the new heuristics. For 30D rPSOleh is best-equivalent in 90\% of cases, with only rPSOlehdd achieving any other best-equivalent results, in 30\% of cases. At 60D rPSOdd comes into prominence, being best-equivalent in 70\% of cases, with 30\% for rPSOleh and 20\% for rPSOlehdd. Finally for 100D rPSOdd is again best-equivalent in 70\% of cases, although now the baseline LEH heuristic is best-equivalent in 50\% of cases, with 10\% for rPSOleh.

Table~\ref{fig:summaryResults} summarises the proportion of the 61 test problem instances for which each robust metaheuristic is identified as the best or statistically equivalent to the best. The three new heuristics lead the order of best to worst results: rPSOleh, rPSOlehdd, rPSOdd, LEH, and d.d.\, with rPSO failing to be the best for any test instance.

\begin{table}[htb]
\begin{center}
\begin{tabular}{r|l}
Heuristic & Best-equiv. \\
\hline
d.d. & 11.48\% \\
LEH & 18.03\% \\
rPSO (\hspace{-0.8mm}~\ref{BaselineRPSO}) & 00.00\% \\ 
\hline
rPSOdd (\hspace{-0.8mm}~\ref{BaselineRPSO},\hspace{-0.8mm}~\ref{ddVelComponent}) & 22.95\% \\
rPSOleh (\hspace{-0.8mm}~\ref{BaselineRPSO},\hspace{-0.8mm}~\ref{partLEHAlgorithm}) & 44.26\% \\
rPSOlehdd (\hspace{-0.8mm}~\ref{BaselineRPSO},\hspace{-0.8mm}~\ref{ddVelComponent},\hspace{-0.8mm}~\ref{partLEHAlgorithm}) & 31.15\%
\end{tabular}
\caption{Summary of the proportion of best or statistically equivalent to the best results for each heuristic. Bracketed numbering on rPSO based heuristics refers to the outer minimisation algorithms used.}
\label{fig:summaryResults}
\end{center}
\end{table}

From the summary Table~\ref{fig:summaryResults} it is not clear how each of the new rPSO heuristics perform individually against the three baseline comparators. Therefore we have conducted a number of further statistical tests. First we compared rPSOdd against the three baselines, rPSOleh against the baselines, and rPSOlehdd against the baselines. The result is that rPSOdd is best-equivalent in 49.2\% of the 61 test instances, rPSOleh is best-equivalent in 75.4\% of instances, and rPSOlehdd is best-equivalent in 65.6\% of instances. In each test the baseline LEH heuristic had the next best performance, with best-equivalent results of 39.4\%, 21.3\% and 29.5\% of the test instances respectively.

Finally we compared each of the new rPSO heuristics individually against each of the three baseline comparators, in a series of one-to-one analyses. The results are summarised in Table~\ref{fig:OneToOneResults}, where each cell comprises three values: top, middle and bottom. Each cell shows the proportion of the 61 test problem instances for which each new heuristic is better than (top), statistically equivalent to (middle), or worse than (bottom) each baseline heuristic. Within each set of three values the highest is highlighted.

In the comparison between the rPSO and rPSOdd methods, specifically considering the impact of augmenting the baseline rPSO velocity equations (Vel1) with a d.d.\ component (Vel2), rPSO is best or statistically equivalent to the best in 23.0\% of the 61 test instances, with rPSOdd best-equivalent in 88.5\%. The dominance of the new rPSO heuristics in one-to-one comparisons against each of the baselines is clear, in particular for rPSOleh and rPSOlehdd. Indeed rPSOleh is individually better than each of the baselines in nearly 80\% of test instances.


\begin{table}[htbp]
\begin{center}
\begin{tabular}{r|lll}
 & d.d. & LEH & rPSO (\hspace{-0.8mm}~\ref{BaselineRPSO}) \\
\hline
 \multirow{3}{*}{rPSOdd (\hspace{-0.8mm}~\ref{BaselineRPSO},\hspace{-0.8mm}~\ref{ddVelComponent})} & \textbf{65.6\%} & \textbf{52.5\%} & \textbf{77.0\%}  \\
 & 6.6\% & 8.2\% & 11.5\%  \\
 & 27.9\% & 39.3\% & 11.5\% \\
 \hline
\multirow{3}{*}{rPSOleh (\hspace{-0.8mm}~\ref{BaselineRPSO},\hspace{-0.8mm}~\ref{partLEHAlgorithm})} & \textbf{86.9\%} & \textbf{78.7\%} & \textbf{93.4\%} \\
 & 3.3\% & 3.3\% & 3.3\%  \\
 & 9.8\% & 18.0\% & 3.3\% \\
 \hline
\multirow{3}{*}{rPSOlehdd (\hspace{-0.8mm}~\ref{BaselineRPSO},\hspace{-0.8mm}~\ref{ddVelComponent},\hspace{-0.8mm}~\ref{partLEHAlgorithm})} & \textbf{86.9\%} & \textbf{70.5\%} & \textbf{91.8\%} \\
 & 3.3\% & 3.3\% & 1.6\%  \\
 & 9.8\% & 26.2\% & 6.6\% \\
\end{tabular}
\caption{Results of one-to-one statistical tests between new and baseline heuristics. Each cell shows the percentage of test problem instances where: (top) the new heuristic is  best, (middle) the new and baseline heuristics are equivalent, and (bottom) the baseline heuristic is best. Bracketed numbering on rPSO based heuristics refers to the outer minimisation algorithms used.}
\label{fig:OneToOneResults}
\end{center}
\end{table}


\section{Conclusions and further work}
\label{sec:concusionsFurtherWork}

We have developed a framework encompassing three new robust metaheuristics for box-constrained, black-box robust optimisation problems under implementation uncertainty. Our robust approaches follow a min max setting, seeking to find solutions that optimise the worst performance. Our new approaches use a baseline robust particle swarm population-based heuristic as a frame, adapting elements of two existing individual-based robust metaheuristics, descent directions \cite{BertsimasNohadaniTeo2010} and largest empty hypersphere \cite{HughesGoerigkWright2019}, and combing them along with new features. The following novel features are introduced here:

\begin{itemize}
	\item An extension of the PSO movement formulation to include particle level, iteration level, d.d.\ vector information: exploiting uncertainty neighbourhood information in order to provide a locally optimal directional movement component.
	\item Efficiency savings in terms of numbers of function evaluations due to the use of a particle level stopping condition.
	\item The introduction of the concept of dormancy, whereby the repeated non-requirement to perform any function evaluations is monitored at a particle level in order to interrupt particles trapped in previously visited regions or outside of the feasible region.
	\item The relocating of dormant particles by an optimal exploration-focussed calculation of the largest empty hypersphere devoid of all high cost points.
\end{itemize}

This results in a framework encompassing three new heuristics:

\begin{itemize}
	\item rPSOdd: augments the baseline rPSO heuristic with an additional component in the standard PSO velocity equation, using the descent direction vector that optimally points away from the worst $\Gamma$-uncertainty neighbourhood points around a candidate point (particle location).
	\item rPSOleh: augments rPSO with both the stopping condition and the determination of the largest hypersphere empty of previously evaluated poor points, from the LEH heuristic. The stopping condition allows efficiencies in the inner maximisation calculations, which are terminated early if any $\Gamma$-neighbourhood point is identified with nominal function value worse than a particle's current personal best information. The calculation of an LEH is used to relocate particles that have become `dormant', either due to repeated movements outside of the feasible region or repeated movements in areas of the feasible region where points with high nominal function value have already been identified.
	\item rPSOlehdd: is simply a combination of both of the rPSOdd and rPSOleh heuristics, incorporating all new features of these approaches.
\end{itemize}

The performance of the new heuristics has been assessed by applying them to 61 test problem instances, covering six dimensions up to 100D, a single 2D problem and 10 multi-dimensional problems. The performance of our new approaches has been compared against three existing baseline approaches, a repeating dd. approach, LEH, and a baseline robust PSO. Our new approaches are shown to outperform the existing approaches across all dimensions. For 10D, 30D and 60D instances the new approaches dominate. They also outperform the baseline heuristic for other dimensions, although both LEH and d.d.\ also produce some good results. 

One potential extension of our new framework is to undertake explicit inner maximisation searches, for example using PSO or GA searches, as opposed to the use of uniform random sampling here. Also, in order to make these techniques more widely applicable other forms of uncertainty, for example model uncertainty, could be accommodated. The current focus on a $\Gamma$-radius uncertainty neighbourhood can also be extended, to the consideration of other descriptions of a point's uncertainty neighbourhood.

\appendix

\section*{Appendices}

\section{Test functions}
\label{sec:testFunctionFormulae}

Functions used in the experimental testing of the enhanced rPSO metaheuristics. These functions are based on \cite{Branke1998, KruisselbrinkEmmerichBack2010, KruisselbrinkReehuisDeutzBackEmmerich2011, Kruisselbrink2012, JamilYang2013, BertsimasNohadaniTeo2010}.

\begingroup
\allowdisplaybreaks
\begin{align*}
\text{Rastrigin:}\quad &  f(\pmb{x}) = 10 n + \sum_{i=1}^{n}{[{(x_i-20)}^{2} - 10 \cos(2\pi (x_i-20))]} \\
& \X = [14.88, 25.12]^n  \\
\text{MultipeakF1:}\quad & f(\pmb{x}) = - \frac{1}{n} \sum_{i=1}^n g(x_i) \\
& g(x_i) =
\begin{cases}
	e^{-2 \ln 2 ( \frac{(x_i+5)-0.1}{0.8} )^2} \sqrt{ \left| \sin(5\pi (x_i+5)) \right| }  &\quad\text{if } 0.4 < x_i+5 \le 0.6 \text{ ,} \\
  e^{-2 \ln 2 ( \frac{(x_i+5)-0.1}{0.8} )^2} \sin^{6}(5\pi (x_i+5))  &\quad\text{otherwise} \\
\end{cases}\\
& \X = [-5, -4]^n \\
\text{MultipeakF2:}\quad & f(\pmb{x}) = \frac{1}{n} \sum_{i=1}^n g(x_i) \\
& g(x_i) = 2\sin(10 \exp (-0.2(x_i-10)) (x_i-10)) \exp (-0.25(x_i-10))\\
& \X = [10, 20]^n \\
\text{Branke's Multipeak:}\quad &f(\pmb{x}) = \max \lbrace c_1, c_2 \rbrace - \frac{1}{n} \sum_{i-1}^n g(x_i) \\
& g(x_i) =
\begin{cases}
	c_1 \Bigg(1-{\frac{4((x_i+5)+\frac{b_1}{2})^2}{b_1^2}} \Bigg)  &\quad\text{if } -b_1 \le (x_i+5) < 0 \text{ ,}
	\\
	c_2 \cdot 16^{\frac{-2 \left| b_2 - 2(x_i+5) \right|}{b_2}} &\quad\text{if } 0 \le (x_i+5) \le b_2 \text{ ,}
	\\
  0  &\quad\text{otherwise} \\
\end{cases} \\
& b_1=2, b_2=2, c_1=1, c_2=1.3 \\ 
& \X = [-7, -3]^n\\
\text{Pickelhaube:}\quad & f(\pmb{x}) = \frac{5}{5-\sqrt{5}} - \max \lbrace g_0(\pmb{x}), g_{1a}(\pmb{x}), g_{1b}(\pmb{x}), g_2(\pmb{x}), \rbrace \\
& g_0(\pmb{x}) = \frac{1}{10} e^{-\frac{1}{2} \| \pmb{x}+30 \|} \\
&g_{1a}(\pmb{x}) = \frac{5}{5-\sqrt{5}} \Bigg( 1-\sqrt{\frac{\| \pmb{x}+30+5 \|}{5\sqrt{n}}} \Bigg) \\
&g_{1b}(\pmb{x}) = c_1 \Bigg( 1-\Bigg(\frac{\| \pmb{x}+30+5 \|}{5\sqrt{n}} \Bigg)^4 \Bigg) \\
&g_2(\pmb{x}) =  c_2 \Bigg( 1-\Bigg(\frac{\| \pmb{x}+30-5 \|}{5\sqrt{n}} \Bigg)^{d_2} \Bigg) \\
& c_1=625/624, c_2=1.5975, d_2=1.1513 \\
& \X = [-40, -20]^n \\
\text{Heaviside Sphere:}\quad & f(\pmb{x}) = \Bigg( 1 - \prod_{i=1}^n g(x_i) \Bigg) +  \sum_{i=1}^n \Bigg(\frac{(x_i+20)}{10} \Bigg)^2 \\
& g(x_i) =
\begin{cases}
	0  &\quad\text{if } 0 < (x_i+20) \text{ ,} \\
  1  &\quad\text{otherwise} \\
\end{cases} \\
& \X =  [-30, -10]^n\\
\text{Sawtooth:}\quad & f(\pmb{x}) = 1 - \frac{1}{n} \sum_{i=1}^n g(x_i) \\
& g(x_i) =
\begin{cases}
	(x_i+5) + 0.8  &\quad\text{if } -0.8 \le (x_i+5) < 0.2 \text{ ,} \\
  0  &\quad\text{otherwise} \\
\end{cases} \\
& \X =  [-6, -4]^n \\
\text{Ackleys:}\quad & f(\pmb{x}) = -20 \exp \Bigg( -0.2 \sqrt{\frac{1}{n} \sum_{i=1}^n (x_i-50)^2} \Bigg) \\ 
& \phantom{f(\pmb{x}) = }- \exp \Bigg( \frac{1}{n} \sum_{i=1}^n \cos(2\pi (x_i-50)) \Bigg) + 20 + \exp(1) \\
& \X = [17.232, 82.768]^n \\
\text{Sphere:}\quad & f(\pmb{x}) = \sum_{i=1}^n (x_i-20)^2 \\
& \X = [15, 25]^n \\
\text{Rosenbrock:}\quad & f(\pmb{x}) = \sum_{i=1}^{n-1} [100((x_{i+1}-10) - (x_i-10)^2)^2 + ((x_i-10)-1)^2]\\
& \X = [7.952, 12.048] \\
\text{2D polynomial:}\quad & f(x, y) = 2x^6 - 12.2x^5 + 21.2x^4 + 6.2x - 6.4x^3 - 4.7x^2 - y^6 - 11y^5 \\
&+  43.3y^4 - 10y - 74.8y^3 + 56.9y^2 - 4.1xy - 0.1y^2x^2 + 0.4y^2x + 0.4x^2y \\
&\X = [-1, 4]
\end{align*}

\endgroup

\section{Box plots}
\label{sec:boxPlots}

Box plots of the results of the experimental testing on our six comparator robust metaheuristcs applied to 60 test problems covering the 10 multi-dimensional problems. Each plot is based on 200 sample runs of each heuristic applied to each problem instance.
\\

\begin{figure}[htbp]
	\centering

			\begin{subfigure}{\textwidth}
				\centering
        \includegraphics[width=1.0\textwidth]{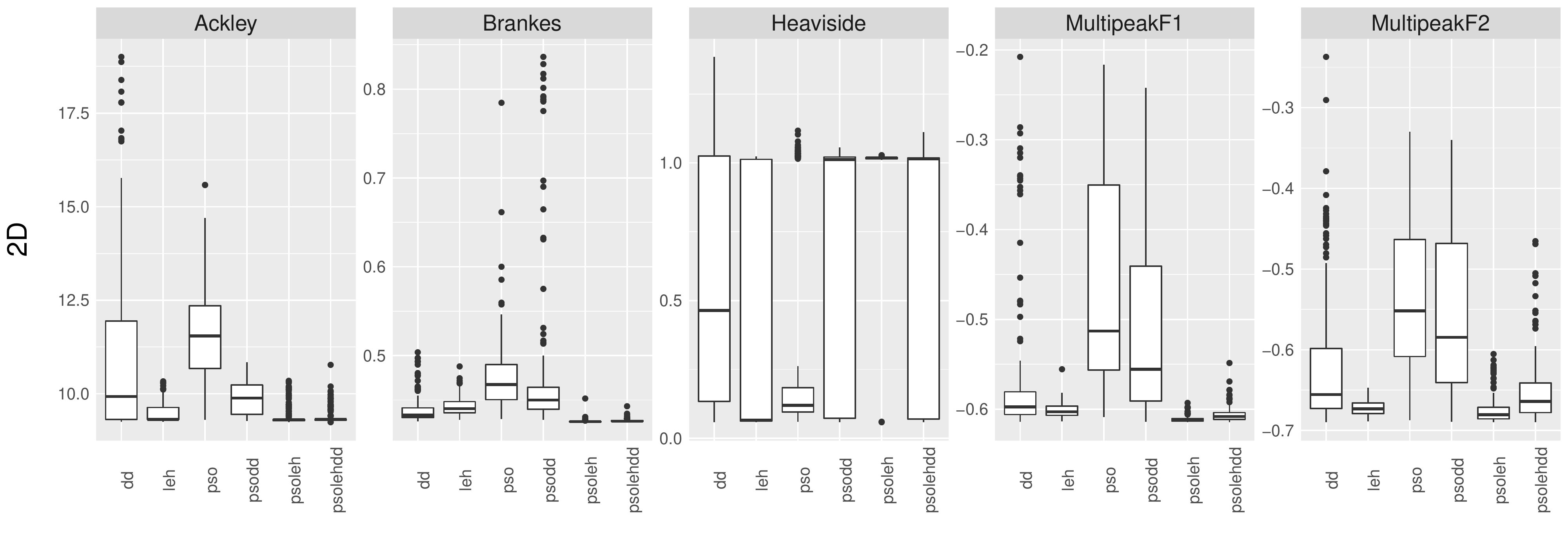}
				\label{fig:boxplots2Df}
				\vspace*{-1.0mm} 
			\end{subfigure}	
			
			\vspace*{-4.0mm} 

			\begin{subfigure}{\textwidth}
				\centering
        \includegraphics[width=1.0\textwidth]{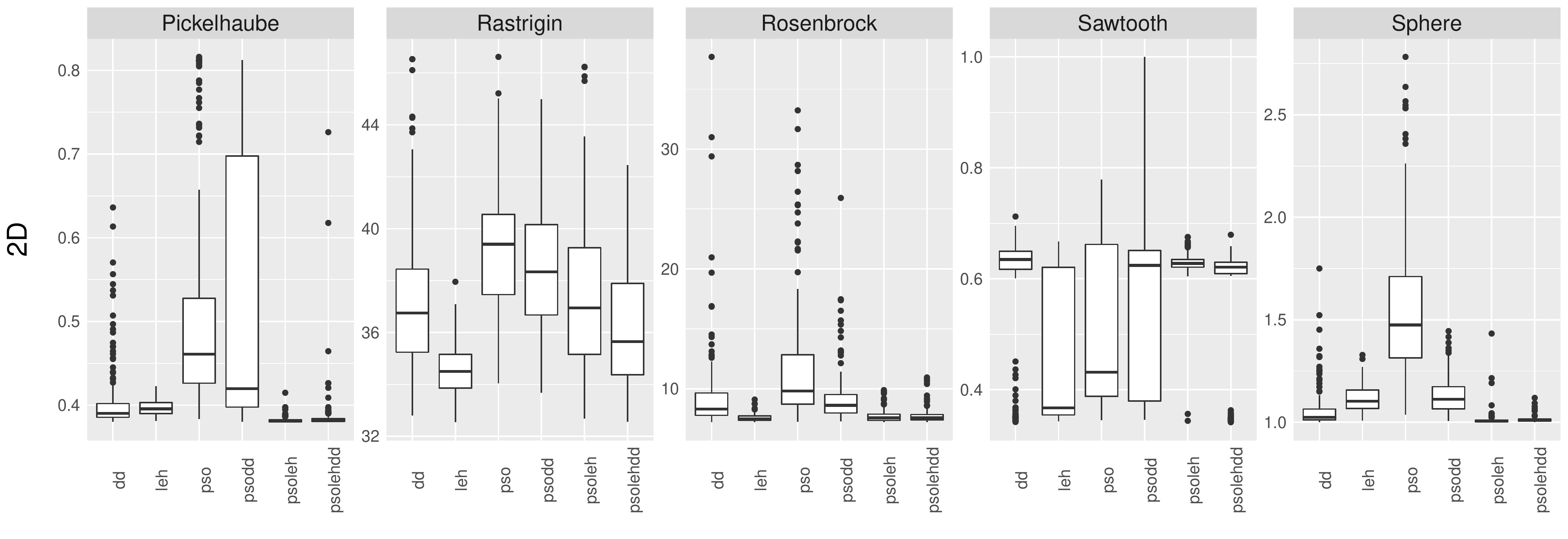}
				\label{fig:boxplots2Db}
				\vspace*{-1.0mm} 
			\end{subfigure}	
			
			\vspace*{-4.0mm} 

			\begin{subfigure}{\textwidth}
				\centering
        \includegraphics[width=1.0\textwidth]{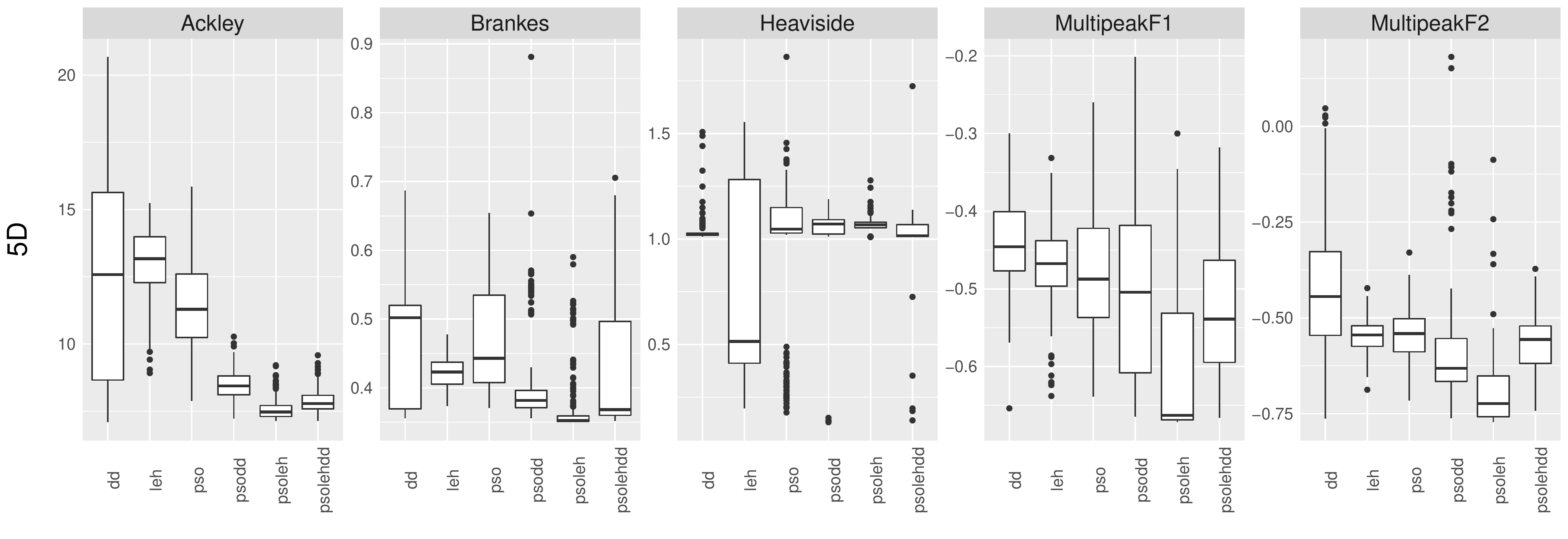}
				\label{fig:boxplots5Df}
				\vspace*{-1.0mm} 
			\end{subfigure}	
			
			\vspace*{-4.0mm} 

			\begin{subfigure}{\textwidth}
				\centering
        \includegraphics[width=1.0\textwidth]{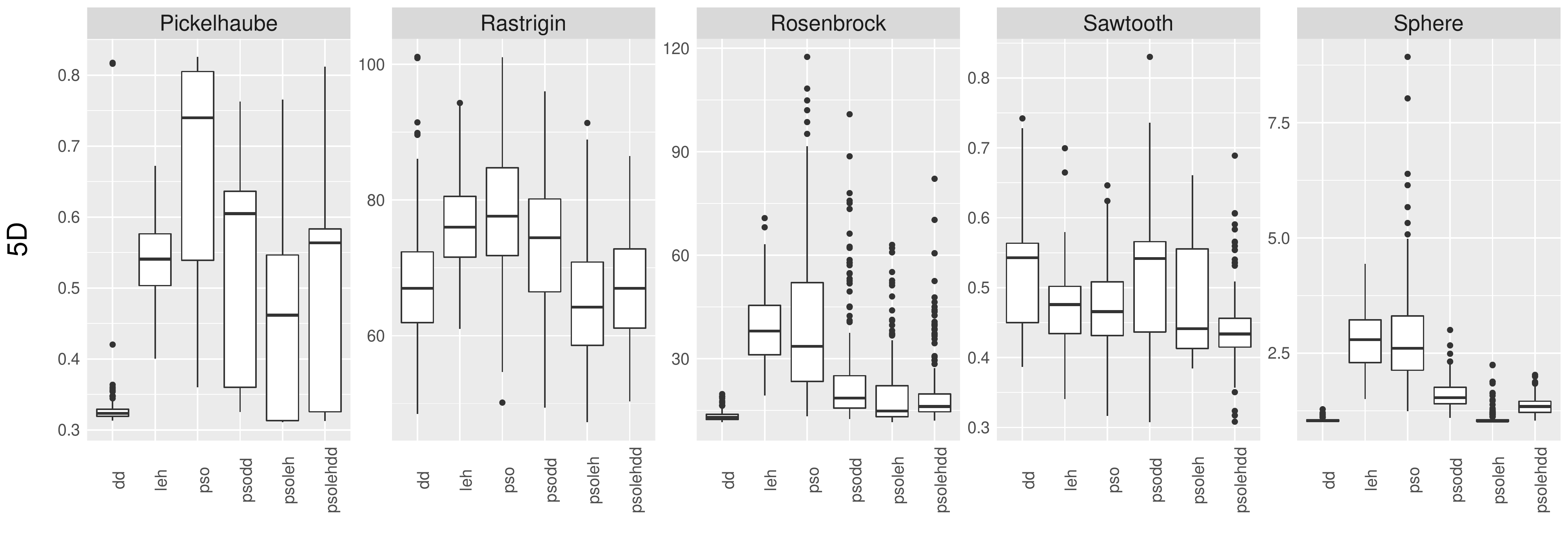}
				\label{fig:boxplots5Db}
				\vspace*{-1.0mm} 
			\end{subfigure}	

	\caption{Box plots of 2D and 5D robust objective values for 200 sample runs.}
	\label{fig:BoxPlotResultsA}
\end{figure}

\begin{figure}[htbp]
	\centering

			\begin{subfigure}{\textwidth}
				\centering
        \includegraphics[width=1.0\textwidth]{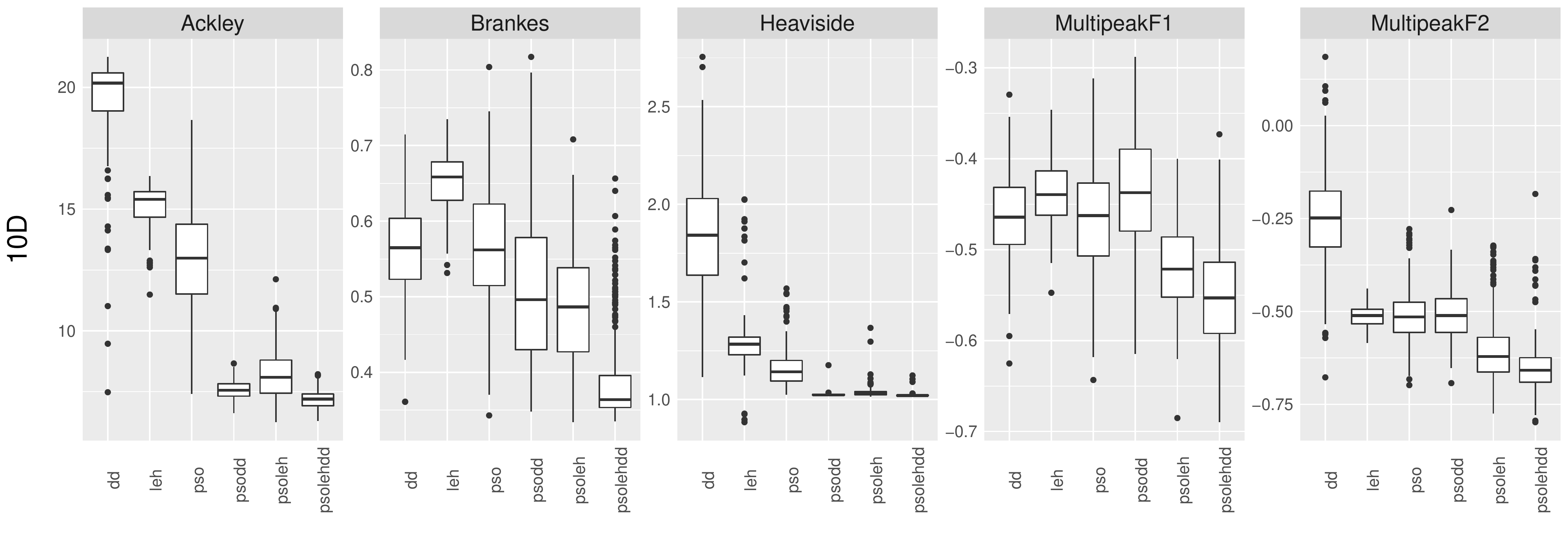}
				\label{fig:boxplots10Df}
				\vspace*{-1.0mm} 
			\end{subfigure}	
			
			\vspace*{-4.0mm} 

			\begin{subfigure}{\textwidth}
				\centering
        \includegraphics[width=1.0\textwidth]{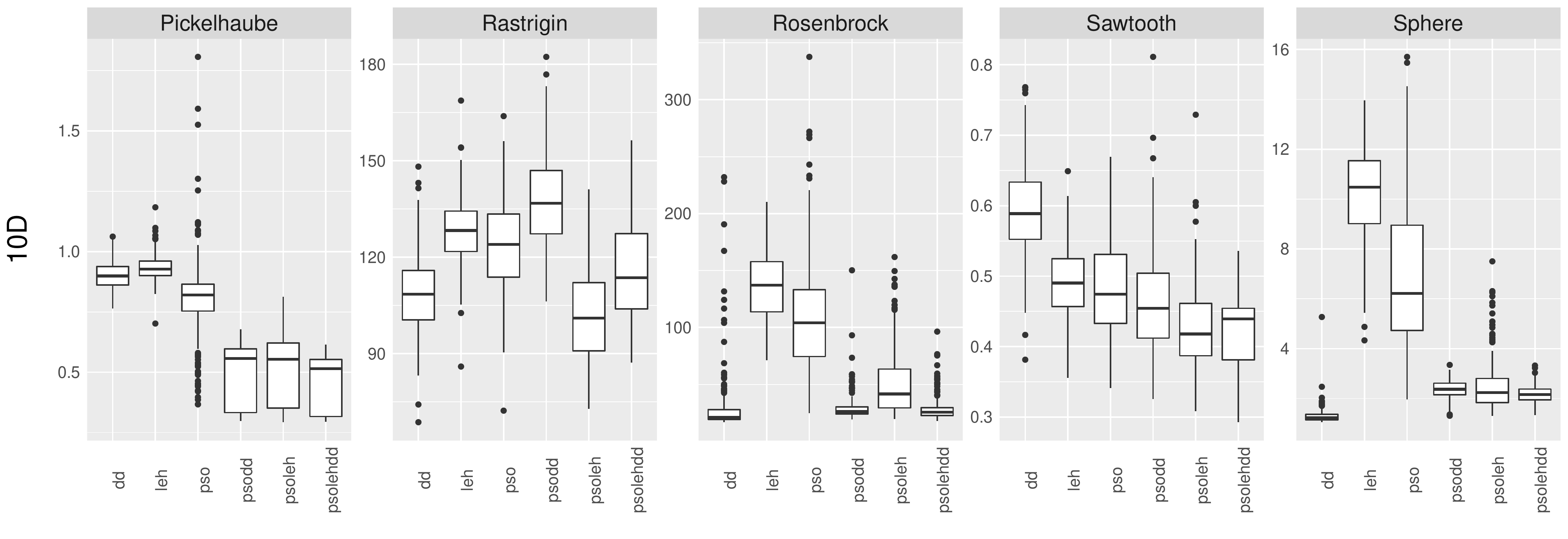}
				\label{fig:boxplots10Db}
				\vspace*{-1.0mm} 
			\end{subfigure}		
			
			\vspace*{-4.0mm} 

			\begin{subfigure}{\textwidth}
				\centering
        \includegraphics[width=1.0\textwidth]{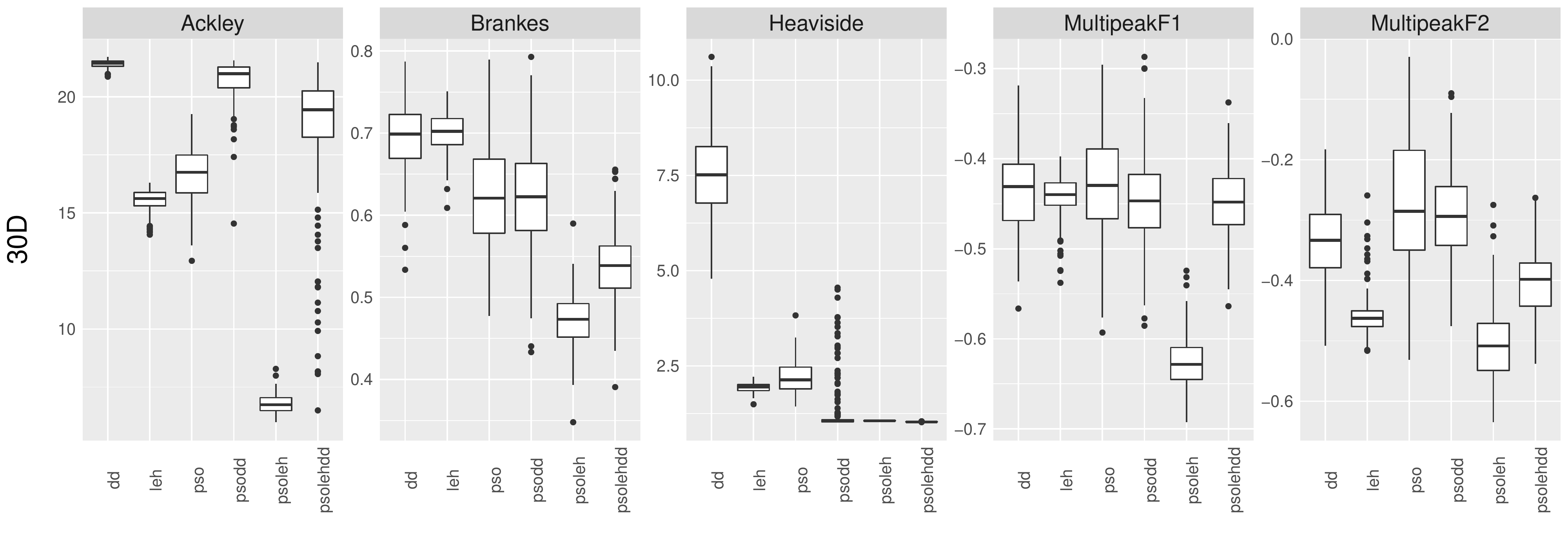}
				\label{fig:boxplots30Df}
				\vspace*{-1.0mm} 
			\end{subfigure}	
			
			\vspace*{-4.0mm} 

			\begin{subfigure}{\textwidth}
				\centering
        \includegraphics[width=1.0\textwidth]{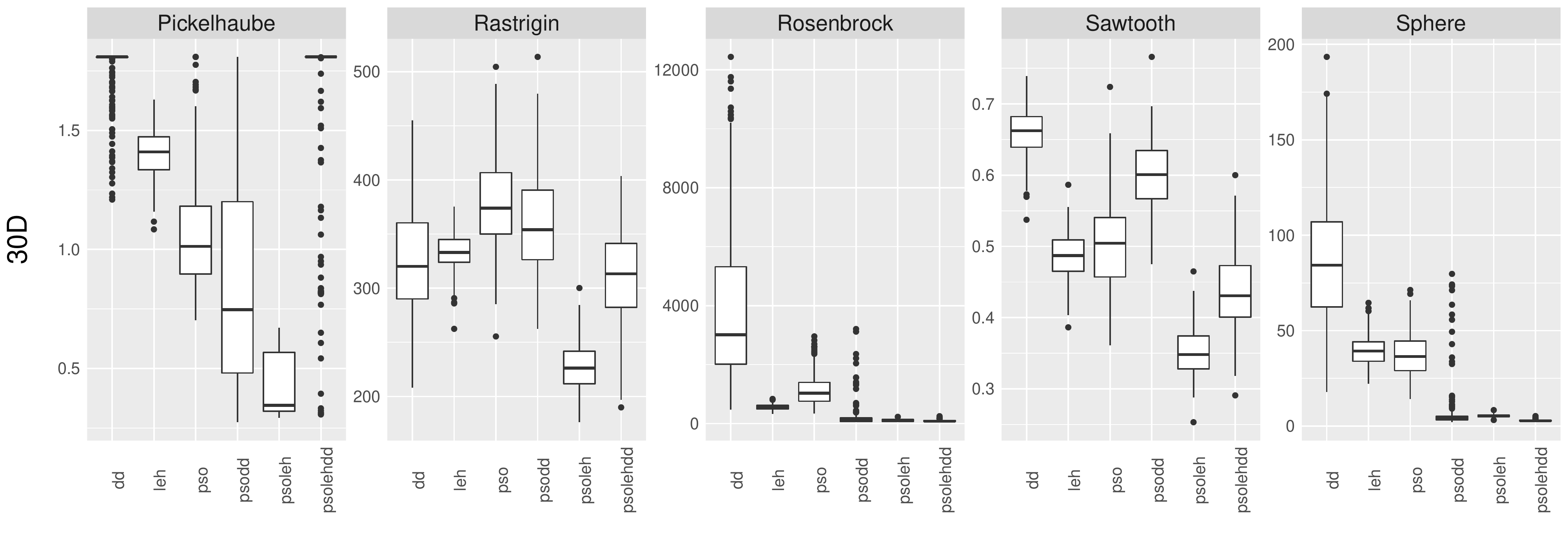}
				\label{fig:boxplots30Db}
				\vspace*{-1.0mm} 
			\end{subfigure}	

	\caption{Box plots of 10D and 30D robust objective values for 200 sample runs.}
	\label{fig:BoxPlotResultsB}
\end{figure}

\begin{figure}[htbp]
	\centering

			\begin{subfigure}{\textwidth}
				\centering
        \includegraphics[width=1.0\textwidth]{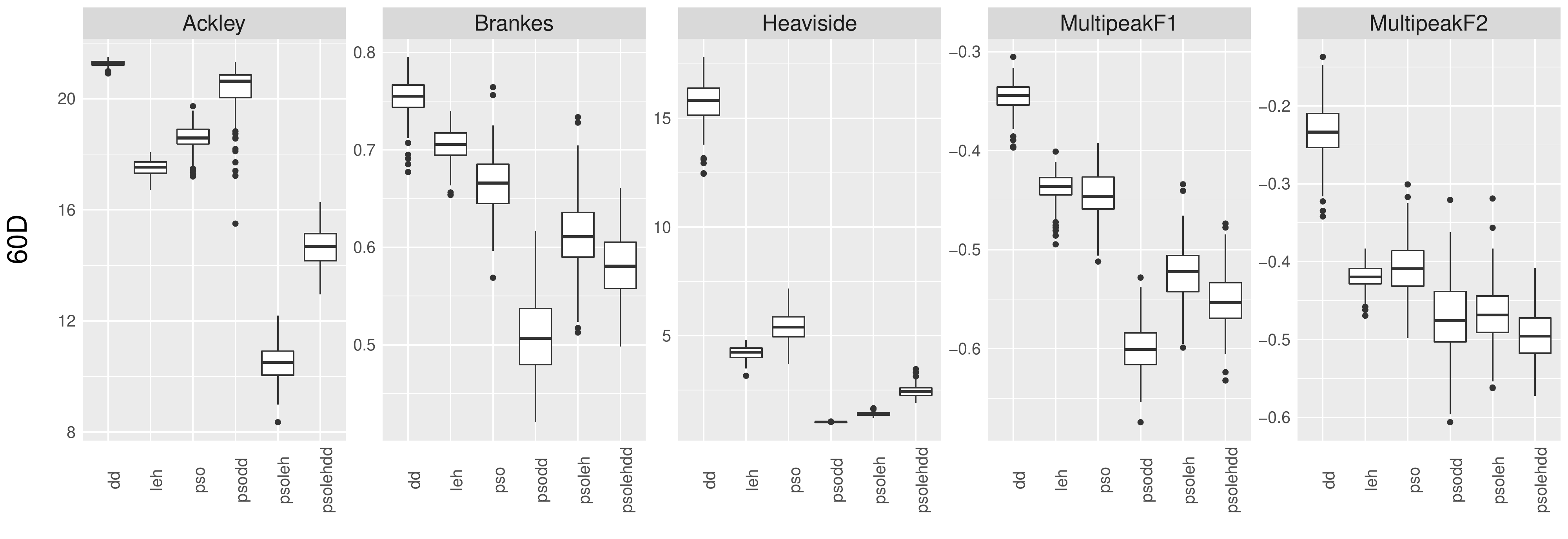}
				\label{fig:boxplots60Df}
				\vspace*{-1.0mm} 
			\end{subfigure}	
			
			\vspace*{-4.0mm} 

			\begin{subfigure}{\textwidth}
				\centering
        \includegraphics[width=1.0\textwidth]{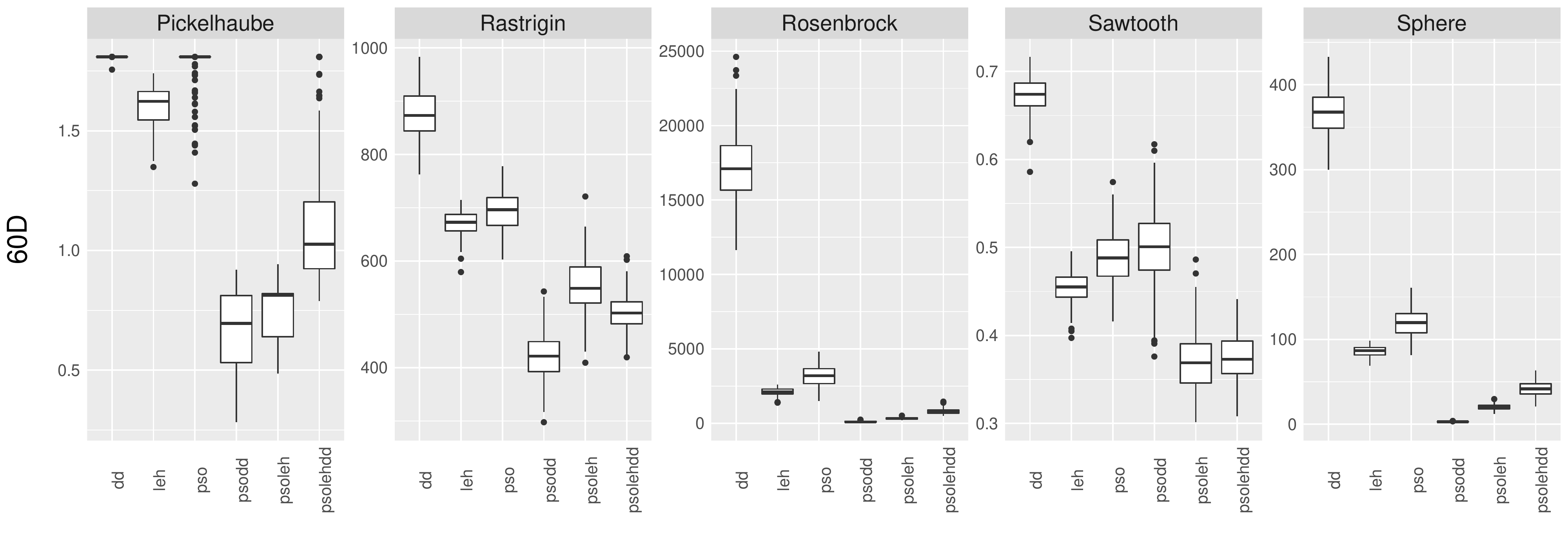}
				\label{fig:boxplots60Db}
				\vspace*{-1.0mm} 
			\end{subfigure}		
			
			\vspace*{-4.0mm} 

			\begin{subfigure}{\textwidth}
				\centering
        \includegraphics[width=1.0\textwidth]{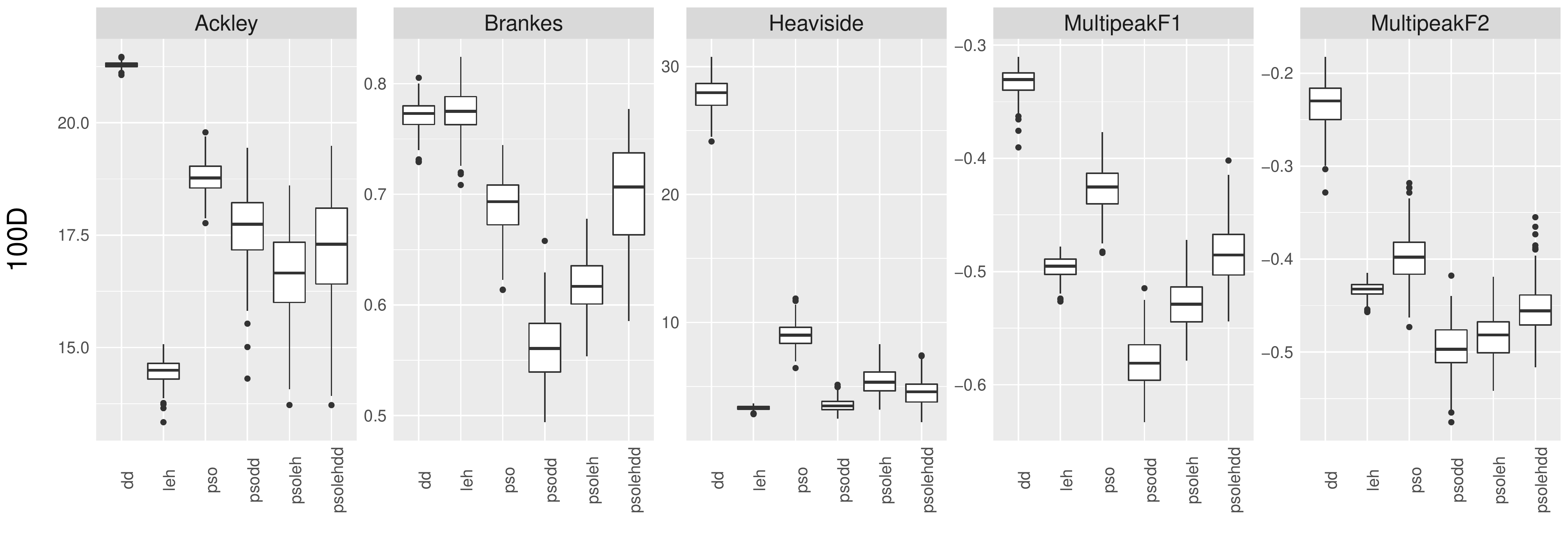}
				\label{fig:boxplots100Df}
				\vspace*{-1.0mm} 
			\end{subfigure}	
			
			\vspace*{-4.0mm} 

			\begin{subfigure}{\textwidth}
				\centering
        \includegraphics[width=1.0\textwidth]{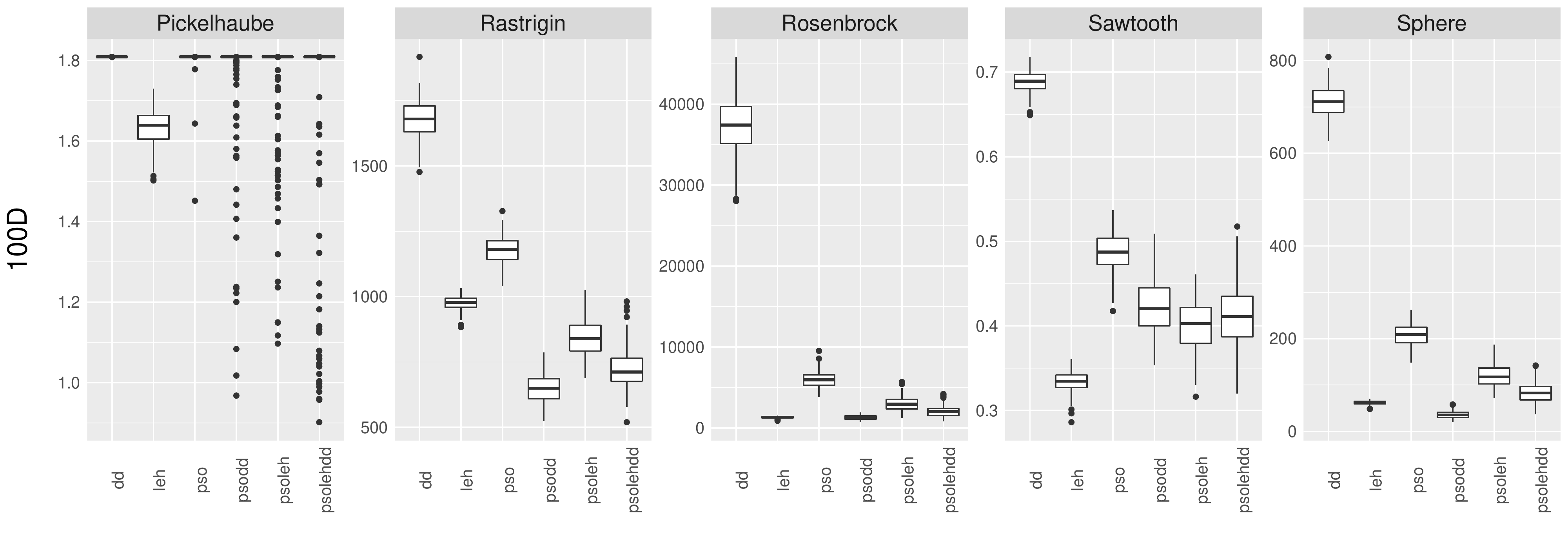}
				\label{fig:boxplots100Db}
				\vspace*{-1.0mm} 
			\end{subfigure}	

	\caption{Box plots of 60D and 100D robust objective values for 200 sample runs.}
	\label{fig:BoxPlotResultsC}
\end{figure}

\newpage

\section{List of Abbreviations}
\label{sec:abbreviations}

\begin{table}[htb]
\begin{center}
\begin{tabular}{p{0.15\textwidth}|p{0.8\textwidth}} 
Abbreviation & Definition \\
\hline
d.d. & descent directions \\
GA & genetic algorithm \\
hcp & high cost point \\
LEH & largest empty hypersphere \\
PSO & particle swarm optimisation \\
rPSO & robust particle swarm optimisation \\
rPSOdd & robust particle swarm optimisation with descent directions \\
rPSOleh & robust particle swarm optimisation with largest empty hypersphere  \\
rPSOlehdd & robust particle swarm optimisation with descent directions and largest empty hypersphere 
\end{tabular}
\caption{Commonly used abbreviations.}
\label{fig:abbrevList}
\end{center}
\end{table}

\clearpage 
\newpage



\newcommand{\etalchar}[1]{$^{#1}$}

\end{document}